\documentclass[12pt]{article}
\usepackage[english]{babel}
\usepackage[latin1]{inputenc}
\usepackage{amsfonts,amssymb,amsmath, epsfig}
\usepackage{color,graphicx,graphics,psfrag}
\usepackage{amsmath,amstext,amssymb,amsfonts, amscd}
\usepackage{subfigure}

\usepackage{hyperref}           

\def\Box{\leavevmode\vbox{\hrule
     \hbox{\vrule\kern4pt\vbox{\kern4pt}%
           \vrule}\hrule}}
\def\blackbox{\leavevmode\vrule height 5pt width 4pt depth 0pt\relax}
\def\endproof{\null\hfill {$\blackbox$}\bigskip}

\newcounter{appendix}
\def\appendix{\advance\c@appendix by 1
   \def\thesection{\Alph{section}}
   \ifnum\c@appendix=1 \setcounter{section}{-1} \fi
   \@startsection {section}{1}{\z@}{-3.5ex plus -1ex minus 
   -.2ex}{2.3ex plus .2ex}{\Large\bf}}

\def\paragraph#1{{\bf #1\ }}

\newtheorem{lemma}{Lemma}[section]  

\newtheorem{theorem}[lemma]{Theorem}

\newtheorem{proposition}[lemma]{Proposition}

\newtheorem{remark}{Remark}[section]

\newtheorem{hypothesis}{Hypothesis}[section]

\title{Kinetic hierarchy and propagation of chaos\\in biological swarm models} 
\author{E. Carlen$^{1}$, R. Chatelin$^{2,3}$, P. Degond$^{2,3}$, B. Wennberg$^{4}$} 
\date{} 
\begin{document}

\maketitle

\vspace{0.5 cm}

\begin{center}
1-
Department of Mathematics,  
Rutgers University \\
110 Frelinghuysen Rd., Piscataway NJ 08854-8019 \\
email: carlen@math.rutgers.edu
\end{center}

\begin{center}
2-Université de Toulouse; UPS, INSA, UT1, UTM ;\\ 
Institut de Mathématiques de Toulouse ; 
F-31062 Toulouse, France. \\
email: robin.chatelin@math.univ-toulouse.fr\\
email: pierre.degond@math.univ-toulouse.fr
\end{center}

\begin{center}
3-CNRS; Institut de Mathématiques de Toulouse UMR 5219 ;\\ 
F-31062 Toulouse, France.
\end{center}

\begin{center}
4-Department of Mathematical Sciences, \\
Chalmers University of Technology, \\
SE41296 Gotheburg \\
email: wennberg@chalmers.se 
\end{center}

\vspace{0.5 cm}
\begin{abstract}
We consider two models of biological swarm behavior. In these models,
pairs of particles interact to adjust their velocities one to each
other. In the first process, called 'BDG', they join their average
velocity up to some noise. In the second process, called 'CL', one of
the two particles tries to join the other one's velocity. This paper
establishes the master equations and BBGKY hierarchies of these two
processes. It investigates the infinite particle limit of the
hierarchies at large time-scale. It shows that the resulting kinetic
hierarchy for the CL process does not satisfy propagation
  of chaos. Numerical simulations indicate that the BDG process has
similar behavior to the CL process.  
\end{abstract}

\medskip
\noindent
{\bf Acknowledgements:} The first author acknowledges support from the 'R\'egion Midi-Pyr\'en\'ees' government in the frame of the 'Chaires Pierre-de-Fermat'. The second author acknowledges support from the ANR under contract 'CBDif-Fr' (ANR-08-BLAN-0333-01). The third author acknowledges support from the Swedish Research Council.


\vskip 0.4cm

\setcounter{equation}{0}
\section{Introduction}
\label{sec_intro}

The derivation of kinetic equations from particle models of swarming
behavior has recently received a great deal of attention. In
biological swarm modeling, the most widely used models are particle
ones (also known as 'Individual-Based Models') \cite{Aoki, Chuang,
  Couzin, Cucker_Smale, Vicsek}. However, to investigate the large
scale behavior of biological systems such as fish schools or insect
swarms, kinetic \cite{BDM, CFRT, DM2, hillen_diffusion_2000,
  painter_modelling_2009} and hydrodynamic \cite{Chuang, Mogilner1,
  TB2} models have proved to be valuable alternatives. The question of
showing a rigorous link between the particle and kinetic levels is
mostly open. In \cite{BCC2}, a mean-field limit of the Vicsek particle
model \cite{Vicsek} is performed and leads to a nonlinear
Fokker-Planck equation proposed in \cite{DM}. A similar program has
been performed for the Cucker-Smale model \cite{BCC1, CFRT}. In
\cite{BDG}, Bertin, Droz and Gr\'egoire propose a binary collision
mechanism  which mimics the Vicsek alignment interaction \cite{Vicsek}
and formally derive a Boltzmann-like kinetic collision operator. There
has been no rigorous justification of this derivation so far. The
present paper is a step in this direction.   

In this work, we investigate two examples of particle systems
representative of swarming behavior, the so-called BDG and CL
processes. These two processes mimic the formation of consensus in
biological groups about the direction of motion to follow. They are
binary processes. In the first process, called 'BDG' (after Bertin,
Droz and Gr\'egoire \cite{BDG}), two interacting particles join their
average velocity up to some noise. In the second process, called 'CL'
(for 'Choose the Leader'), one of the two particles tries to join the
other one's velocity up to some noise. In this paper, we focus on
space-homogeneous problems and ignore the spatial
variables. Consequently, interactions may happen among any pair of
individuals in the pool with a certain probability. We also assume
that the individuals move in a two-dimensional space with unit
speed. The state of each particle is described by its velocity vector
$v$ on the one-dimensional sphere ${\mathbb S}^1$.  

The state of an $N$-particle system can be described by its
$N$-particle probability $F_N$. In the present framework, $F_N$ is a
function of the $N$ velocity coordinates $(v_1, \ldots, v_N)$ on the
torus ${\mathbb T}^N = ({\mathbb S}^1)^N$ and of time. The particle
dynamics translates  into a time-evolution equation for $F_N$ called
the 'master equation'. In a previous work \cite{CDW}, we have
investigated the class of 'pair-interaction driven' master equations,
of which the BDG and CL master equations are members. We have shown
that, as $N \to \infty$, propagation of chaos
  holds. A propagation of chaos result states that the solution
  $F_N(t)$ can be approximated (in a sense to be defined below) by an
  $N$-fold tensor product of the single-particle distribution $F_1(t)$
  provided that this property is true initially. This means that the
  particles become nearly independent and that the system can be
  described by its single-particle distribution $F_1(t)$ instead of
  the $N$-particle distribution $F_N$. The dimension of the problem is
  therefore considerably reduced.  

To investigate the large $N$ limit, it is difficult to work with $F_N$ alone. Indeed, the limit of $F_N$ as $N \to \infty$ is literally a function of an infinite number of variables. The functional treatment is simplified by considering the $k$-particle marginal $F_{N,k}$, which is the joint probability of any subset of $k$ particles. The number of variables involved in $F_{N,k}$ is $k$ and stays fixed as $N \to \infty$. The drawback of this method is that the equation satisfied by $F_{N,k}$ depends on the other marginals $F_{N,k'}$ in general. Thus, the equations for the $(F_{N,k})_{k \in \{1, \ldots, N\} }$ are all coupled together, forming the so-called BBGKY hierarchy \cite{CIP}. When $N = \infty$, the hierarchy involves an infinite number of coupled equations and is called the kinetic hierarchy. Showing a propagation result in the limit $N \to \infty$ involves breaking the coupling between the equations in the kinetic hierarchy in some way. 

Consensus formation in swarm models should be associated with the
build-up of correlations between the particles over time. The fact
that the BDG and CL models, as a result of \cite{CDW},  satisfy a
propagation of chaos result is counter-intuitive. The resolution of
this paradox lies in the investigation of time scales. Indeed, the
result of \cite{CDW} is only valid on finite time intervals at the
kinetic scale. On this time scale, the number of collisions undergone
by each particle is bounded independently of $N$. The present paper
investigates whether correlation build-up happens at larger time
scales.  

Large time scales are investigated thanks to an appropriate time
rescaling, i.e. a change in the time time. However, the dynamics must
also be rescaled in some way to keep the leading order terms in the
hierarchy finite. Here, the appropriate rescaling consists in letting
the variance of the noise involved in the process tend to zero. In the
BDG dynamics, this is not enough. For this reason, we introduce a
'Biased BDG' dynamics, where the collision probability depends on the
relative velocities of the particles. Then, the rescaling also
involves a grazing collision limit, i.e. having the collision occur
only if the relative velocities of the two particles are small.  

The main objective of this paper is two establish the kinetic
hierarchies for the rescaled BDG and CL processes. We then investigate
whether these hierarchies possess solutions which satisfy propagation
of chaos. For the CL hierarchy, we show that it is never the case. In
\cite{CDW}, it was already established that the invariant densities
(i.e. the stationary solutions) do not satisfy propagation of
chaos. This was done by looking at the single and two-particle
marginals only. Here, we extend \cite{CDW} by showing that the
time-dependent solution of the CL hierarchy never satisfies
propagation of chaos either. We also provide a general formula for the
$k$-particle marginal invariant density.  

Concerning the BDG dynamics, the situation is unclear, in spite of the
apparent simplicity of the hierarchy equations. We notice that uniform
densities are stationary solutions of the BDG hierarchy. However, the
question of uniqueness of stationary solutions for this hierarchy is
open. There might exist other solutions which do not satisfy the chaos
property. In a companion paper \cite{CDW2}, we investigate the kinetic
equation associated to the BDG dynamics. From \cite{CDW}, we know that
propagation of chaos is true and that the kinetic equation is valid on
the kinetic time scale. The uniform distribution is clearly a
stationary solution of this kinetic equation. However, in \cite{CDW2},
we show that this equilibrium is linearly unstable if the noise level
is small enough. This suggests the existence of a second class of
anisotropic equilibria (similar to the Von-Mises equilibria of
\cite{DFL, FL}).  
The existence of multiple equilibria may be a clue that  propagation
of chaos is broken at infinite time. However, these considerations are
pure conjectures at this level. 

To improve our understanding, we use numerical experiments. We
generate the stationary one and two particle marginals by running a
large number of independent time-dependent runs of the particle
dynamics. The experimental results concerning the CL dynamics
consolidate the theoretical findings. In particular, the theoretical
and numerical $2$-particle correlations show remarkably good
agreement. The experimental study of the BDG dynamics shows a similar
behavior to the CL dynamics. For this reason, it should be expected
that the BDG dynamics lacks chaos property on the large time
scale. However, a rigorous result in this direction is not available
yet.    

In the literature,  propagation of chaos has been mainly investigated
in the context of the Boltzmann equation and its caricature proposed
by Kac. Early works involve the names of Kac, Lanford, McKean and
others \cite{Kac, Lanford3, McKean}. They have initiated a
considerable activity \cite{GM3, IP1, Sznitman}. A new approach has
been recently developed in \cite{MMW}. Rates of convergence towards
chaos have been investigated mainly in the context of the Kac model
recently in \cite{CCL2, CCLRV, DSC, Janv, Maslen}.  

Lachowicz~\cite{Lachowicz} has recently considered a class of biologically
motivated Markov jump processes, and proves a propagation of chaos
result as part of the derivation of kinetic and macroscopic
equations. Some of the  basic ideas going into his proof, as well as
into ours, are present in  the original works of Kac~\cite{Kac}, but the
models, and hence the implementation of these ideas, are quite
different.

The outline of the paper is as follows. In section
\ref{sec_particle_models}, we present the two particle processes that
we are interested in. In section \ref{sec_master_eq}, we derive the
master equations and BBGKY hierarchies of these processes. Section
\ref{sec_grazing} is the core of the paper. It performs the
$N\to\infty$ limit in the rescaled hierarchies and develop the
consequences that result from it. Section \ref{sec_numerics} reports
on the numerical experiments. A conclusion is drawn in section
\ref{sec_conclusion}. Finally, two appendices collect the technical
proofs of the main results of the paper.

\setcounter{equation}{0}
\section{Particle models on the circle ${\mathbb S}^1$}
\label{sec_particle_models}

We consider systems of $N$ particles on the circle ${\mathbb S}^1$. Each configuration of the system corresponds to a $N$-tuple $\omega = (v_1, \ldots , v_N) \in {\mathbb T}^N$ with ${\mathbb T}^N$ being the $N$-dimensional torus ${\mathbb T}^N = ({\mathbb T}^1)^N$. This state space can be seen as corresponding to the velocities of a system of mutually interacting swarming agents (see \cite{CDW} for details). 

The dynamics is first defined as a time-discrete dynamics. Let $\omega^n$ denote the value at the $n$-th iterate. We now consider different rules for passing from $\omega^n$ to $\omega^{n+1}$.

\subsection{The BDG dynamics}
\label{subsec_def_BDG}

This dynamics is named after Bertin, Droz and Grégoire \cite{BDG} who introduced it as a model for the Vicsek dynamics \cite{Vicsek}. Given the state $\omega^n \in {\mathbb T}^N$ at time $n$, finding $\omega^{n+1}$ consists of the following steps:

\begin{itemize}
\item[-]  pick an un-ordered pair $(i,j)$ ($i<j$) randomly according to a uniform distribution, i.e. with probability 
\begin{equation} P_{i,j} = \frac{2}{N(N-1)} . \label{eq_uniform_unordered} \end{equation}
and compute an 'average velocity' 
\begin{equation} \hat v_{ij}^n = \frac{v_i^n + v_j^n}{|v_i^n + v_j^n|} . \label{hat_vij} \end{equation}

\item[-] Then define
\begin{equation}
v_i^{n+1} = w_i^n \hat v_{ij}^n,  \quad v_j^{n+1} = w_j^n \hat v_{ij}^n , 
\label{eq_rule_BDG}
\end{equation}
where $w_i^n$ and $w_j^n$ are two independent identically distributed random variables on ${\mathbb S}^1$ distributed according to the probability $g(w)$. The notations use the multiplicative group structure of ${\mathbb S}^1$. We suppose that $g$ is symmetric: 
$$ g(w^*) = g(w), \quad \forall w \in {\mathbb S}^1. $$

\end{itemize}

\subsection{The biased BDG dynamics}
\label{subsec_def_BBDG}

In the sequel, we will consider a 'biased' version of the BDG dynamics defined as follows: 
Let $H(w)$ be a function $w \in {\mathbb S}^1_+ \to H(w) \in [0,1]$ be given,  where ${\mathbb S}^1_+ = \{ z \in {\mathbb S}^1 \, | \, \mbox{Re} \, z \geq 0 \}$. We also assume that $H$ is symmetric: $H(w) = H(w^*)$. The biased BDG dynamics is similar to the BDG dynamics except for an acceptance-rejection procedure based on $H$. Namely, the procedure consists of the following steps:

\begin{itemize}
\item[-]  pick an un-ordered pair $(i,j)$ ($i<j$) randomly according to a uniform distribution, i.e. with probability (\ref{eq_uniform_unordered}) and compute $\hat v_{ij}^n$ according to (\ref{hat_vij}).

\item[-] With probability $H(\hat v_{ij}^{n*} v_i^n)$, perform the collision according to (\ref{eq_rule_BDG}), 
where again, $w_i^n$ and $w_j^n$ are independent identically distributed random variables distributed according to the probability $g(w)$. 

\item[-] With probability $1- H(\hat v_{ij}^{n*} v_i^n)$, ignore the collision, i.e. define
$$ v_i^{n+1} = v_i^n ,  \quad v_j^{n+1} = v_j^n . $$

\end{itemize}

The phase of the quantity $\hat v_{ij}^{n*} v_i^n$ is the angle between $\hat v_{ij}^{n}$ and $v_i^n$. The biased BDG dynamics performs the collision with a probability depending on this angle. For instance, one can imagine that the collision will be performed more frequently if this angle is small than if it is large. It is a straightforward matter to notice that $\mbox{Re}(\hat v_{ij}^{n*} v_i^n) \geq 0$. This point will be proved below. This is why the function $H$ needs only be defined on ${\mathbb S}^1_+$. It is also easy to see that $\hat v_{ij}^{n*} v_j^n = ( \hat v_{ij}^{n*} v_i^n)^*$. Since $H$ is symmetric, the probability $H(\hat v_{ij}^{n*} v_i^n)$ is invariant under exchange of $i$ and $j$. We also assume that the time unit is chosen such that 
$$ \max_{w \in {\mathbb S}^1_+} H(w) = 1. $$
In the sequel, BDG will refer to the biased-BDG dynamics except explicitly mentioned.

\subsection{The CL dynamics}
\label{subsec_def_CLD}

This dynamics is named after the acronym 'Choose the Leader'. It consists of the following steps

\begin{itemize}
\item[-]  pick an ordered pair $(i,j)$, ($i \not = j$) randomly according to a uniform distribution, i.e. with probability \begin{equation} P_{i,j} = \frac{1}{N(N-1)} . \label{eq_uniform_ordered} \end{equation}

\item[-] Define 
$$ v_i^{n+1} = w^n v_j^n,  \quad v_j^{n+1} = v_j^n . $$
where $w^n$ is a random variables on ${\mathbb S}^1$ distributed according to the probability $g(w)$. 

\end{itemize}

\setcounter{equation}{0}
\section{Master equations and kinetic hierarchy}
\label{sec_master_eq}

\subsection{Master equation}
\label{subsec_master_principles}

In this section, we recall the definition of the master equation of the particle system. We first introduce the Markov transition operator $Q_N$. It acts on functions $\Phi (\omega)$, with $\omega \in {\mathbb T}^N$, as follows:
$$ Q_N \Phi (\omega) = {\mathbb E} ( \Phi (\omega^{n+1}) \, | \, \omega^n = \omega) , $$
where ${\mathbb E}$ is the expectation over all stochastic processes involved in the dynamics sending $\omega^n$ to $\omega^{n+1}$. 

Let now $F_N^n(\omega)$ be the $N$-particle probability distribution function at iterate $n$. Then, by definition, $F_N^n$ is such that 
$$ {\mathbb E} ( \Phi (\omega^{n}) ) = \int F_N^n(\omega) \Phi (\omega) \, d \omega . $$
Therefore, by the standard properties of conditional expectations, we have: 
\begin{eqnarray} 
{\mathbb E} ( \Phi (\omega^{n+1}) ) &=& \int F_N^{n+1}(\omega) \Phi (\omega) \, d \omega \nonumber \\
&=& {\mathbb E} \big( {\mathbb E} ( \Phi (\omega^{n+1}) \, | \, \omega^n) \big) \nonumber \\
&=& \int (Q_N \Phi) (\omega) \, F_N^{n}(\omega) \, d \omega \nonumber \\
&=& \int \Phi (\omega) \,  (Q_N^* F_N^{n}) (\omega) \, d \omega, \label{eq_cond_exp}
\end{eqnarray}
where $Q_N^*$ is the adjoint operator to $Q_N$. 
The dynamics of $F_N^n$ is thus: 
$$ F_N^{n+1} = Q_N^* F_N^{n} . $$

To pass to a continuous-in-time dynamics, we assume that the collision times for a given particle occur according to a Poisson stream with rate $\nu$. Since there are $N$ particles, the total collision frequency is of the order of $N \nu$. Then, the time-dependent $N$-particle distribution function $F_N = F_N(\omega,t)$ satisfies the following master equation:
\begin{equation} \frac{\partial}{\partial t} F_N = \nu N \, L_N^* F_N, \quad \quad  L_N^* = Q_N^* - \mbox{Id} . 
\label{eq_master_eq}
\end{equation}
We refer e.g. to \cite{CDW} for details. The weak form of the master equation is given for any test function 
$\Phi \in C^0({\mathbb T}^N)$ by: 
\begin{eqnarray} 
& & \hspace{-1cm}  \frac{\partial}{\partial t} \int_{{\mathbb T}^N} F_N(\omega,t) \, \Phi(\omega) \, d\omega = \nu N \, \int_{{\mathbb T}^N}   F_N(\omega,t) \, L_N \Phi(\omega) \, d\omega , \quad \quad  L_N = Q_N - \mbox{Id} . 
\label{eq_master_eq_weak}
\end{eqnarray}
 Because the particles are identical, and they cannot be
  ordered in a natural way, it is natural to assume that the initial
  distribution is invariant under permuation of the indices, and this
  invariance is then preserved by the dynamics.

The derivations of the master equations for the BDG and CL dynamics
are performed below. Before this, we recall the notion of hierarchy
and  propagation of chaos in the next section.

\subsection{Kinetic hierarchy and propagation of chaos}
\label{subsec_chaos}

We first define the $k$-particle marginal $F_{N,k}$ of $F_N$. For $k \in \{1, \ldots, N\}$, $F_{N,k} (v_1, \ldots, v_k)$ is defined by
$$ 
F_{N,k} (v_1, \ldots, v_k) = \int_{(v_{k+1}, \ldots, v_N) \in {\mathbb T}^{N-k}} F_N(v_1, \ldots, v_k, v_{k+1}, \ldots, v_N) \, dv_{k+1}, \ldots, dv_N . 
$$
By permutation invariance, without loss of generality, we can choose
to integrate out the last $N-k$ variables only and the resulting
$F_{N,k}$ is also permutationally invariant. The equation for
$F_{N,k}$ is found by integrating (\ref{eq_master_eq}) over $(v_{k+1},
\ldots, v_N)$. In general, the right-hand side of the resulting
equation involves higher order marginals, such as
$F_{N,k+1}$. Therefore, the equations for the marginals are all
coupled to each other, forming the so-called BBGKY hierarchy (see
\cite{CIP} for instance). The BBGKY hierarchy is a key ingredient in
the investigation of the $N \to \infty$ limit because the number of
variables involved in a given marginal is fixed. By contrast, the
number of variables involved in $F_N$ equals $N$ and goes to infinity
with $N$, which makes the functional treatment more complex. As long
as $N$ is finite, the BBGKY hierarchy is finite (the number of
equations in the hierarchy is equal to $N$) and does not carry more
information than the master equation itself. However, in the limit
$N\to \infty$, the hierarchy becomes infinite and is called the
Boltzmann hierarchy or kinetic hierarchy. There is no master equation
any more and the kinetic hierarchy is the only object that carries the
information about the process.  

Of course, it is desirable to break the hierarchy into a finite number
of equations. For this purposes, one says that the sequence $F_{N,k}$
satisfies the chaos property (or is $f$-chaotic) if there exists a
function $f(v)$ of the single variable $v$ such that

\begin{equation}
  \label{eq:propc}
  F_{N,k} (v_1, \ldots, v_k) \rightarrow \prod_{j=1}^k f(v_j), \quad
\mbox{ as } \quad N \to \infty,
\end{equation}
in the weak star topology of measures, for each $k \in {\mathbb
  N}$. This expresses that the $N$-particle probability $F_N$
approaches a product probability as $N$ becomes large, and translates
the fact that the particles become nearly independent in this limit.  

For a solution of the master equation (\ref{eq_master_eq}), one can
only expect this property to be true if at least the initial condition
satisfies it. By ``propagation of chaos'' we mean that if
(\ref{eq:propc}) holds for the initial data, {\em i.e.} that
$(F_N)|_{t=0}$ is $f_0$ chaotic for some function $f_0$, then for all
$t>0$, there is a function $f(\cdot,t)$ such that $F_N(\cdot,t)$ is
$f(\cdot,t)$-chaotic. In general, the rate of convergence in
(\ref{eq:propc}) depends on $k$ and on $t$, and only in particular
cases can one hope for uniform in time estimates
(see \cite{MischlerMouhot},\cite{MMW}). For a similar class of Markov processes
with applications to biology, Lachowicz has proven $L^1$-convergence
with bounds of the form $C N^{-\eta(T)}$, with a function $\eta(T)$
could be decaying exponentially fast with $T$  \cite{Lachowicz}.

%
%
%
%
%
%

If  propagation of chaos holds, then, as $N \to \infty$, one can replace $F_{N,k}(t)$ by the product
$\prod_{j=1}^k f(v_j,t)$ in the hierarchy and get a closed equation
for $f(v,t)$. The resulting equation for $f(v,t)$ is a kinetic
equation. Combining the BBGKY hierarchy and a  propagation of chaos
result is one of the ways one can derive kinetic equations from
N-particle systems (see e.g. \cite{Lanford3} in the case of the
Boltzmann equation or \cite{Kac} for Kac's equation).  

In the next sections, we derive the hierarchy for both the BDG and CL
processes.

\subsection{Master equation and hierarchy for the BDG dynamics}
\label{subsec_master_BBDG}

The following proposition establishes the master equations for the BDG dynamics. The master equation for the BDG dynamics (i.e. when $H$ is identically equal to $1$) as been previously established in \cite{CDW}. We assume that the Lebesgue measure $dv$ on ${\mathbb S}^1$ is normalized so that $\int_{{\mathbb S}^1} dv = 1$. 

\begin{proposition}
The master equation for the BDG dynamics is given by
\begin{eqnarray} 
& & \hspace{-1cm}   \frac{\partial F_N}{\partial t} (v_1, \ldots , v_N,t)= \frac{2 \nu}{N-1} \sum_{i<j}   \nonumber\\
& & \hspace{-0.5cm}  \left\{ 2 \int_{(u,z) \in {\mathbb S}^1 \times {\mathbb S}^1_+} H(z) \,   g(u^* v_i ) \, g(u^* v_j) \,  F_N (v_1, \ldots , uz, \ldots ,  uz^*, \ldots , v_N)\, du \, dz \right. \nonumber \\
& & \hspace{7cm} \left. \phantom{\int_{{\mathbb T}^2}}   -   H(\hat v_{ij}^* v_i)  \,  F_N (v_1, \ldots , v_N)  \right\} , 
\label{eq_mas_BBDG_strong}
\end{eqnarray}
where $uz$ and $uz^*$ are on the $i$-th and $j$-th positions respectively. 
\label{prop_mas_eq_BBDG}
\end{proposition}

\medskip
\noindent
The proof of this Lemma is given in Appendix A, section \ref{sec_appA_proof_master_BBDG}. We now turn to the BBGKY hierarchy and state the:

\medskip
\begin{proposition}
Let $k \in \{1, \ldots, N\}$. The $k$-particle marginal $F_{N,k}$ of the solution of the BDG master equation (\ref{eq_mas_BBDG_strong}) satisfies: 
\begin{eqnarray} 
& & \hspace{-1cm}   \frac{\partial F_{N,k}}{\partial t} (v_1, \ldots , v_k,t)= \frac{2 \nu}{N-1} \left[ \phantom{\int_{{\mathbb T}^2}} \right.  \nonumber\\
& & \hspace{-0.5cm} \sum_{i<j\leq k} \left\{ 2 \int_{(u,z) \in {\mathbb S}^1 \times {\mathbb S}^1_+} H(z)  \,   g(u^* v_i ) \, g(u^* v_j) \,  F_{N,k} (v_1, \ldots , uz, \ldots ,  uz^*, \ldots , v_k)\, du \, dz \right. \nonumber \\
& & \hspace{5cm} \left. \phantom{\int_{{\mathbb T}^2}}   -    H(\hat v_{ij}^* v_i) \,  F_{N,k} (v_1, \ldots , v_k)  \right\} \\
& & \hspace{-0.5cm} + (N-k) \sum_{i\leq k} \left\{ 2 \int_{(u,z) \in {\mathbb S}^1 \times {\mathbb S}^1_+} H(z)  \,   g(u^* v_i ) \,  F_{N,k+1} (v_1, \ldots , uz, \ldots , v_k, uz^*)\, du \, dz \right. \nonumber \\
& & \hspace{5cm} \left. \left. - \int_{{\mathbb S}^1}  H(\hat v_{i\, k+1}^* v_i) \,  F_{N,k+1} (v_1, \ldots , v_{k+1})  \, dv_{k+1} \right\} \, \right] .
\label{eq_hierarchy_BBDG_strong}
\end{eqnarray}
\label{prop_hierarchy_BBDG}
\end{proposition}

\medskip
\noindent
As examples, we write the first two elements of the hierarchy. For the one-particle marginal equation, the first sum is empty and the only remaining term corresponds to the choice $j=1$ in the second sum. Therefore, the equation is written: 
\begin{eqnarray} 
& & \hspace{-1cm}   \frac{\partial F_{N,1}}{\partial t} (v_1,t) = 2 \nu \left\{ 2 \int_{(u,z) \in {\mathbb S}^1 \times {\mathbb S}^1_+} H(z)  \,   g(u^* v_1 ) \,  F_{N,2} (uz, uz^*)\, du \, dz \right. \nonumber \\
& & \hspace{5cm} \left. - \int_{{\mathbb S}^1}  H(\hat v_{1\, 2}^* v_1) \,  F_{N,2} (v_1, v_2) \, dv_2 \right\} \, .
\label{eq_hierarchy_BBDG_strong_1}
\end{eqnarray}
For the two-particle marginal, there is only one term from the first sum, corresponding to $(i,j) = (1,2)$ and two terms from the second sum corresponding to the choices $i=1$ or $i=2$. This leads to 
\begin{eqnarray*} 
& & \hspace{-1cm}   \frac{\partial F_{N,2}}{\partial t} (v_1, v_2,t)= \frac{2 \nu}{N-1} \left[ \phantom{\int_{{\mathbb T}^2}}  2 \int_{(u,z) \in {\mathbb S}^1 \times {\mathbb S}^1_+} H(z)  \,   g(u^* v_1 ) \, g(u^* v_2) \,  F_{N,2} ( uz, uz^*)\, du \, dz \right. \nonumber \\
& & \hspace{8cm}    -    H(\hat v_{12}^* v_1) \,  F_{N,2} (v_1, v_2)   \\
& & \hspace{1.5cm} + (N-2) \left\{ 2 \int_{(u,z) \in {\mathbb S}^1 \times {\mathbb S}^1_+} H(z)  \,   g(u^* v_1 ) \,  F_{N,3} (uz, v_2, uz^*)\, du \, dz \right. \nonumber \\
& & \hspace{6cm}  - \int_{{\mathbb S}^1}  H(\hat v_{1\, 3}^* v_1) \,  F_{N,3} (v_1, v_2, v_3)  \, dv_{3}  . \\
& & \hspace{3.5cm} + 2 \int_{(u,z) \in {\mathbb S}^1 \times {\mathbb S}^1_+} H(z)  \,   g(u^* v_2 ) \,  F_{N,3} (v_1, uz, uz^*)\, du \, dz \nonumber \\
& & \hspace{6cm} \left. \left. - \int_{{\mathbb S}^1}  H(\hat v_{2\, 3}^* v_2) \,  F_{N,3} (v_1, v_2, v_3)  \, dv_{3} \right\} \, \right] .
\end{eqnarray*}
As anticipated, the hierarchy is not closed. Each level $k$ requires the knowledge of the next level $k+1$. If  propagation of chaos holds, i.e. if 
\begin{equation} 
F_{N,2} (v_1,v_2) \approx  F_{N,1} (v_1) \, F_{N,1}(v_2), 
\label{eq_approx_chaos}
\end{equation}
then, (\ref{eq_approx_chaos}) can be substituted into (\ref{eq_hierarchy_BBDG_strong_1}) and leads to 
\begin{eqnarray} 
& & \hspace{-1cm}   \frac{\partial F_{N,1}}{\partial t} (v_1,t) = 2 \nu \left\{ 2 \int_{(u,z) \in {\mathbb S}^1 \times {\mathbb S}^1_+} H(z)  \,   g(u^* v_1 ) \,  F_{N,1} (uz) \, F_{N,1}( uz^*)\, du \, dz \right. \nonumber \\
& & \hspace{5cm} \left. - \left( \int_{{\mathbb S}^1}  H(\hat v_{1\, 2}^* v_1)  \, F_{N,1}(v_2) \, dv_2 \right) \,  F_{N,1} (v_1) \right\} \, .
\label{eq_BBDG_kinetic}
\end{eqnarray}
This is the kinetic equation proposed in \cite{BDG}. The question to be investigated  is whether the approximation (\ref{eq_approx_chaos}) can be used. 

It can be seen from (\ref{eq_mas_BBDG_strong}) that the master equation can be put in the form 
\begin{eqnarray} 
& & \hspace{-1cm}   \frac{\partial F_N}{\partial t} =  \frac{2 \nu}{N-1} \sum_{i<j} (Q_{(i,j)}^*  - I) F_N, 
\label{eq_binary}
\end{eqnarray}
where $I$ is the identity and $Q_{(i,j)}^*$ is the following binary collision operator: 
\begin{eqnarray*} 
& & \hspace{-1cm}   Q_{(i,j)}^* F_N (v_1, \ldots , v_N)= \\
& & \hspace{-0.5cm}  =2 \int_{(u,z) \in {\mathbb S}^1 \times {\mathbb S}^1_+} H(z) \,   g(u^* v_i ) \, g(u^* v_j) \,  F_N (v_1, \ldots , uz, \ldots ,  uz^*, \ldots , v_N)\, du \, dz + \\
& & \hspace{7.5cm}
+ (1 - H(\hat v_{ij}^{*} v_i)) \, F_N (v_1, \ldots , v_N).
\end{eqnarray*}
Its adjoint given by
\begin{eqnarray*} 
& & \hspace{-1cm}  
Q_{(i,j)}^* \Phi(v_1, \ldots , v_N)  = \\
& & \hspace{-0.5cm}  
=  H(\hat v_{ij}^{*} v_i) \int_{(v'_i,v'_j) \in {\mathbb T}^2}  g(\hat v^{*}_{ij} v'_i ) \, g(\hat v^{*}_{ij} v'_j)  \, \Phi (v_1, \ldots , v'_i, \ldots ,  v'_j, \ldots , v_N) \, dv'_i \, dv'_j 
+ \\
& & \hspace{7.5cm}
+ (1 - H(\hat v_{ij}^{*} v_i)) \, \Phi (v_1, \ldots , v_N),
\end{eqnarray*}
is a Markovian operator operating on $\Phi$ through $v_i$ and $v_j$ alone. In \cite{CDW}, the general framework of master equations of the type (\ref{eq_binary}), called pair-interaction driven master equations, is investigated. It is proved that  propagation of chaos holds on any finite time interval $[0,T]$. However, it is not known if  propagation of chaos holds uniformly in time and in particular, if the invariant measure (i.e. the equilibrium $F_{N,\infty}$ corresponding to $\frac{\partial F_N}{\partial t} = 0$ in (\ref{eq_binary})) is chaotic. The case $H \equiv 1$ corresponds to the unbiased BDG model and has been investigated in \cite{CDW}. 

\begin{remark}
For computational purposes, the weak forms of the master equation and hierarchy are more convenient. The weak form of the master equation (\ref{eq_mas_BBDG_strong}) is : 
\begin{eqnarray} 
& & \hspace{-1cm}  \frac{\partial}{\partial t} \int_{{\mathbb T}^N} F_N(v_1, \ldots , v_N,t) \, \Phi(v_1, \ldots , v_N) \, dv_1 \ldots dv_N = \frac{2 \nu}{N-1}  \int_{(v_1, \ldots , v_N) \in {\mathbb T}^n}  \nonumber \\
& & \hspace{-0.5cm}
\sum_{i<j} H(\hat v_{ij}^{*} v_i) \left\{ \int_{(v'_i,v'_j) \in {\mathbb T}^2}  g(\hat v^{*}_{ij} v'_i ) \, g(\hat v^{*}_{ij} v'_j)  \, \Phi (v_1, \ldots , v'_i, \ldots ,  v'_j, \ldots , v_N) \, dv'_i \, dv'_j \right. \nonumber  \\
& & \hspace{2cm}  \left. \phantom{\int_{(v'_i,v'_j) \in {\mathbb T}^2}} - \, \Phi (v_1, \ldots , v_N) \,  \right\} F_N (v_1, \ldots , v_N) \, dv_1 \ldots  dv_N  ,
\label{eq_mas_BBDG_weak}
\end{eqnarray}
for any continuous test function $\Phi(v_1, \ldots , v_N)$ on ${\mathbb T}^N$. The weak form of the hierarchy (\ref{eq_hierarchy_BBDG_strong}) is written as follows: 
\begin{eqnarray} 
& & \hspace{-1cm}  \frac{\partial}{\partial t} \int_{{\mathbb T}^k} F_{N,k}(v_1, \ldots , v_k,t) \, \Phi(v_1, \ldots , v_k) \, dv_1 \ldots dv_k = \frac{2 \nu}{N-1} \int_{(v_1, \ldots , v_k) \in {\mathbb T}^k} \left[ \phantom{\int_{{\mathbb T}^2}} \right. \nonumber \\
& & \hspace{-0.5cm}
 \sum_{i<j\leq k} H(\hat v_{ij}^{*} v_i) \left\{  \int_{(v'_i,v'_j) \in {\mathbb T}^2}  g(\hat v^{*}_{ij} v'_i ) \, g(\hat v^{*}_{ij} v'_j) \, \Phi (v_1, \ldots , v'_i, \ldots ,  v'_j, \ldots , v_k) \, dv'_i \, dv'_j  \right. \nonumber  \\
& & \hspace{1cm} \left. \phantom{\int_{(v'_i,v'_j) \in {\mathbb T}^2}} - \Phi (v_1, \ldots , v_k) \, \right\}   F_{N,k} (v_1, \ldots , v_k) \nonumber \\
& & \hspace{-0.5cm}
+ (N-k) \int_{v_{k+1} \in {\mathbb S^1}} \sum_{i\leq k} H(\hat v_{i \, k+1}^{*} v_i) \left\{  \int_{v'_i \in {\mathbb S}^1}  g(\hat v^{*}_{i \, k+1} v'_i ) \, \Phi (v_1, \ldots , v'_i, \ldots , v_k) \, dv'_i   \right. \nonumber  \\
& & \hspace{1cm} \left. \left. \phantom{\int_{v'_i \in {\mathbb S}^1}} - \Phi (v_1, \ldots , v_k) \, \right\}   F_{N,k+1} (v_1, \ldots , v_k, v_{k+1}) \, dv_{k+1} \right] \, dv_1 \ldots  dv_k 
,
\label{eq_hierarchy_BBDG_weak}
\end{eqnarray}
for any continuous test function $\Phi(v_1, \ldots , v_k)$ on ${\mathbb T}^k$. 
\label{rem_BBDG_weak}
\end{remark}

\subsection{Master equation and hierarchy for the CL dynamics}
\label{subsec_master_CLD}

Before stating the result, we introduce some notations: We write $(v_1, \ldots , \hat v_i, \ldots , v_N)$ for $(v_1, \ldots , v_{i-1}, v_{i+1}, \ldots , v_N)$, i.e. we mean that $v_i$ is absent from the list. We also define:
\begin{eqnarray} 
& & \hspace{-1cm}  [ F_N ]_{\hat i} (v_1, \ldots , \hat v_i, \ldots , v_N) = \int_{v_i \in {\mathbb S}^1} F(v_1, \ldots , v_i, \ldots , v_N) \, dv_i.
\label{eq_def_[FN]}
\end{eqnarray}

\begin{proposition}
The master equation for the CL dynamics is given by
\begin{eqnarray} 
& & \hspace{-1cm}   \frac{\partial F_N}{\partial t}(v_1, \ldots , v_N,t) =  \nonumber \\
& & \hspace{-0.5cm}  = \frac{2 \nu}{N-1} \sum_{i<j}  \left\{  \frac{1}{2} g(v_i^* v_j ) \, \left( [ F_N ]_{\hat j} (v_1, \ldots , \hat v_j, \ldots , v_N) +  \phantom{\frac{1}{2}}  \right. \right. \nonumber \\
& & \hspace{3.5cm}   \left. \left. \phantom{\frac{1}{2}} +  [ F_N ]_{\hat i} (v_1, \ldots , \hat v_i, \ldots , v_N) \right) - F_N (v_1, \ldots , v_N) \phantom{\frac{1}{2}} \hspace{-0.3cm} \right\}.
\label{eq_mas_CLD_strong}
\end{eqnarray}
 
\label{prop_mas_eq_CLD}
\end{proposition}

\noindent
The proof of this proposition can be found in \cite{CDW}. We reproduce it in Appendix A, section \ref{sec_appA_proof_master_CLD} for the reader's convenience. We now consider the BBGKY  hierarchy. We have the:

\begin{proposition}
Let $k \in \{1, \ldots, N\}$. The $k$-particle marginal $F_{N,k}$ of the solution of the CL master equation (\ref{eq_mas_CLD_strong}) satisfies: 
\begin{eqnarray} 
& & \hspace{-1cm}   \frac{\partial F_{N,k}}{\partial t}(v_1, \ldots , v_k,t) =  \nonumber \\
& & \hspace{-0.5cm}  = \frac{2 \nu}{N-1} \left[ \sum_{i<j\leq k}  \left\{  \frac{1}{2} g(v_i^* v_j ) \, \left( F_{N,k-1}  (v_1, \ldots , \hat v_i, \ldots , v_j, \ldots,  v_k) + \phantom{\frac{1}{2}} \right. \right. \right. \nonumber \\
& & \hspace{3cm} \left. \left. \phantom{\frac{1}{2}}  +  F_{N,k-1}  (v_1, \ldots , v_i, \ldots , \hat v_j, \ldots,  v_k) \right) - F_{N,k} (v_1, \ldots , v_k) \right\} \nonumber \\
& & \hspace{1.2cm}  + (N-k) \sum_{i \leq k}   \frac{1}{2} \left\{ \int_{v_{k+1} \in {\mathbb S}^1} g(v_{k+1}^* v_i ) \, F_{N,k}  (v_1, \ldots , \hat v_i, \ldots ,  v_k, v_{k+1}) \, dv_{k+1} \right.  \nonumber \\
& & \hspace{9cm}  \left. \left. - F_{N,k} (v_1, \ldots , v_k)  \phantom{\int_{v_{k+1} \in {\mathbb S}^1}} \hspace{-1.3cm} \right\} \, \,  \right]
.
\label{eq_hierarchy_CLD_strong}
\end{eqnarray}
\label{prop_hierarchy_CLD}
\end{proposition}

\noindent
As examples, we write the first two elements of the hierarchy (they have been previously established in \cite{CDW}). For the one-particle marginal equation, we get: 
\begin{eqnarray} 
& & \hspace{-1cm}   \frac{\partial F_{N,1}}{\partial t} (v_1,t) = \nu \left( \int_{v_{2} \in {\mathbb S}^1} g(v_{2}^* v_1 ) \, F_{N,1}  (v_2) \, dv_2  - F_{N,1} (v_1) \right), \label{eq_hierarchy_CLD_strong_1}
\end{eqnarray}
and for the two-particle marginal, we have: 
\begin{eqnarray*} 
& & \hspace{-0.5cm}   \frac{\partial F_{N,2}}{\partial t} (v_1, v_2,t) = \frac{2 \nu}{N-1} \left[ \left\{  \frac{1}{2} g(v_1^* v_2 ) \, \big( F_{N,1}  (v_1) + F_{N,1}  (v_2) \big) - F_{N,2} (v_1, v_2) \right\} \right. \nonumber \\
& & \hspace{0.cm}  + (N-2) \left(    \frac{1}{2} \left\{ \int_{v_{3} \in {\mathbb S}^1} g(v_{3}^* v_1 ) \, F_{N,2}  (v_2, v_3) \, dv_3 + \int_{v_{3} \in {\mathbb S}^1} g(v_{3}^* v_2 ) \, F_{N,2}  (v_1, v_3) \, dv_3  \right\}  \right. \nonumber \\
& & \hspace{11cm}  \left. \left.  - F_{N,2} (v_1,  v_2)  \phantom{\int_{v_3 \in {\mathbb S}^1}} \hspace{-1cm}  \right) \,  \right]. 
\end{eqnarray*}

By contrast to the BDG hierarchy, the CL hierarchy is closed at any order. This is a very remarkable feature of this model, due to the fact that the pair interaction only acts on one of the variables. The CL master equation (\ref{eq_mas_CLD_strong}) can be put in the frame of pair-interaction driven master equations (\ref{eq_binary}), with 
\begin{eqnarray*} 
& & \hspace{-1cm}   Q_{(i,j)} F_N (v_1, \ldots , v_N) =  \nonumber \\
& & \hspace{0.5cm}  =   \frac{1}{2} g(v_i^* v_j ) \, \left( [ F_N ]_{\hat j} (v_1, \ldots , \hat v_j, \ldots , v_N) +  [ F_N ]_{\hat i} (v_1, \ldots , \hat v_i, \ldots , v_N) \phantom{\frac{1}{2}} \hspace{-0.3cm} \right) .
\end{eqnarray*}
Indeed, its adjoint
\begin{eqnarray*} 
& & \hspace{-1cm}  Q_{(i,j)}^* \Phi(v_1, \ldots , v_N)  = \nonumber \\
& & \hspace{-0.5cm} =  \int_{z \in {\mathbb S}^1}  \frac{1}{2} \left( \Phi (v_1, \ldots , v_i, \ldots ,  z v_i, \ldots , v_N) +  \Phi (v_1, \ldots , z v_j, \ldots ,  v_j, \ldots , v_N) \phantom{\frac{1}{2}} \hspace{-0.3cm} \right)\, g(z) \, dz  ,
\end{eqnarray*}
is a Markovian operator acting through $v_i$ and $v_j$ alone. Therefore, the result of \cite{CDW} applies and  propagation of chaos is true on any finite time interval. However, again, it is not known if propagation of chaos is valid uniformly in time or breakdowns at large times. 

In \cite{CDW}, thanks to the closed hierarchy, an analytical formula for the marginals of the equilibrium density $F_{N,\infty}$ is given. It is shown that, if the noise $g$ is  properly rescaled with $N$, the equilibrium density is {\bf not} chaotic. In sections \ref{sec_grazing} and \ref{sec_numerics}, we revisit this example with a special choice of the noise rescaling and we illustrate the loss of chaos numerically. This counter-example is not in contradiction with the previous result of \cite{CDW} because of the noise rescaling on the one hand and of of the large time scales on the other hand. 

\begin{remark}
Again, we give the weak forms of the master equation and hierarchy of the CL dynamics, which are useful for computational purposes. The master equation (\ref{eq_mas_CLD_strong}) is given in weak form: 
\begin{eqnarray} 
& & \hspace{-0.5cm}  \frac{\partial}{\partial t} \int_{{\mathbb T}^N} F_N(v_1, \ldots , v_N,t) \, \Phi(v_1, \ldots , v_N) \, dv_1 \ldots dv_N = \frac{2 \nu}{N-1} \sum_{i<j} \int_{(v_1, \ldots , v_N) \in {\mathbb T}^n}  
\left\{ \phantom{\int_{z \in {\mathbb S}^1}} \right. \nonumber \\
& & \hspace{0.cm}  \int_{z \in {\mathbb S}^1}  \frac{1}{2} \left( \Phi (v_1, \ldots , v_i, \ldots ,  z v_i, \ldots , v_N) +  \Phi (v_1, \ldots , z v_j, \ldots ,  v_j, \ldots , v_N) \phantom{\frac{1}{2}} \hspace{-0.3cm} \right)\, g(z) \, dz \nonumber \\
& & \hspace{3cm} \left. \phantom{\int_{z \in {\mathbb S}^1}  \frac{1}{2}}  - \Phi(v_1, \ldots, v_N) \right\}  \, F_N (v_1, \ldots , v_N) \, dv_1 \ldots  dv_N  ,
\label{eq_mas_CLD_weak}
\end{eqnarray}
for any continuous test function $\Phi(v_1, \ldots , v_N)$ on ${\mathbb T}^N$. The weak form of the hierarchy (\ref{eq_hierarchy_CLD_strong}) is as follows: 
\begin{eqnarray} 
& & \hspace{-0.5cm}  
\frac{\partial}{\partial t} \int_{{\mathbb T}^N} F_{N,k}(v_1, \ldots , v_k,t) \, \Phi(v_1, \ldots , v_k) \, dv_1 \ldots dv_k = \frac{2 \nu}{N-1}  \left[ \sum_{i<j \leq k} \int_{(v_1, \ldots ,\hat v_i, \ldots, \hat v_j, \ldots, v_k) \in {\mathbb T}^{k-2}} \left\{
\phantom{\sum_{i<j \leq k}} \right. \right.   \nonumber \\
& & \hspace{0cm} 
\frac{1}{2} \int_{(v_i,z) \in {\mathbb T}^2} \Phi (v_1, \ldots , v_i, \ldots , z v_i, \ldots , v_k) F_{N,k-1} (v_1, \ldots , v_i, \ldots, \hat v_j, \ldots ,  v_k) \, g(z) \, dz \, dv_i
\nonumber \\
& & \hspace{0cm} 
+ \frac{1}{2} \int_{(v_j,z) \in {\mathbb T}^2}  \Phi (v_1, \ldots , z v_j, \ldots ,  v_j, \ldots , v_k) F_{N,k-1} (v_1, \ldots, \hat v_i, \ldots , v_j, \ldots ,  v_k) \, g(z) \, dz \, dv_j \nonumber \\
& & \hspace{6cm} \left. - \int_{(v_i,v_j) \in {\mathbb T}^2} \Phi(v_1, \ldots, v_k)  \, F_{N,k} (v_1, \ldots , v_k) \, dv_i \, dv_j \right\} \nonumber \\
& & \hspace{0cm} 
+ (N-k) \sum_{i\leq k} \int_{(v_1, \ldots ,\hat v_i, \ldots, v_k) \in {\mathbb T}^{k-1}} \frac{1}{2} \left\{ \phantom{\int_{(v_{k+1},z) \in {\mathbb T}^2}} \right. \nonumber \\
& & \hspace{0.5cm} 
\int_{(v_{k+1},z) \in {\mathbb T}^2} \Phi (v_1, \ldots , z v_{k+1}, \ldots , v_k) \, F_{N,k} (v_1, \ldots ,\hat v_i, \ldots, v_k, v_{k+1}) \, g(z) \, dz \, dv_{k+1} \nonumber \\
& & \hspace{5.5cm} \left. \left. \phantom{\int_{z \in {\mathbb S}^1}}  - \int_{v_i \in {\mathbb S}^1} \Phi(v_1, \ldots, v_k)  \, F_{N,k} (v_1, \ldots , v_k) \, dv_i \,  \right\} \, \,  \right] ,
\label{eq_hierarchy_CLD_weak}
\end{eqnarray}
for any continuous test function $\Phi(v_1, \ldots , v_k)$ on ${\mathbb T}^k$. 
\label{rem_CLD_weak}
\end{remark}

\setcounter{equation}{0}
\section{Rescaled hierarchies and the limit $N \to \infty$}
\label{sec_grazing}

\subsection{Noise rescaling}
\label{subsec_grazing_into}

Large time scales are not covered by the  propagation of chaos result of \cite{CDW}. The goal of this section is to investigate whether  propagation of chaos is still valid for the BDG and CL dynamics at large time scales or not. To do so, it is necessary to rescale the noise distribution $g$ (and, in the case of the BDG dynamics, the bias function $H$). Indeed, if we rescale time to large time scales, we simultaneously need to rescale the collision operators in order to keep the leading order terms finite. The appropriate scalings of $g$ and $H$ correspond to respectively a small noise intensity and grazing collision asymptotics. 

We suppose that the noise probability distribution $g(v)$ depends on a small parameter $\varepsilon$ and we denote it by $g_\varepsilon$. The parameter $\varepsilon$ will be linked to $N$ in such a way that $\varepsilon \to 0$ as $N \to \infty$. Similarly, we assume that $H=H_\varepsilon(v)$ depends on $\varepsilon$ and we introduce the probability distribution 
\begin{equation} 
h_\varepsilon(v) = \frac{H_\varepsilon(v)}{{\mathcal H}_\varepsilon}, \quad \quad {\mathcal H}_\varepsilon = \int_{{\mathbb S}^1_+} H_\varepsilon(v) \, dv . 
\label{eq_rescale_h}
\end{equation}
We keep the assumption that 
\begin{equation} \max_{w \in {\mathbb S}^1_+} H_\varepsilon(w) = 1, \forall \varepsilon >0.  
\label{eq_ass_Heps}
\end{equation}

It will be more convenient to introduce the phases of the velocities, i.e. we will write 
\begin{equation}
v = e^{i \theta}, \quad \mbox{ with } \quad \theta \in {\mathbb R}/2 \pi {\mathbb Z}.
\label{eq_phase_not}
\end{equation} 
We note that $dv = d\theta/(2 \pi)$. 
We assume that $g_\varepsilon$ and $h_\varepsilon$ are deduced from probability distribution functions $\gamma$ and $\eta$ defined on ${\mathbb R}$ by the following scaling relations: 

\begin{hypothesis}
We assume that 
\begin{eqnarray}
& & \hspace{-1cm} g_\varepsilon (\theta) = \frac{1}{\varepsilon \bar \gamma_\varepsilon} \gamma (\frac{\theta}{\varepsilon}), \quad \forall \theta \in [-\pi,\pi], 
\label{eq_def_gamma} \\
& & \hspace{-1cm} h_\varepsilon (\theta) = \frac{1}{\varepsilon \bar \eta_\varepsilon} \eta (\frac{\theta}{\varepsilon}), \quad \forall \theta \in [-\pi/2,\pi/2], 
\end{eqnarray}
where 
$$ \bar \gamma_\varepsilon = \int_{-\pi/\varepsilon}^{\pi/\varepsilon} \gamma(x) \, \frac{dx}{2 \pi}, \quad \bar \eta_\varepsilon = \int_{-\pi/(2\varepsilon)}^{\pi/(2\varepsilon)} \eta(x) \, \frac{dx}{\pi}, $$
and where $\gamma$ and $\eta$ are probability densities on ${\mathbb R}$ (for the measures $\frac{dx}{2 \pi}$ and $\frac{dx}{\pi}$ respectively) which are even and have finite second order moments $\sigma^2$ and $\tau^2$: 
\begin{equation}
\sigma^2 = \int_{-\infty}^{\infty} \gamma(x)\, x^2 \, \frac{dx}{2 \pi} <\infty, \quad  \tau^2 = \int_{-\infty}^{\infty}  \eta(x) \, x^2 \, \frac{dx}{\pi}<\infty. 
\label{eq_def_sigma_tau}
\end{equation}
\label{hyp_geps_heps}
\end{hypothesis}

\noindent
The limit $\varepsilon \to 0$ in $g_\varepsilon$ corresponds to a small noise intensity limit. In $h_\varepsilon$, $\varepsilon \tau$ represents the typical relative velocity at which collisions may happen. With relative velocities larger than $\varepsilon \tau$, the particles have very little probability to collide. Therefore, the limit $\varepsilon \to 0$ in $h_\varepsilon$ represents a grazing collision limit, in a way similar to the grazing collision limit of the Boltzmann equation \cite{DL, Desvil}. The hypothesis \ref{hyp_geps_heps} can be weakened but since the main purpose of this paper is illustrative, we do not seek the broadest generality. We stress that the statements given in the following sections are formal.

\begin{remark}
With hypothesis (\ref{hyp_geps_heps}), we have $ {\mathcal H}_\varepsilon = O(\varepsilon)$. Indeed, with the assumption (\ref{eq_ass_Heps}), we have 
$$ H_\varepsilon (\theta)  = \frac{1}{\max_{[-\pi/(2 \varepsilon),\pi/(2 \varepsilon)]} \eta} \, \, \eta (\frac{\theta}{\varepsilon}). $$
Therefore, 
$$ {\mathcal H}_\varepsilon \sim C_0 \varepsilon , \quad C_0 = \frac{1}{\max_{{\mathbb R}} \eta} \, . $$
\label{remark_Heps}
\end{remark}

\subsection{Rescaled BDG hierarchy and the limit $N \to \infty$}
\label{subsec_rescaled_BBDG}

We introduce the following operator: 
\begin{eqnarray} 
& & \hspace{-1cm} \bar D_{ij} = \frac{\partial^2}{\partial \theta_i^2} +  \frac{\partial^2}{\partial \theta_j^2} + 2 \frac{\partial^2}{\partial \theta_i \, \partial \theta_j} = \left( \frac{\partial}{\partial \theta_i} + \frac{\partial}{\partial \theta_j} \right)^2.
\label{eq_def_bar_D_ij}
\end{eqnarray} 

\medskip
\begin{theorem}
We assume that $\varepsilon=\varepsilon_N$ is linked to $N$ in such a way that $\varepsilon_N \to 0$ as $N \to \infty$. We rescale time in such a way that $t' = 2 {\mathcal H}_{\varepsilon_N} \varepsilon_N^2 t$. We assume that for any fixed $k$, $F_{N,k}$ converges in the weak star topology of measures as $N \to \infty$ towards a probability measure $F_{[k]}$, uniformly on any finite time interval. Additionally, we assume that $F_{N,k}$ remains bounded in $C^3(({\mathbb R}/(2 \pi {\mathbb Z}))^k)$ uniformly with respect to $N$, and with respect to time on any finite interval. Then, $F_{[k]}$ is a solution of the following infinite hierarchy: 
\begin{eqnarray} 
& & \hspace{-1cm} \frac{\partial F_{[k]}}{\partial t} (\theta_1, \ldots, \theta_k,t)  =  \nu \, (\sigma^2 - \tau^2) \, \sum_{i \leq k}
(\bar D_{i \, k+1} F_{[k+1]} ) (\theta_1, \ldots, \theta_i , \ldots , \theta_k, \theta_i).
\label{eq_BBGKY_BBDG_grazing}
\end{eqnarray} 
\label{thm_hierarchy_grazing_BBDG}
\end{theorem}

\noindent
The proof of this result is given in Appendix B, section \ref{sec_appB_proof_grazing_BBDG}. Of course, this theorem is formal because it is not known it $F_{N,k}$ satisfies the assumptions.

\medskip
In the limit $N \to \infty$, the hierarchy takes the form of an inductive sequence of heat-like equations on the $k$-dimensional torus. When $\sigma = \tau$, the evolution of $F_{[k]}$ takes place at a longer time scale. To find this evolution, one must expand the master equations to higher order terms in $1/N$. Such higher order expansions are beyond the scope of the present paper.

If the chaos assumption
\begin{equation}
F_{[k]}(\theta_1, \ldots \theta_k) = \prod_{i=1}^k F_{[1]}(\theta_i) ,
\label{eq_chaos_ass}
\end{equation}
is true, then $F_{[1]}=F_{[1]}(\theta,t)$ satisfies the nonlinear diffusion equation:
\begin{eqnarray}
& & \hspace{-1cm} \frac{\partial F_{[1]}}{\partial t} (\theta,t) = 2 (\sigma^2 - \tau^2) \frac{\partial}{\partial \theta} \left( F_{[1]} \frac{\partial F_{[1]}}{\partial \theta} \right).
\label{eq_BBDG_chaos}
\end{eqnarray}
This diffusion equation is of forward type (and thus, well-posed in the classical sense) if and only if $\sigma > \tau$. If $\sigma < \tau$, the diffusion equation is of backwards type and is only well-posed for specific initial conditions. In the case $\sigma > \tau$, the noise added after the interaction (measured by $\sigma$) is larger than the typical distance between colliding particles (measured by $\tau$). Thus, in average, the particles are further to each other after the collision than before it. The dynamics is then of diffusive type. Conversely, if $\sigma < \tau$ the particles are in average closer to each other after the collision than before. This dynamics produces concentrations at a rate which depends on the regularity of the initial data. However, it is not known if  propagation of chaos holds for this model and the validity of (\ref{eq_BBDG_chaos}) (even in the case $\tau > \sigma$) is subject to caution. 

We can look at the invariant densities, i.e. the stationary solutions of the hierarchy (\ref{eq_BBGKY_BBDG_grazing}). The uniform density
\begin{eqnarray} 
& & \hspace{-1cm} F_{[k], \mbox{\scriptsize eq}} = 1, \qquad \forall k \in {\mathbb N}^*, 
\label{eq_BBGKY_BBDG_grazing_equi}
\end{eqnarray} 
is an obvious equilibrium solution. It satisfies the chaos assumption, i.e. is a $k$-fold tensor product of the uniform single-particle marginal $F_{[1], \mbox{\scriptsize eq}} = 1$. It is not easy to see if this is the unique stationary solution of the hierarchy. The right-hand side of (\ref{eq_BBGKY_BBDG_grazing}) cancels functions of the form $\phi(\theta_i - \theta_k)$. Therefore, it is tempting to think that invariant densities similar to those of the CL dynamics (\ref{eq_Feq_1}) (see next section) exist. However, a quick check shows that it is not the case, unless for uniform distributions. If the uniform densities (\ref{eq_BBGKY_BBDG_grazing_equi}) are the unique equilibria of the hierarchy, this could be a hint that  propagation of chaos could be true for the BDG hierarchy. This would be in marked contrast with the CL hierarchy which is examined in the next section. However, the numerical simulations performed in section \ref{sec_numerics} seem to indicate that the two kinds of dynamics have a quite similar behavior. Therefore, we cannot make any conjecture whether  propagation of chaos holds for the hierarchy (\ref{eq_BBGKY_BBDG_grazing}).

\subsection{Rescaled CL hierarchy and the limit $N \to \infty$}
\label{subsec_grazing_BBGKY_CLD}

\begin{theorem}
We assume that $\varepsilon=\varepsilon_N$ is linked to $N$ by 
\begin{equation} \varepsilon_N  = \frac{1}{\sqrt N} . 
\label{eq_scaling_CLD}
\end{equation}
We also rescale time according to $t' = t/N$. We assume that, for any fixed $k$, $F_{N,k}$ converges in the weak star topology of measures towards a probability measure on $[0,2 \pi]^k$, as $N \to \infty$ uniformly on any finite time interval. Then, the limit $F_{[k]}$ satisfies the infinite hierarchy:
\begin{eqnarray} 
& & \hspace{-1cm}  \frac{\partial F_{[k]}}{\partial t} (\theta_1, \ldots, \theta_k,t)=  \sum_{i<j\leq k} \left[ \left( F_{[k-1]} (\theta_1, \ldots, \hat \theta_i, \ldots, \theta_j, \ldots, \theta_k) + \right. \right.
\nonumber \\
& & \hspace{0cm}  
\left. \left. 
+ F_{[k-1]} (\theta_1, \ldots, \theta_i, \ldots, \hat \theta_j, \ldots, \theta_k)  \right) \delta(\theta_i - \theta_j) - 2 F_{[k]}(\theta_1, \ldots, \theta_k) \right] + \nonumber \\
& & \hspace{8cm}  + \frac{\sigma^2 }{2} \sum_{i \leq k } 
\frac{\partial^2 F_{[k]}}{\partial \theta_i^2}(\theta_1, \ldots, \theta_k)  ,
\label{eq_BBGKY_CLD_grazing}
\end{eqnarray}
If $k=1$, only the second term (in factor of $\sigma^2$) remains. 
\label{thm_hierarchy_grazing_CLD}
\end{theorem}

\noindent
We notice that the assumptions on $F_{N,k}$ are weaker than in the BDG case and actually close to be satisfied. Indeed, for any given time, we can extract a subsequence $F_{N,k}$ which converges in the weak star topology of measures. What is lacking is some uniform time estimate which would allow the extraction of a single sequence on a whole time interval and the uniform weak convergence on this interval. 

The link between $\varepsilon$ and $N$ expressed by (\ref{eq_scaling_CLD}) is tighter than in the BDG case. This scaling allows to keep the largest possible number of non-zero terms in the hierarchy. Each level of the hierarchy involves a damped heat equation on the torus with a delta source term involving lower order terms of the hierarchy. Like in the finite $N$ case, the hierarchy is closed at any order. This remarkable feature allows us to show that this hierarchy does not satisfy the chaos property (\ref{eq_chaos_ass}). This is expressed in the following

\begin{theorem}
The CL hierarchy (\ref{eq_BBGKY_CLD_grazing}) does not satisfy the chaos property (\ref{eq_chaos_ass}).
\label{thm_chaos_CLD}
\end{theorem}

\medskip
\noindent
{\bf Proof:} The equations for the first and second marginals are respectively: 
\begin{eqnarray} 
& & \hspace{-1cm}  \frac{\partial F_{[1]}}{\partial t} =  \frac{\sigma^2 }{2} \frac{\partial^2 F_{[1]}}{\partial \theta_1^2}  ,
\label{eq_CLD_F1} \\
& & \hspace{-1cm}  \frac{\partial F_{[2]}}{\partial t} =  (F_{[1]} (\theta_1) + F_{[1]} (\theta_2) ) \delta(\theta_1 - \theta_2) - 2 F_{[2]} +  \frac{\sigma^2 }{2} \left( 
\frac{\partial^2 F_{[2]}}{\partial \theta_1^2} + \frac{\partial^2 F_{[2]}}{\partial \theta_2^2} \right)
.
\label{eq_CLD_F2} 
\end{eqnarray}
Now, let us suppose that  propagation of chaos holds. This means that for all solutions $F_{[1]}$ of (\ref{eq_CLD_F1}), the function $F_{[2]}(\theta_1,\theta_2,t) = F_{[1]}(\theta_1,t) F_{[1]}(\theta_2,t)$ must be a solution of (\ref{eq_CLD_F2}). However, for such a $F_{[2]}$, we have: 
\begin{eqnarray*} 
\frac{\partial F_{[2]}}{\partial t}(\theta_1,\theta_2,t) &=&  F_{[1]} (\theta_1,t) \frac{\partial F_{[1]}}{\partial t}(\theta_2,t) + F_{[1]} (\theta_2,t) \frac{\partial F_{[1]}}{\partial t}(\theta_1,t)\\
&=& \frac{\sigma^2 }{2} \left( F_{[1]} (\theta_1,t) \frac{\partial^2 F_{[1]}}{\partial \theta_1^2} (\theta_2,t) + F_{[1]} (\theta_2,t)  \frac{\partial^2 F_{[1]}}{\partial \theta_1^2} (\theta_1,t) \right) \\
&=& \frac{\sigma^2 }{2} \left( 
\frac{\partial^2 F_{[2]}}{\partial \theta_1^2} + \frac{\partial^2 F_{[2]}}{\partial \theta_2^2} \right) . 
\end{eqnarray*}
This implies that 
$$ (F_{[1]} (\theta_1,t) + F_{[1]} (\theta_2,t) ) \delta(\theta_1 - \theta_2) = 2 F_{[1]}(\theta_1,t) F_{[1]}(\theta_2,t) $$
which has $F_{[1]}(\theta,t) = 0$ for only solution. This shows that in general, the CL dynamics does not satisfy the chaos assumption. \endproof

\medskip
\noindent
We can precise what the first and second marginals of the equilibrium are.

\begin{proposition}
(i) The only stationary solution of (\ref{eq_CLD_F1}) is the isotropic measure 
\begin{equation}
F_{[1], \mbox{\scriptsize eq}} (\theta_1) = 1. 
\label{eq_equi_F1}
\end{equation}

\noindent
(ii) The only stationary solution of (\ref{eq_CLD_F2}) corresponding to the one-particle marginal (\ref{eq_equi_F1}) is of the form 
\begin{equation}
F_{[2], \mbox{\scriptsize eq}} (\theta_1, \theta_2) = {\mathcal M}(\theta_1 - \theta_2), 
\label{eq_equi_F2}
\end{equation}
where ${\mathcal M} (\theta)$ is the solution of 
\begin{equation}
- \frac{\sigma^2}{2} \frac{\partial {\mathcal M}}{\partial \theta^2} (\theta) + {\mathcal M}(\theta) = \delta(\theta). 
\label{eq_M_Helmholtz}
\end{equation}
${\mathcal M} (\theta)$ has the expression: 
\begin{equation}
{\mathcal M}(\theta) = \frac{1}{2 \bar \sigma (e^{2 \pi/\bar \sigma} - 1)} \left( e^{\theta/\bar \sigma} + e^{(2 \pi - \theta)/\bar \sigma} \right), \quad \theta \in [0,2 \pi], \quad \mbox{ with } \quad \bar \sigma = \frac{\sigma}{\sqrt 2}.  
\label{eq_express_M}
\end{equation}
\label{prop_equi_CLD}
\end{proposition}

\medskip
\noindent
{\bf Proof:} $F_{[2], \mbox{\scriptsize eq}}$ satisfies 
\begin{eqnarray*} 
& & \hspace{-1cm} - \frac{\sigma^2 }{4} \left( 
\frac{\partial^2 F_{[2], \mbox{\scriptsize eq}}}{\partial \theta_1^2} + \frac{\partial^2 F_{[2], \mbox{\scriptsize eq}}}{\partial \theta_2^2} \right) + F_{[2], \mbox{\scriptsize eq}} =   \delta(\theta_1 - \theta_2) .
\end{eqnarray*}
Taking Fourier series, it is easy to see that the Fourier coefficients $\hat F_{[2], \mbox{\scriptsize eq}} (n_1,n_2)$ such that $n_1 + n_2 \not = 0$ are identically zero. Therefore, $F_{[2], \mbox{\scriptsize eq}}$ is the sum of a Fourier series which only involves wave-numbers $(n_1,n_2)$ such that $n_1 + n_2 = 0$. It follows that $F_{[2], \mbox{\scriptsize eq}}$ is of the form (\ref{eq_equi_F2}) with ${\mathcal M} (\theta)$ satisfying (\ref{eq_M_Helmholtz}). Formula (\ref{eq_express_M}) is a simple calculation which is left to the reader. \endproof

\begin{remark}
We again verify that the stationary solution of the CL hierarchy does not satisfy the chaos property. Indeed, we have 
$$ F_{[2], \mbox{\scriptsize eq}} (\theta_1, \theta_2) = {\mathcal M}(\theta_1 - \theta_2) \not = 1 = F_{[1], \mbox{\scriptsize eq}} (\theta_1) \, F_{[1], \mbox{\scriptsize eq}} (\theta_2). $$
That the equilibrium distribution is not chaotic was already proved in \cite{CDW}. Here, we have shown that the whole time-dynamics of the hierarchy does not satisfy the chaos property. Again, this is not in contradiction to the  propagation of chaos result of \cite{CDW}, which applies to the unscaled dynamics at a shorter time-scale. 
\label{rem_chaos_equil_CLD}
\end{remark}

\begin{remark}
Because of (\ref{eq_scaling_CLD}), we can denote the scaled noise probability by $g_N$. In \cite{CDW}, the scaling of $g_N$ was defined in such a way that 
\begin{equation}
\lim_{N\to\infty} (N-2)(\widehat g_N(n)-1) := G(n)
\label{eq_def_G}
\end{equation}
exists and is finite, 
where $\widehat g_N(n)$ is the $n$-th Fourier coefficient of $g_N$: 
$$ \widehat g_N(n) = \int_{[0,2\pi]} e^{-ik\theta} \, g_N(\theta) \, \frac{d \theta}{2 \pi}, \quad n \in {\mathbb Z}. $$
Here, assuming that $\gamma$ defined by (\ref{eq_def_gamma}) decays at infinity fast enough, we find that $G(n)$ resulting from (\ref{eq_def_G}) is finite and given by
\begin{equation}
G(n) = - \frac{\sigma^2 \, n^2}{2}. 
\label{eq_G_here}
\end{equation}
In \cite{CDW}, the correlation function ${\mathcal M}$ is found from its Fourier series:
$$ \widehat {\mathcal M} (n) = \frac{1}{1 - G(n)} , $$
which, with (\ref{eq_G_here}), yields 
\begin{equation}
\widehat {\mathcal M} (n) = \frac{1}{1 + \bar \sigma^2 \, n^2}, \qquad \bar \sigma = \frac{\sigma}{\sqrt 2} , 
\label{eq_hatM_here}
\end{equation}
for our choice of the noise distribution $g_N$. It is a simple computation to check that the Fourier series of (\ref{eq_express_M}) precisely gives (\ref{eq_hatM_here}), showing the consistency between the present results and those of \cite{CDW}. 
\label{rem_compar_CDW}
\end{remark}

\medskip
\noindent
We can say more about the marginals of the invariant density, namely:

\begin{proposition}
The $k$-particle marginal of the equilibrium $F_{[k], \mbox{\scriptsize eq}}$ is of the form 
\begin{eqnarray} 
&&\hspace{-1cm}
F_{[k], \mbox{\scriptsize eq}} (\theta_1, \ldots, \theta_k) = \phi_{k-1} (\theta_1-\theta_k, \ldots, \theta_{k-1} - \theta_k) 
\label{eq_Feq_1}
\\
&&\hspace{-0.5cm}
= \phi_{k-1} (\theta_1-\theta_i, \ldots, \theta_{i-1}-\theta_i, \theta_{i+1}-\theta_i, \ldots, \theta_k - \theta_i), \quad \forall i \in \{ 1, \dots, k \}, \label{eq_Feq_2}
\end{eqnarray}
where $ \phi_{k-1} (\xi_1, \ldots, \xi_{k-1})$ is $2 \pi$-periodic in each coordinate and permutationally invariant. 
Moreover, 
\begin{eqnarray} 
&&\hspace{-1cm}
\phi_{k-1} (\xi_1, \ldots, \xi_{k-1}) = \int_{\xi_k \in [0,2\pi]} \phi_{k-1} (\xi_1, \ldots, \xi_{k-1}, \xi_k) \, \frac{d \xi_k}{2 \pi}. 
\label{eq_Feq_int}
\end{eqnarray}
\label{prop_marginals_equilibrium}
\end{proposition}

\medskip
\noindent
Proposition \ref{prop_marginals_equilibrium} shows that the $k$-particle marginal is just a function of the relative phases of the particles with respectively to any one given particle. Therefore, up to the phase of this given particle, the equilibrium statistics is perfectly known. 

\medskip
\noindent
{\bf Proof.} We prove by induction that the Fourier coefficients $ \hat F_{[k], \mbox{\scriptsize eq}} (n_1, \ldots,  n_k)$ of $_{[k], \mbox{\scriptsize eq}}$ are of the form 
\begin{equation} 
\hat F_{[k], \mbox{\scriptsize eq}} (n_1, \ldots,  n_k) = c(n_1, \ldots, n_k) \, \delta(n_1 + \ldots + n_k), 
\label{eq_Fourier_phi}
\end{equation}
with permutationally invariant coefficients $c(n_1, \ldots, n_k)$. The notation $\delta(n)$ stands for the Kronecker symbol $\delta_{n \, 0}$. Indeed, if this is the case, we can write 
\begin{eqnarray} 
&&\hspace{-1cm}
F_{[k], \mbox{\scriptsize eq}} (\theta_1, \ldots, \theta_k)  = \sum_{(n_1, \ldots, n_k) \, \in \,  {\mathbb Z}^k \atop \scriptstyle n_1 + \ldots + n_k = \, 0 }
c(n_1, \ldots, n_k) \, e^{i (n_1 \theta_1 + \ldots + n_k \theta_k)} \nonumber \\
&&\hspace{-1.cm}
= \hspace{-0.3cm} \sum_{(n_1, \ldots, n_{k-1}) \, \in \,  {\mathbb Z}^{k-1} }
c(n_1, \ldots, n_{k-1}, - (n_1 + \ldots + n_{k-1})) \, e^{i (n_1 (\theta_1-\theta_k) + \ldots + n_{k-1} (\theta_{k-1} - \theta_k))} 
\label{eq_def_phi_1} \\
&&\hspace{-1cm}
= \hspace{-0.3cm} \sum_{(n_1, \ldots, \hat n_i, \ldots n_k) \, \in \,  {\mathbb Z}^{k-1} }
c(n_1, \ldots, n_{i-1},  - (n_1 + \ldots + \hat n_i + \ldots + n_k), n_{i+1}, \ldots, n_k) \nonumber \\
&&\hspace{3.5cm}
e^{i (n_1 (\theta_1-\theta_i) + \ldots + n_{i-1} (\theta_{i-1} - \theta_i) + n_{i+1} (\theta_{i+1} - \theta_i) + \ldots + n_k (\theta_k - \theta_i))} \label{eq_def_phi_2} .
\end{eqnarray}
Thanks to (\ref{eq_def_phi_1}),  we can write (\ref{eq_Feq_1}) and thanks to (\ref{eq_def_phi_2}),  we have (\ref{eq_Feq_2}), with 
\begin{eqnarray} 
&&\hspace{-1cm}
\phi_{k-1} (\xi_1, \ldots, \xi_{k-1})  = \nonumber \\
&&\hspace{-0.cm}
= \hspace{-0.3cm} \sum_{(n_1, \ldots, n_{k-1}) \, \in \,  {\mathbb Z}^{k-1} }
c(n_1, \ldots, n_{k-1}, - (n_1 + \ldots + n_{k-1})) \, e^{i (n_1 \xi_1 + \ldots + n_{k-1} \xi_{k-1})} 
\label{eq_phi_3} \\
&&\hspace{-0.cm}
= \hspace{-0.3cm} \sum_{(n_1, \ldots, \hat n_i, \ldots n_k) \, \in \,  {\mathbb Z}^{k-1} }
c(n_1, \ldots, n_{i-1},  - (n_1 + \ldots + \hat n_i + \ldots + n_k), n_{i+1}, \ldots, n_k) \nonumber \\
&&\hspace{6.5cm}
e^{i (n_1 \xi + \ldots + n_{i-1} \xi_{i-1} + n_{i+1} \xi_i + \ldots + n_k \xi_{k-1})} .
\label{eq_phi_4}
\end{eqnarray}
We note that, thanks to the permutation symmetry of $c$, the two definitions (\ref{eq_phi_3}) and (\ref{eq_phi_4}) are consistent. Additionally, the so-defined $\phi_{k-1}$ is $2 \pi$-periodic and permutationally invariant. Finally, (\ref{eq_Feq_int}) follows from integrating (\ref{eq_Feq_2}) with respect to $\theta_k$. 

We now prove (\ref{eq_Fourier_phi}). Taking the Fourier series of (\ref{eq_BBGKY_CLD_grazing}) (with $\frac{\partial F_{[k]}}{\partial t} = 0$ since we are interested in the equilibrium, we find:
\begin{eqnarray} 
& & \hspace{-1cm}  \hat F_{[k], \mbox{\scriptsize eq}} (n_1, \ldots, n_k)  =  \frac{1}{\frac{\sigma^2 }{2} \big( \sum_{i \leq k } n_i^2 \big) + k(k-1)} \sum_{i<j\leq k} 
\nonumber \\
& & \hspace{0cm}  
\big[  \hat F_{[k-1], \mbox{\scriptsize eq}} (n_1, \ldots, \hat n_i, \ldots, n_i + n_j, \ldots, n_k) + \nonumber \\
& & \hspace{5cm}  
+ \hat F_{[k-1], \mbox{\scriptsize eq}} (n_1, \ldots, n_i + n_j, \ldots, \hat n_j, \ldots, n_k) \big].
\label{eq_BBGKY_CLD_grazing_Fourier}
\end{eqnarray}
Now, we proceed by induction, starting from $k=2$. In this case, from (\ref{eq_equi_F1}), we have $\hat F_{[1], \mbox{\scriptsize eq}} (n) = \delta(n)$ and we get 
\begin{eqnarray*} 
& & \hspace{-1cm}  \hat F_{[2], \mbox{\scriptsize eq}} (n_1, n_2)  =  \frac{2}{\frac{\sigma^2 }{2} \big( n_1^2 + n_2^2 \big) + 2} \, \delta(n_1 + n_2), 
\end{eqnarray*}
which is of the form (\ref{eq_Fourier_phi}). Next, inserting (\ref{eq_Fourier_phi}) at level $k-1$ into (\ref{eq_BBGKY_CLD_grazing_Fourier}), we get: 
\begin{eqnarray*} 
& & \hspace{-1cm}  \hat F_{[k], \mbox{\scriptsize eq}} (n_1, \ldots, n_k)  =  \frac{1}{\frac{\sigma^2 }{2} \big( \sum_{i \leq k } n_i^2 \big) + k(k-1)} \sum_{i<j\leq k} 
\nonumber \\
& & \hspace{0cm}  
\big[  c (n_1, \ldots, \hat n_i, \ldots, n_i + n_j, \ldots, n_k) + \nonumber \\
& & \hspace{3cm}  
+ c (n_1, \ldots, n_i + n_j, \ldots, \hat n_j, \ldots, n_k) \big] \, \delta(n_1 + \ldots + n_k), 
\end{eqnarray*}
which is also of the form (\ref{eq_Fourier_phi}). Therefore, (\ref{eq_Fourier_phi}) is proven by induction on $k$, which ends the proof of proposition \ref{prop_marginals_equilibrium}. \endproof

\medskip
\noindent
The particularly simple form of the CL hierarchy leads to many analytical formulas. These formulas can be used to compare the theoretical predictions to numerical experiments. This is performed in the next section.

\setcounter{equation}{0}
\section{Numerical simulations}
\label{sec_numerics}

\subsection{Operation mode}
\label{subsec_simulations_principles}

The discrete CL and BDG dynamics (see section \ref{sec_particle_models}) have been simulated. The experimental protocol is as follows. For a given  number $N$ of particles, $M$ independent simulations are performed (in the experiments, $M=1000$). Each simulation is run until an equilibrium is reached. The detection of the equilibrium is detailed below. Once an equilibrium has been reached, one particle is randomly picked and its corresponding $\theta$ is collected and added to an histogram. Similarly, a pair of particles is randomly picked and the corresponding pair $(\theta_1,\theta_2)$ is added to a 2-dimensional histogram. Then, a new simulation is started and the procedure is continued until the $M$ simulations have been performed. The histograms of the $M$ samples of $\theta$ or $(\theta_1,\theta_2)$ gives experimental access to the steady-state one and two particle distribution functions of the process. 

To quantify if the equilibrium state is reached, three macroscopic quantities are observed:
\begin{itemize}
 \item[-] The \emph{average velocity}: $\bar v^n = \frac{1}{N} \sum_{i=1}^N v_i^n $
 \item[-] The \emph{average direction}: $\omega^n = \frac{\bar v^n}{\lvert \bar v^n\rvert}$
 \item[-] The \emph{order parameter}: $\alpha^n = \lvert \bar v^n\rvert^2$
\end{itemize}
It is assumed that the equilibrium is reached when the relative difference $\frac{\lvert \varphi^{n+\kappa} - \varphi^n \rvert}{\lvert \varphi^n \rvert}$ (where $\varphi$ is one of the above defined macroscopic quantities) at two iterates $n$ and $n+\kappa$ (where $\kappa $ is a constant (generally $\kappa = 1.2 N$)) is smaller than a fixed tolerance $\epsilon \ll 1$, i.e.  
$$\frac{\lvert \varphi^{n+\kappa} - \varphi^n \rvert}{\lvert \varphi^n \rvert}<\epsilon . $$ 
Therefore, the process is supposed to have reached an equilibrium if during $\kappa$ time steps the relative change in the macroscopic quantities is small. However a lower bound for the number of iterations $\mathcal{N}$ is set in order to have all the particles interact at least once:
\begin{equation}
 \mathcal{N} > \lambda N
\end{equation}
Typically the value $\lambda = 3$ has been chosen.

Finally, we use a periodized Gaussian with standard deviation $\sigma_n$ as noise distribution: 
\begin{equation}
g_N(\theta) = 2\pi \sum_{k=-\infty}^{+\infty} \frac{1}{\sigma_N \sqrt{2\pi}} e^\frac{- (\theta+2k\pi)^2}{2\sigma_N^2}, \quad \theta \in [-\pi, \pi]. 
\label{E-Normal}
\end{equation} 
We use 
\begin{equation}
\sigma_N = \frac{2 \pi \gamma}{\sqrt N} , 
\label{eq_def_sigmaN}
\end{equation} 
to be consistent with the rescaling of proposition \ref{thm_hierarchy_grazing_CLD}. Specifically, we have 
\begin{equation}
\sigma = 2 \pi \gamma,
\label{eq_rel_sigma_gamma}
\end{equation} 
where $\sigma$ is the standard deviation defined at (\ref{eq_def_sigma_tau}). $\gamma$ is a measure of the intensity of the noise if the length of the unit circle is fixed to one. We will use $\gamma = 1/2$, $\gamma = 1/20$, $\gamma = 1/200$, meaning that we expect that the width of the noise distribution (which is equal to $2 \gamma$) be comparable to the length of the unit circle or equal to $10\%$ or $1\%$ of it respectively.   

For the BDG dynamics, we use a similar law (\ref{E-Normal}) for $h_N$ with $\tau_N$ related to $\gamma'$ by a similar relation as (\ref{eq_def_sigmaN}) and $\gamma'$ to $\tau$ defined at (\ref{eq_def_sigma_tau}) by a similar relation as (\ref{eq_rel_sigma_gamma}).

\subsection{Results}
\label{subsec_simulations_results}

\subsubsection{Unbiased BDG dynamics}
\label{subsec_simulations_results_BDG}

\begin{figure}[h]
\centering
\subfigure[$N=10^2$]{\includegraphics[width=0.5\textwidth]{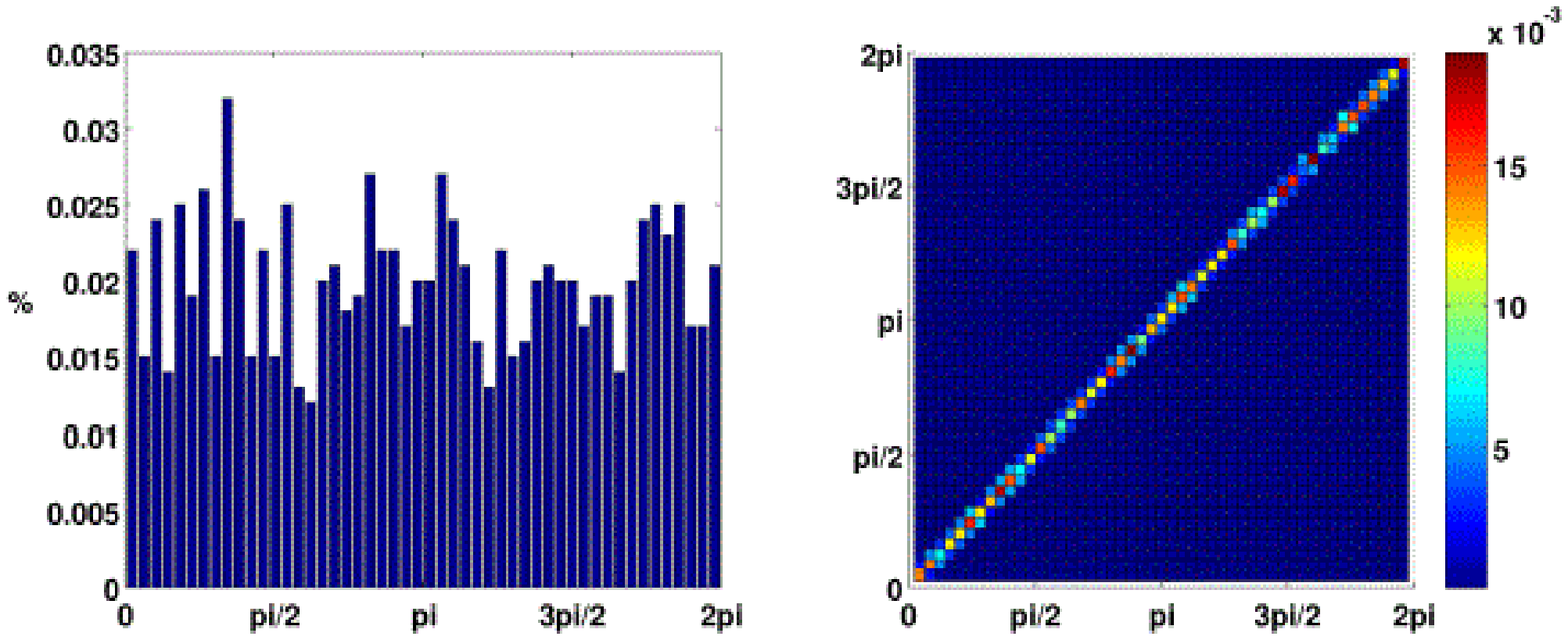}}
\subfigure[$N=10^3$]{\includegraphics[width=0.5\textwidth]{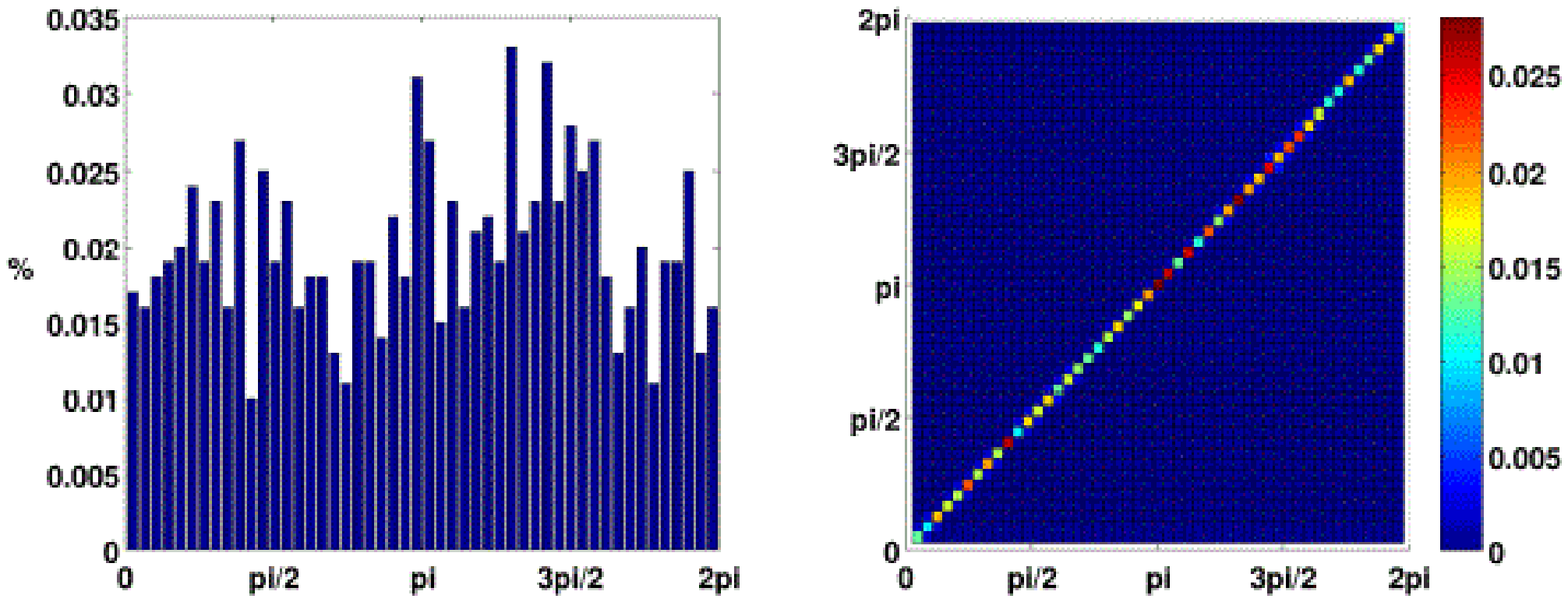}}
\subfigure[$N=10^4$]{\includegraphics[width=0.5\textwidth]{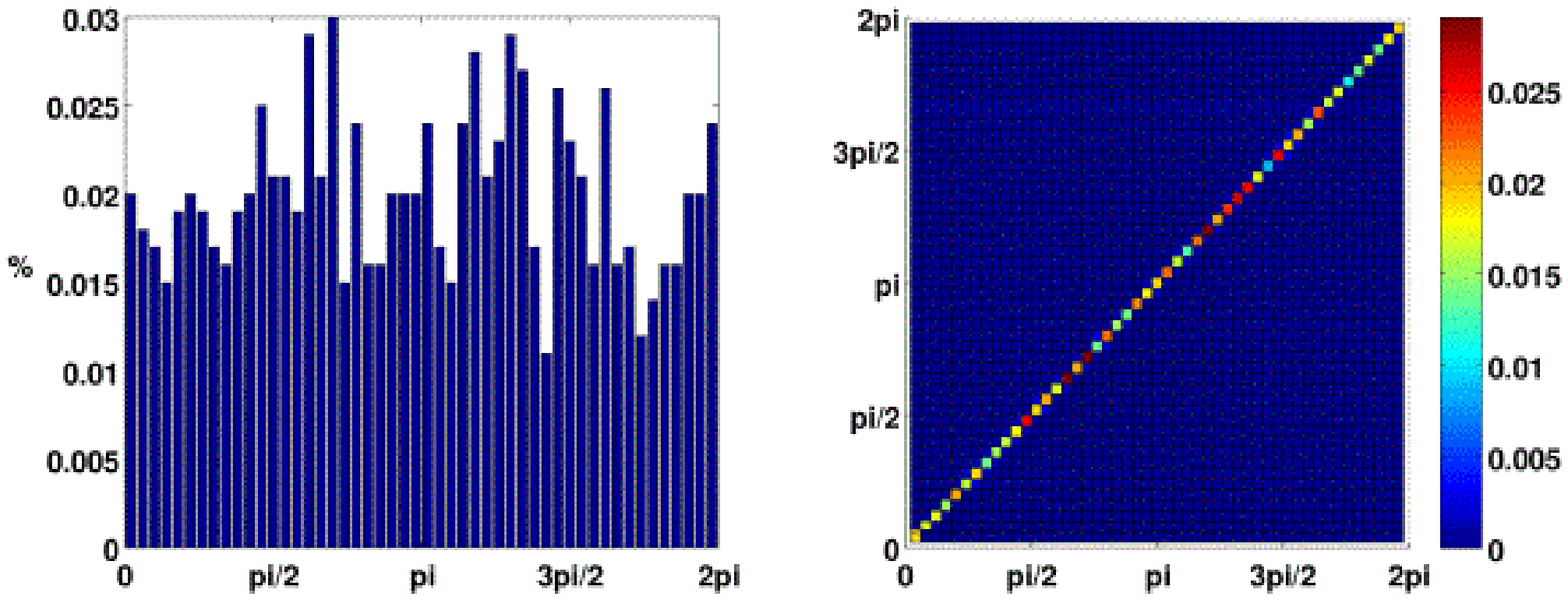}}
\subfigure[$N=10^5$]{\includegraphics[width=0.5\textwidth]{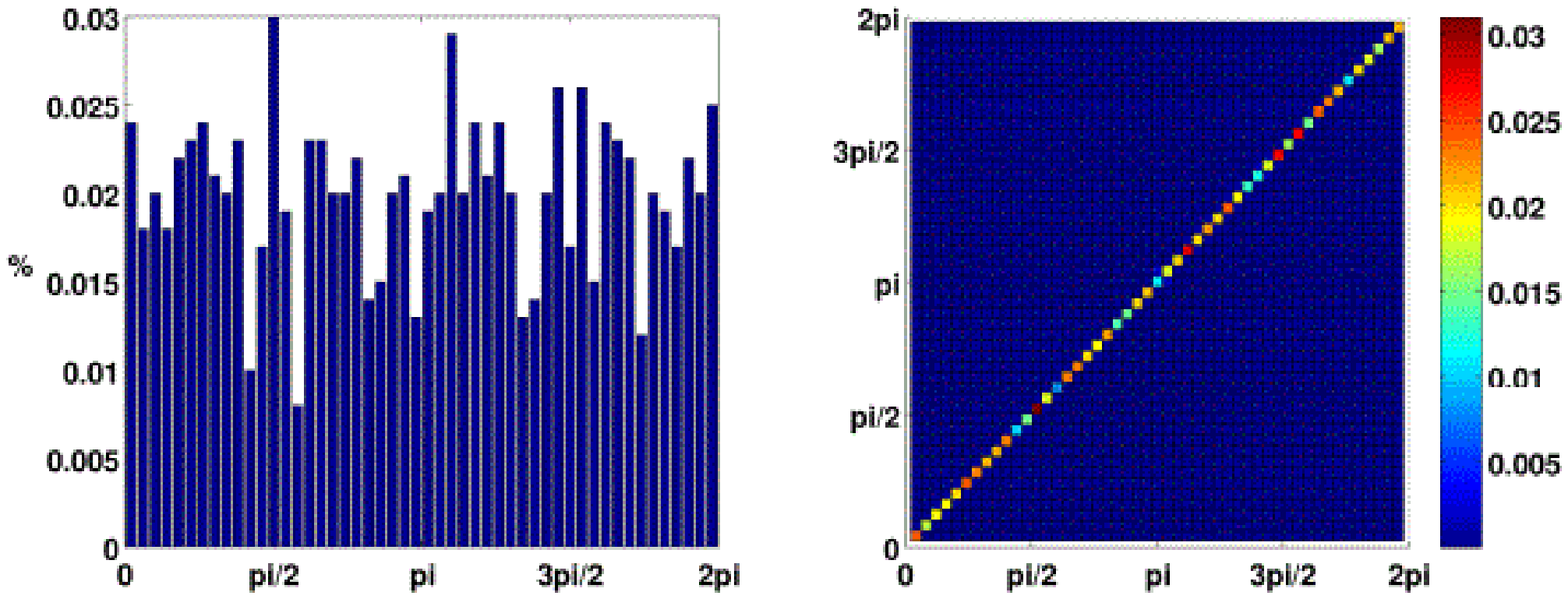}}
\caption{Unbiased BDG dynamics with $\gamma = 0.02$ and four values of the particle number (from top to bottom: $N=10^2$, $N=10^3$, $N=10^4$, $N=10^5$). One and two particle equilibrium distributions (left and right respectively). Left: the horizontal axis is $\theta \in [0,2\pi]$, the vertical axis is the distribution of $\theta$. Right: the square is the domain of $(\theta_1,\theta_2) \in [0,2\pi]^2$. The color coding indicates the 2-particle distribution. }
\label{Correlations-BDG-02}
\end{figure}

The equilibrium one and two particle distributions for the unbiased BDG dynamics are presented on Fig. \ref{Correlations-BDG-02} as a function of the number of particles for $\gamma = 0.02$. The histogram of the one-particle distribution (left) clearly shows that the distribution is isotropic. Indeed, in spite of a relatively large noise level, we cannot distinguish any structure. The 2-dimensional histogram of the two-particle distribution (right) shows a high level of correlation along the diagonal, due to the small value of $\gamma$. Clearly, the 2-particle distribution is not a product of two copies of the one-particle distribution and the chaos property is not satisfied.

\subsubsection{CL dynamics}
\label{subsec_simulations_results_CLD}

\begin{figure}[h]
\centering
\subfigure[$N=10^2$]{\includegraphics[width=0.5\textwidth]{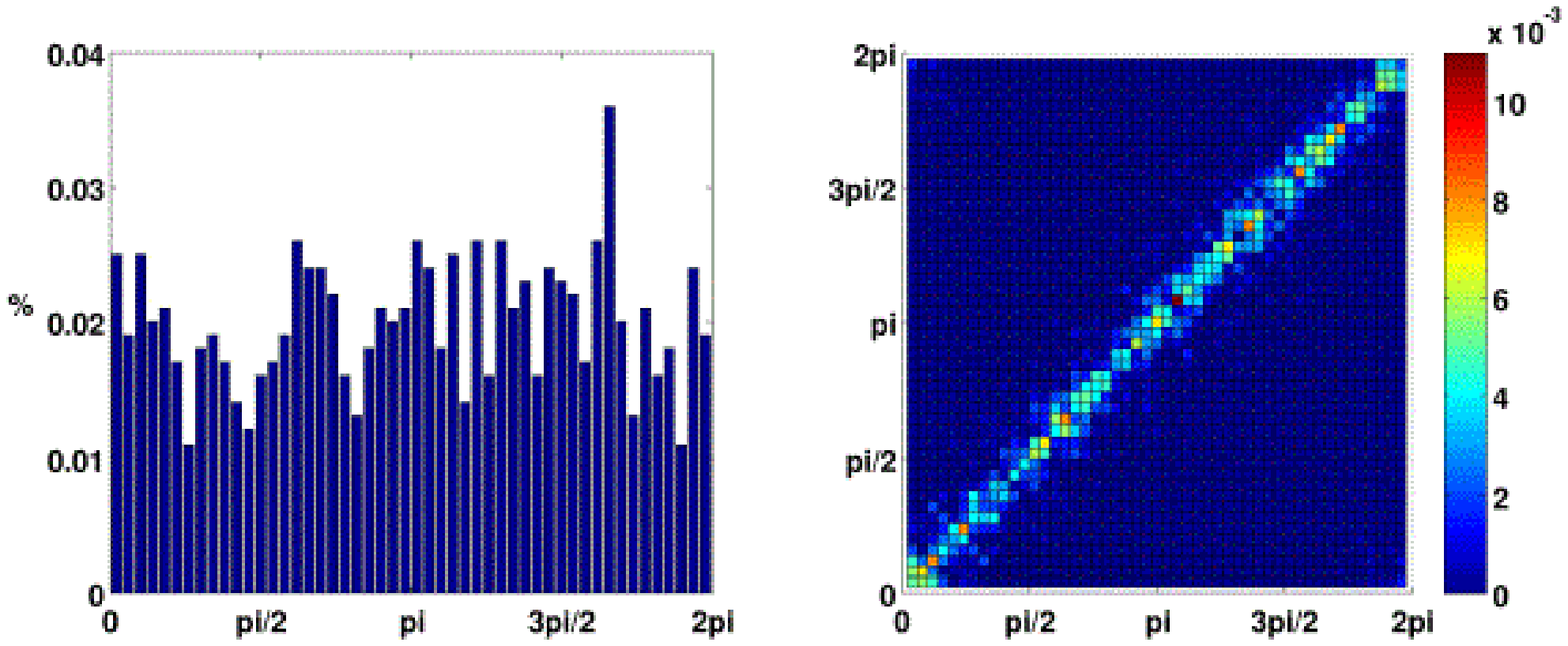}}
\subfigure[$N=10^3$]{\includegraphics[width=0.5\textwidth]{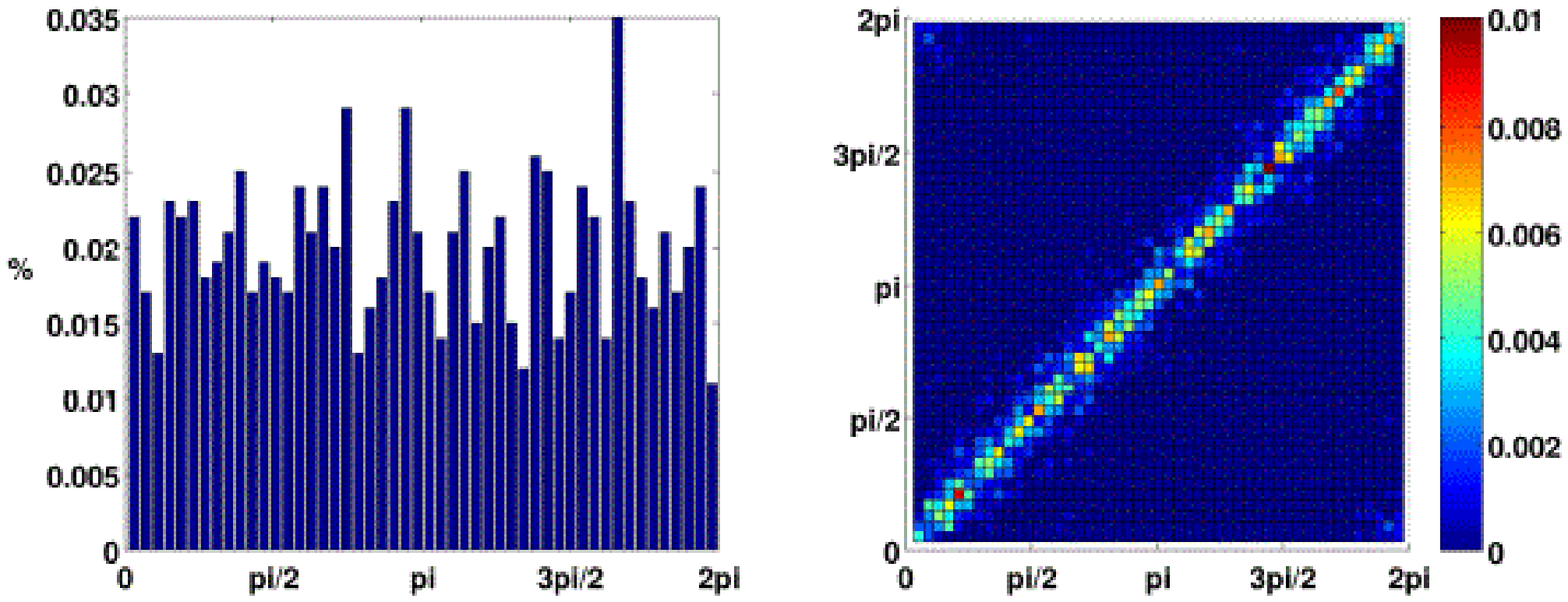}}
\subfigure[$N=10^4$]{\includegraphics[width=0.5\textwidth]{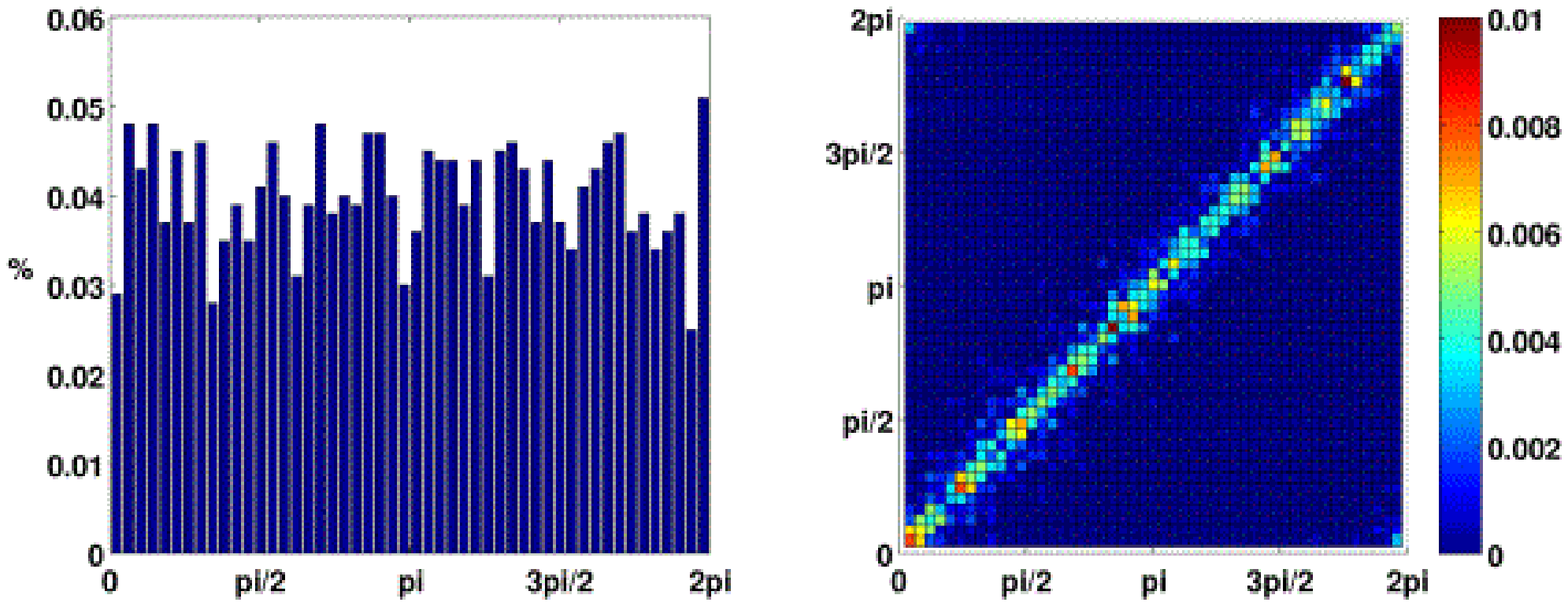}}
\subfigure[$N=10^5$]{\includegraphics[width=0.5\textwidth]{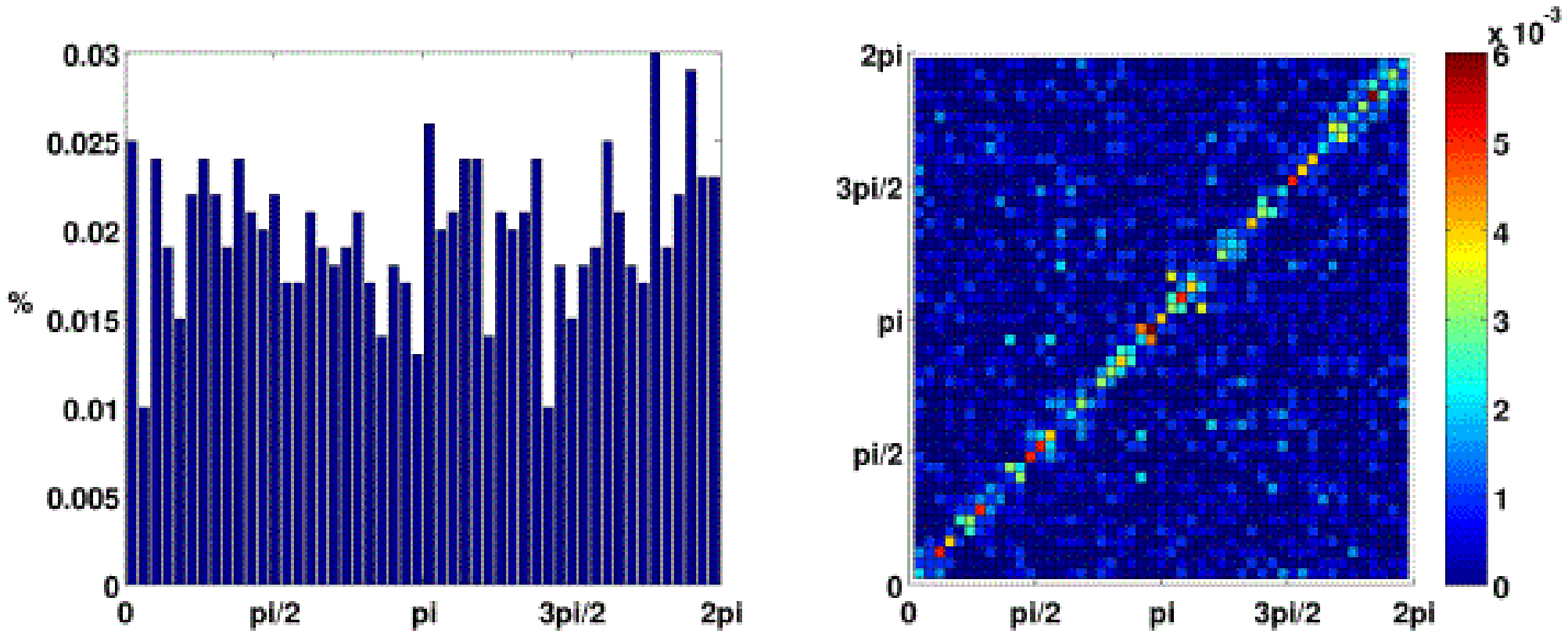}}
\caption{CL dynamics with $\gamma = 0.02$ and four values of the particle number (from top to bottom: $N=10^2$, $N=10^3$, $N=10^4$, $N=10^5$). One and two particle equilibrium distributions (left and right respectively). Left: the horizontal axis is $\theta \in [0,2\pi]$, the vertical axis is the distribution of $\theta$. Right: the square is the domain of $(\theta_1,\theta_2) \in [0,2\pi]^2$. The color coding indicates the 2-particle distribution. }
\label{Correlations-Carlen-02}
\end{figure}

The equilibrium one and two particle distributions are presented on Fig. \ref{Correlations-Carlen-02} as a function of the number of particles. Similar conclusions as for the unbiased BDG dynamics can be drawn, except for a slightly larger spread of the 2-particle distribution function about the diagonal. We also see that the spread is about the same for all particle numbers, which confirms the relevance of the choice of  scaling the noise standard deviation like the inverse square root of the particle number.  The statistics of the 2-particle distribution function for $N=10^5$ is noisy because the convergence to the equilibrium is very slow and the simulation was stopped before the equilibrium was reached. Yet, correlations appear with  a spread of the two-particle distribution about the diagonal of the same order of magnitude as for smaller number of particles.

\begin{figure}[h]
\centering
\subfigure[$\gamma = \frac{1}{2}$]{\includegraphics[width=0.5\textwidth]{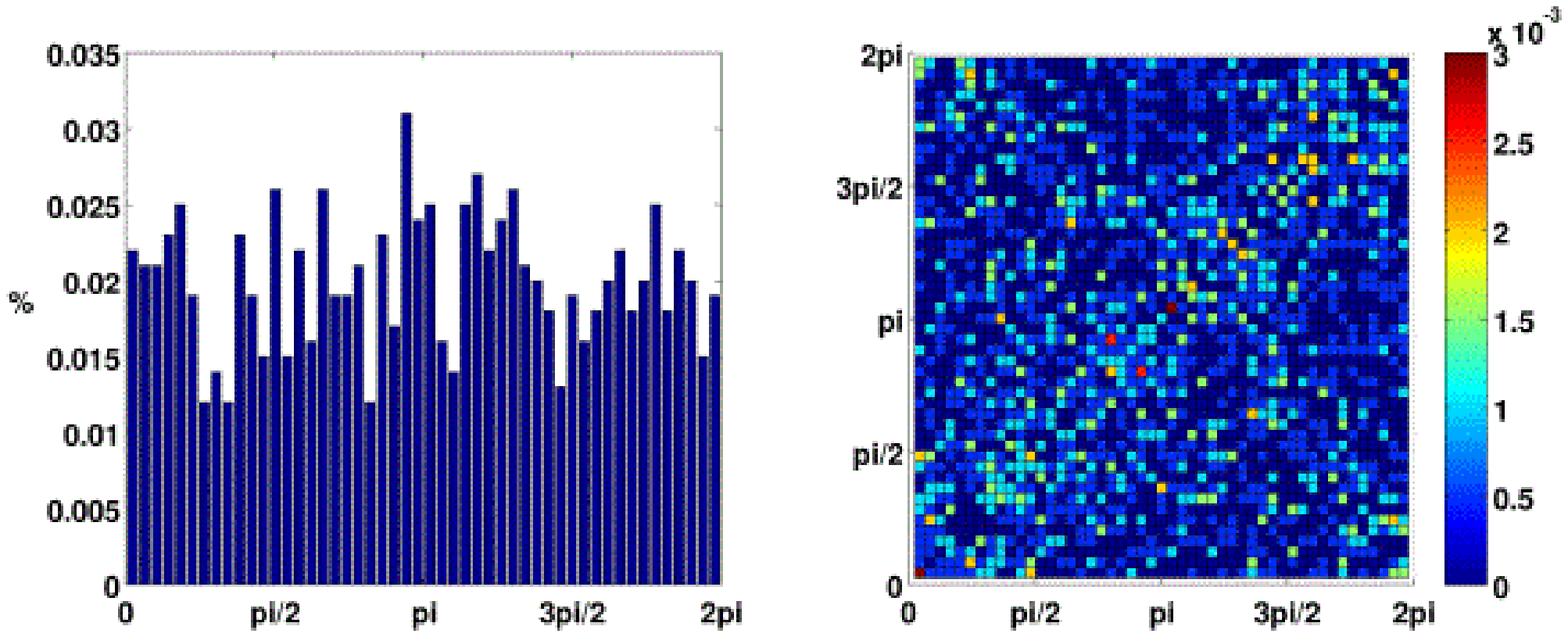}}
\subfigure[$\gamma = \frac{1}{20}$]{\includegraphics[width=0.5\textwidth]{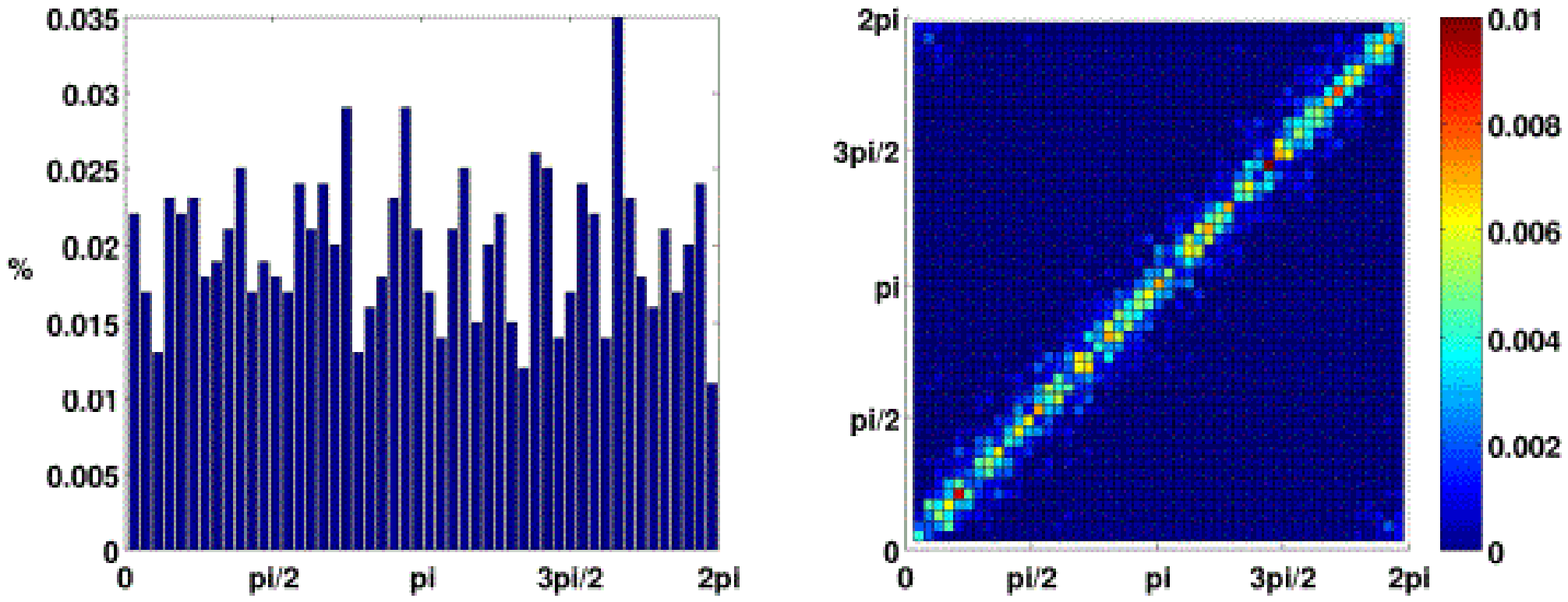}}
\subfigure[$\gamma = \frac{1}{200}$]{\includegraphics[width=0.5\textwidth]{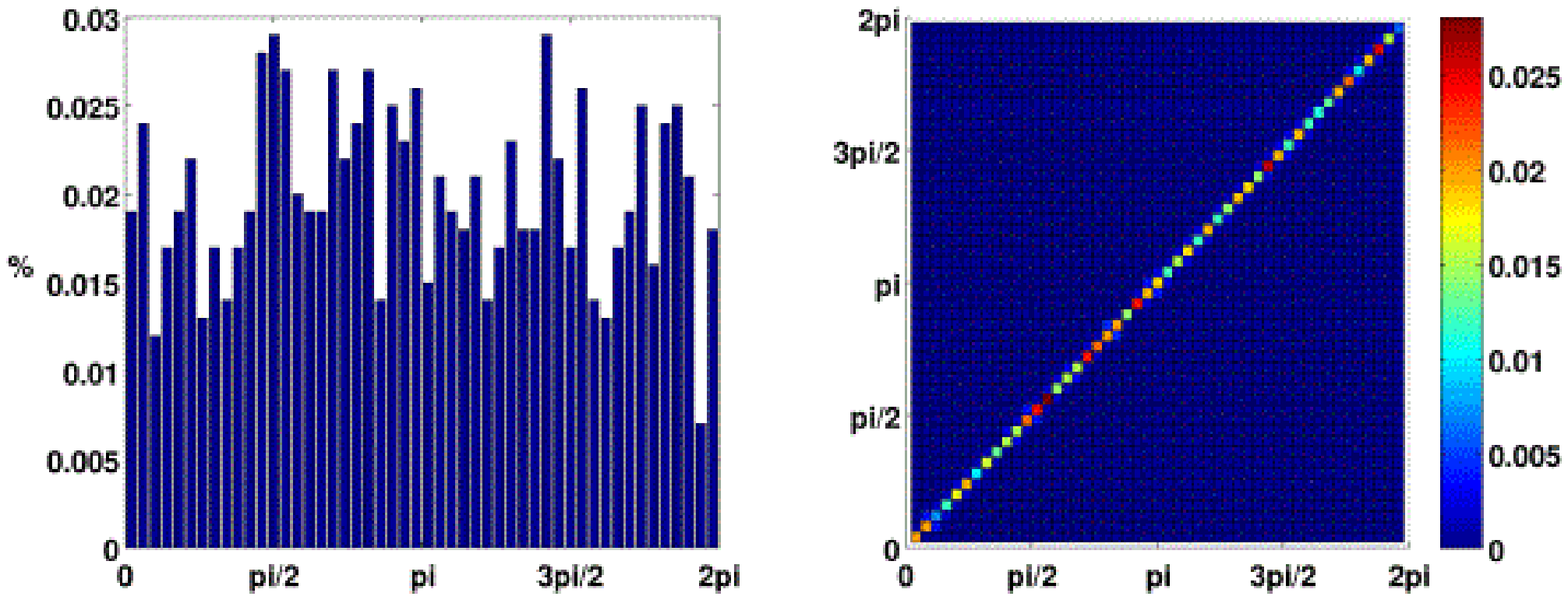}}
\caption{CL dynamics with $N=10^3$ and three values of the parameter $\gamma$ (from top to bottom: $\gamma = 1/2$, $\gamma = 1/20$, $\gamma = 1/200$). One and two particle equilibrium distributions (left and right respectively). Left: the horizontal axis is $\theta \in [0,2\pi]$, the vertical axis is the distribution of $\theta$. Right: the square is the domain of $(\theta_1,\theta_2) \in [0,2\pi]^2$. The color coding indicates the 2-particle distribution. }
\label{Correlations-Carlen-en-fonction-gamma}
\end{figure}

Fig. \ref{Correlations-Carlen-en-fonction-gamma} explores the dependence of the one and two-particle distribution functions as a function of $\gamma$, for a given particle number $N=10^3$. For a large value of $\gamma$ ($\gamma = 1/2$), almost no structure appears in the two-particle distribution function. When $\gamma$ is reduced to $\gamma = 1/20$, a clear correlation along the diagonal builds up. When $\gamma$ is further reduced to $\gamma = 1/200$, the two-particle distribution exhibits a delta-like behavior, with a total concentration on the diagonal.

\subsubsection{BDG dynamics}
\label{subsec_simulations_results_BBDG}

\begin{figure}[h]
\centering
\subfigure[$N=10^2$]{\includegraphics[width=0.33\textwidth]{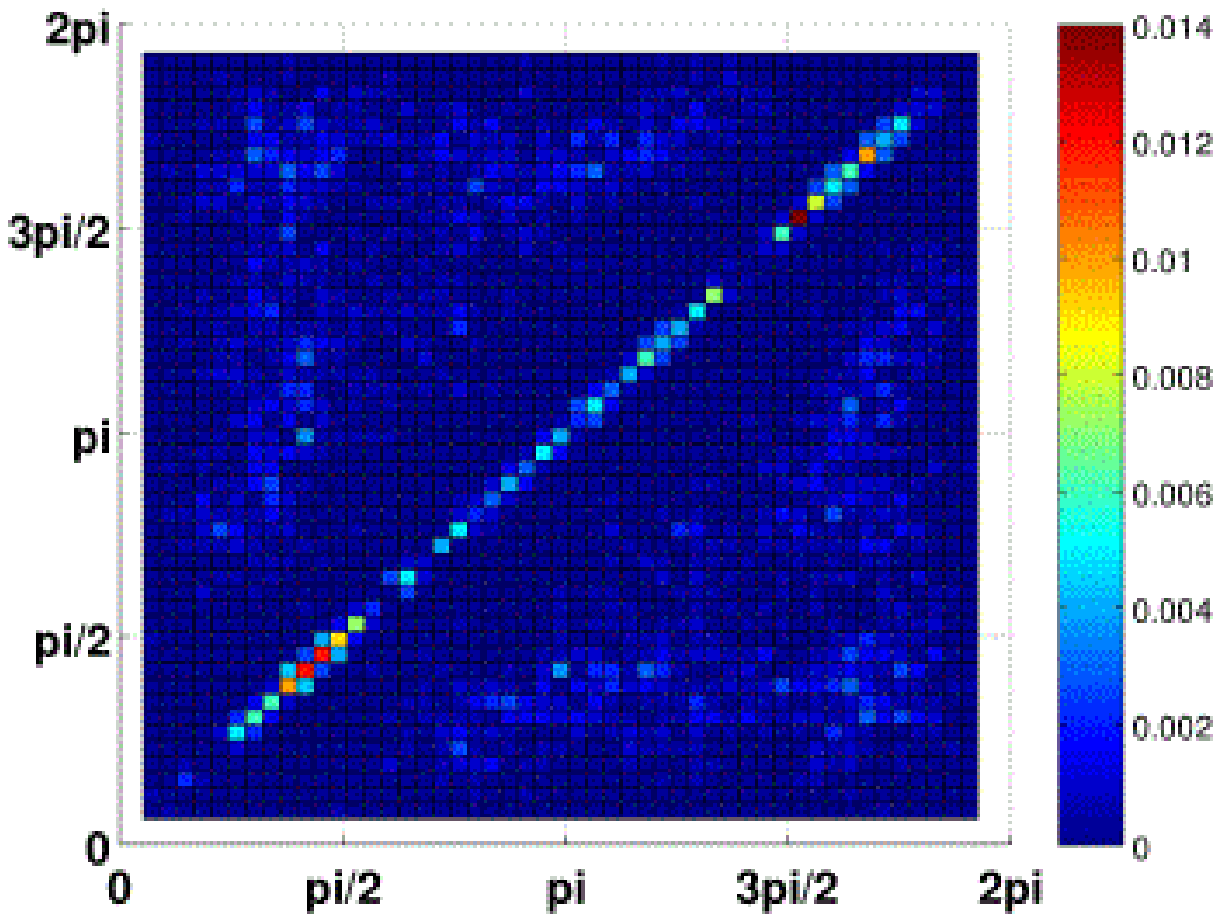}}
\subfigure[$N=10^3$]{\includegraphics[width=0.33\textwidth]{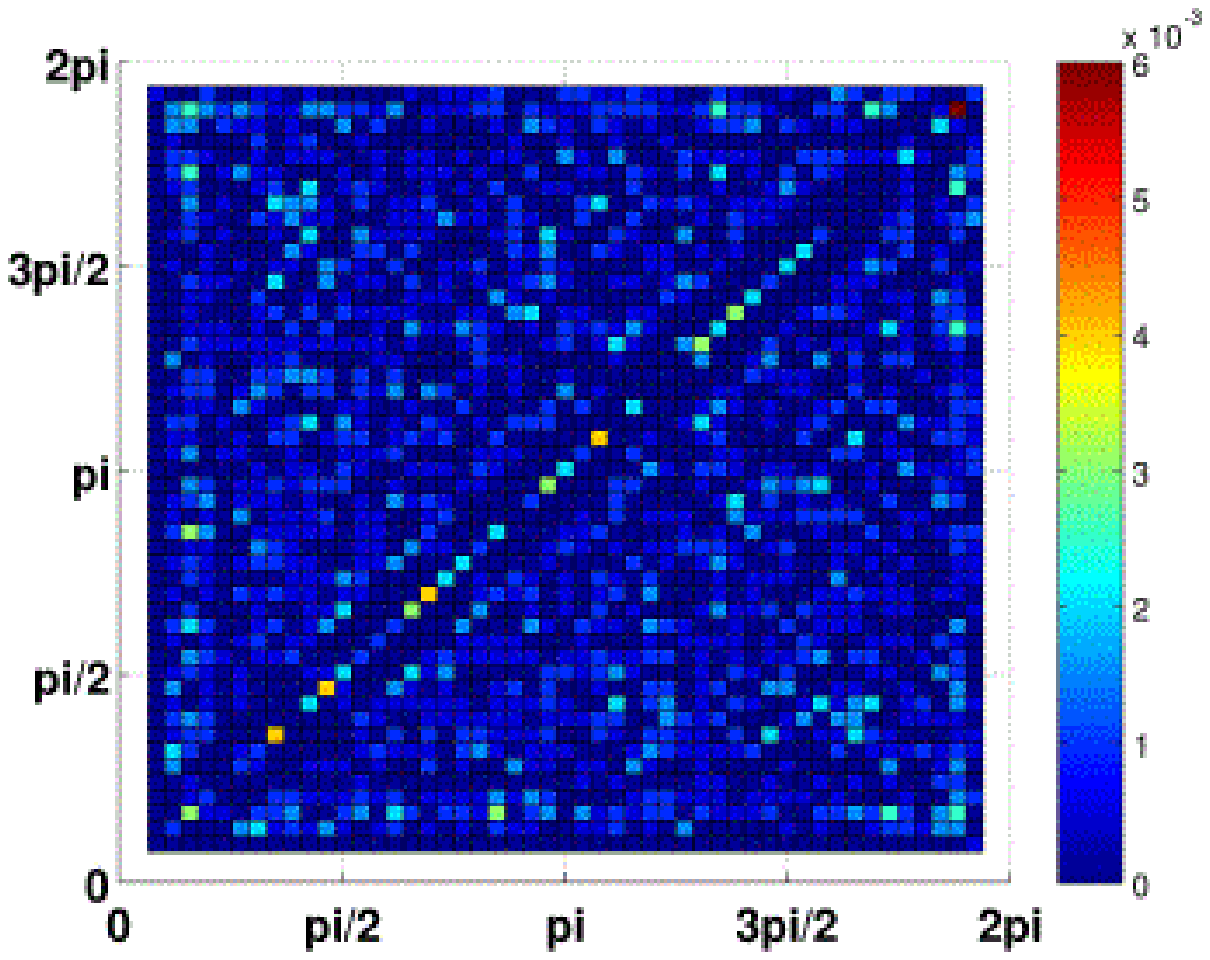}}
\subfigure[$N=10^4$]{\includegraphics[width=0.33\textwidth]{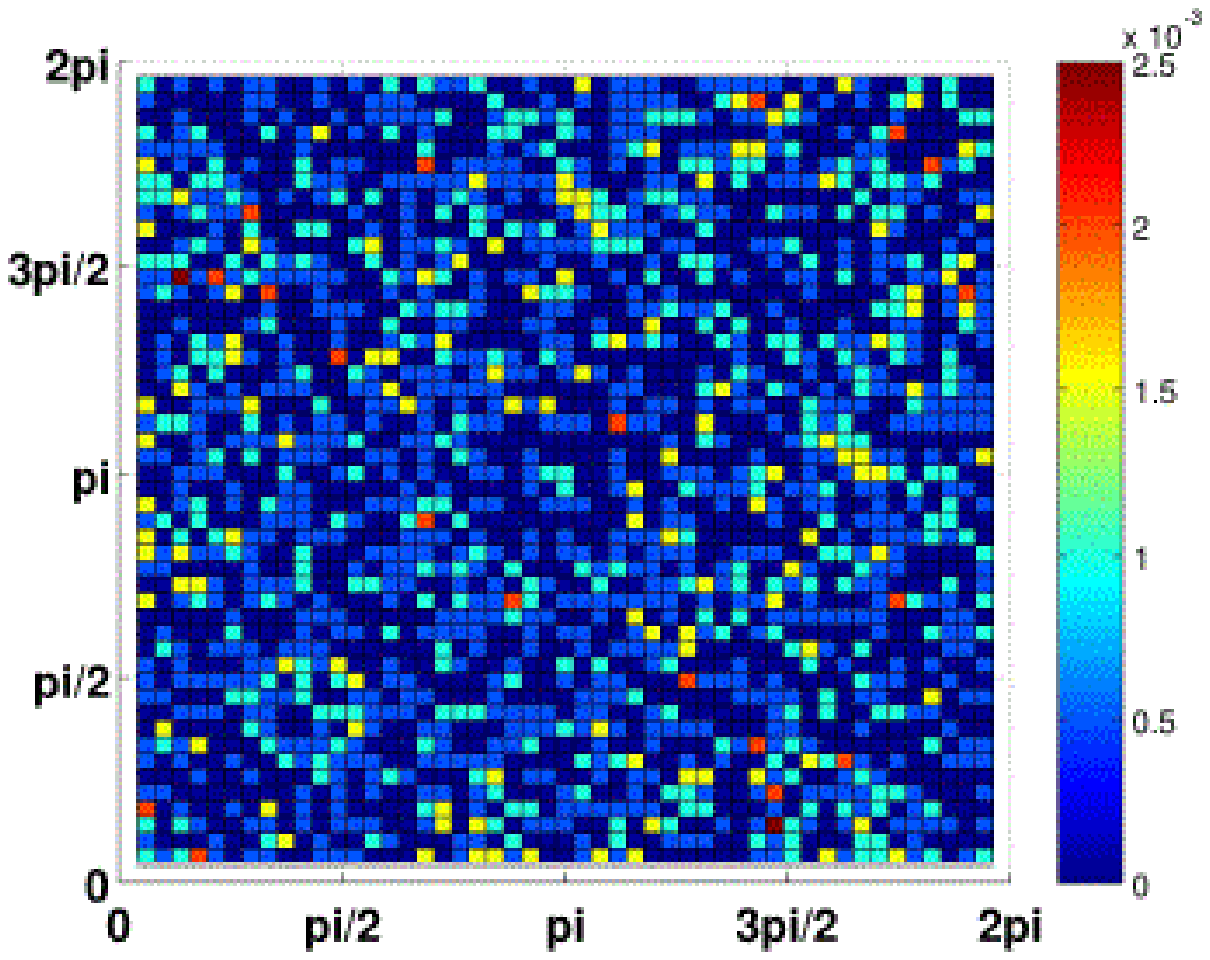}}
\caption{BDG dynamics with $\gamma = 0.02$ and $\gamma' = 10 \gamma$ and three values of the particle number (from top to bottom: $N=10^2$, $N=10^3$, $N=10^4$). The two particle equilibrium distribution only is plotted. The square is the domain of $(\theta_1,\theta_2) \in [0,2\pi]^2$. The color coding indicates the 2-particle distribution. }
\label{Correlations-BBDG}
\end{figure}

On figure~\ref{Correlations-BBDG}, the two-particle distribution law of the BDG dynamics is presented for $\frac{\gamma'}{\gamma} = 10$. This means that on the average, the particles are $10$ times closer after performing a collision than before. We are in the case of backwards diffusion, with $\sigma^2 - \tau^2 <0$. The one-particle distribution is not plotted as it looks very similar to the other cases. For a small particle number $N=10^2$, two-particle correlations along the diagonal are formed. Around the diagonal, a hollow zone reflects the fact that if particles are close they collide a lot and tend to correlate their velocities. On the other hand, if they are too far away, they cannot collide and correlation does not appear. This hollow zone seems a good indicator to measure the ratio $\frac{\gamma'}{\gamma}$, if it is not known. For a larger number of particles $N=10^3$ or $N=10^4$, the results are more noisy, but this is due to the fact that we have to stop the simulation before an equilibrium is reached, due to very large CPU times. Indeed, most of the collisions are rejected in the acceptance-rejection procedure. The acceptance-rejection procedure should be biased in order to restore reasonable simulation times. This has not been done yet and will be the subject of future works.

\subsection{Comparing the numerical and analytical solutions}

\begin{figure}[h]
\subfigure[]{\includegraphics[width=0.9\textwidth]{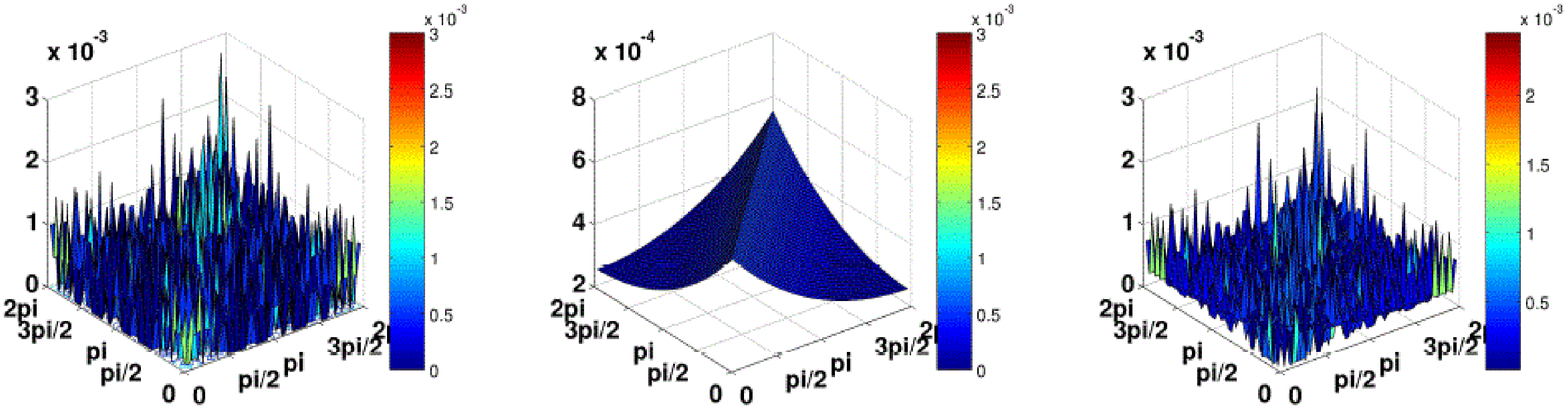}}
\subfigure[]{\includegraphics[width=0.9\textwidth]{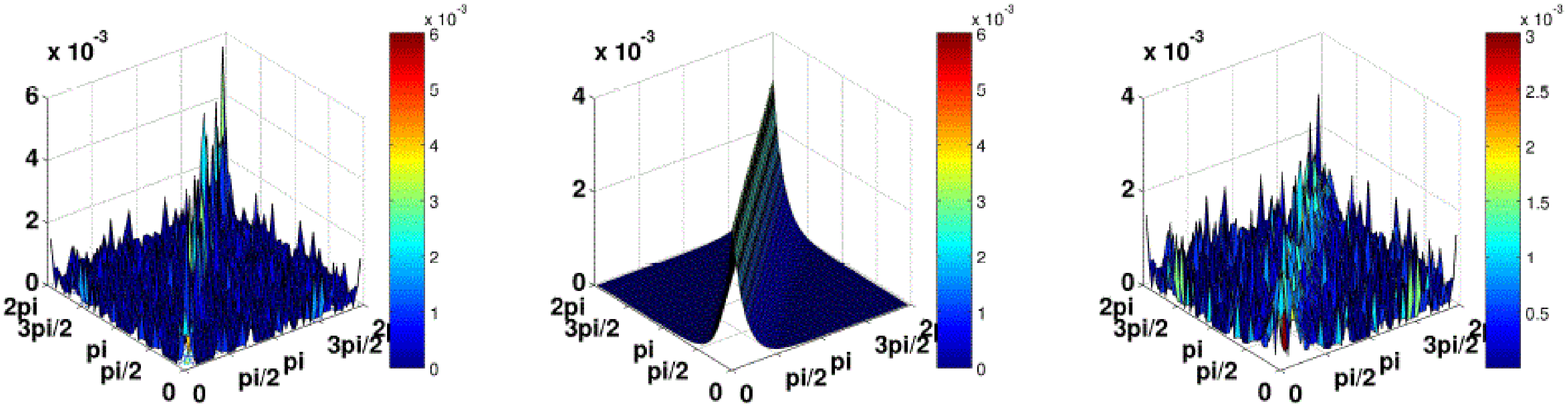}}
\subfigure[]{\includegraphics[width=0.9\textwidth]{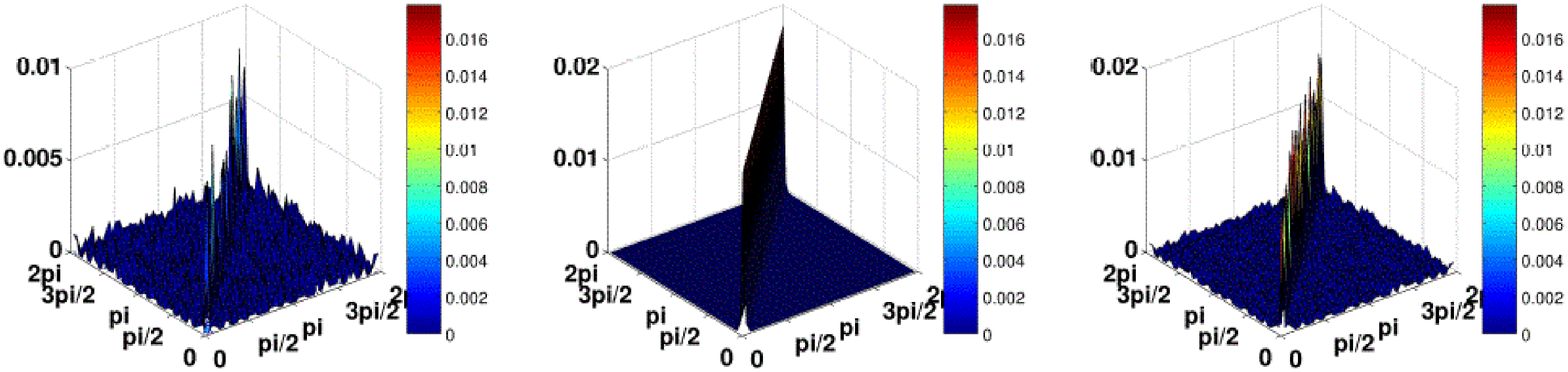}}
\caption{Perspective view of the equilibrium two-particle distribution of the CL dynamics as functions of $(\theta_1,\theta_2) \in [0,2\pi]^2$, with (a) $\gamma=\frac{1}{2}$ ; (b)  $\gamma=\frac{1}{20}$ ; (c) $\gamma=\frac{1}{200}$. Left: numerical solution. Middle: analytical solution. Right: Difference between the numerical and analytical solution.  }
\label{m2-cld-nuw-vs-an}
\end{figure}

Using (\ref{eq_express_M}) the analytical solution for the equilibrium two-particle distribution function of the CL dynamics can be visualized and compared to the numerical solution. Fig. \ref{m2-cld-nuw-vs-an} shows this comparison for three values of $\gamma$: $\gamma = 1/2$, $\gamma = 1/20$, $\gamma = 1/200$. Despite the inherent fluctuations due to the noise and the relatively limited number of independent experiments ($M=1000$), the agreement appears quite convincing.

\setcounter{equation}{0}
\section{Conclusion}
\label{sec_conclusion}

We have considers two models of biological swarm behavior, the BDG and CL models. They describe consensus formation in a collection of moving agents about the direction of motion to follow. We have first formed the master equations and BBGKY hierarchies of these two processes. Then, we have investigated their large particle number limit at the large time scale. To this aim, we had to simultaneously rescale the noise variance of the processes and make it small in order to keep the leading order of the dynamics finite. We have examined the resulting kinetic hierarchies and shown that, in the case of the CL dynamics,  propagation of chaos does not hold, neither at the stationary state, nor in the time-dynamical process. We could exhibit the general form of the $k$-particle correlations at the stationary state. The study leaves the related questions for the BDG dynamics mostly open. However, numerical simulations indicate that the BDG dynamics has a similar behavior as the CL dynamics. 

More theoretical investigations are on the way to understand the behavior of the BDG dynamics closely. In the case of the CL model, the knowledge of the general form of the correlations opens the way for the derivation of a generalized kinetic model. Such a model would replace the standard kinetic model which does not hold because of the violation of the chaos property. It would involve a kinetic-like equation for the correlations.

\clearpage
\setcounter{equation}{0}
\section{Appendix A. Establishment of master equations and hierarchies}
\label{sec_appA_proof_master}

\subsection{BDG master equation and hierarchy}
\label{sec_appA_proof_master_BBDG}

{\bf Proof of Proposition \ref{prop_mas_eq_BBDG}.} The Markov transition operator for the BDG dynamics is clearly given by:
\begin{eqnarray*} 
& & \hspace{-0cm}  Q_N \Phi =  \frac{2}{N(N-1)} \sum_{i<j} \left( \phantom{\int_{(w,z) \in {\mathbb T}^2}} \right. \\
& & \hspace{-0.5cm} H(\hat v_{ij}^* v_i) \, \int_{(w,z) \in {\mathbb T}^2}  \Phi (v_1, \ldots , w \hat v_{ij}, \ldots ,  z \hat v_{ij}, \ldots , v_N) \, g(w) \, g(z) \, dw \, dz  \\
& & \hspace{4cm}  \left. \phantom{\int_{(w,z) \in {\mathbb T}^2}  \Phi } +  (1-H(\hat v_{ij}^* v_i)) \, \Phi (v_1, \ldots , v_i , \ldots ,  v_j , \ldots , v_N)  \right) ,
\end{eqnarray*}
where $w \hat v_{ij}$ and $z \hat v_{ij}$ are on the $i$-th and $j$-th positions respectively. Here, we recall that $\int_{{\mathbb S}^1} dv = 1$. We then compute: 
\begin{eqnarray*} 
& & \hspace{-1cm}  I_\Phi := \int_{(v_1, \ldots , v_N) \in {\mathbb T}^n} Q_N \Phi (v_1, \ldots , v_N) \,  F_N (v_1, \ldots , v_N) \, dv_1 \ldots  dv_N  \\
& & \hspace{-0.3cm} = \frac{2}{N(N-1)} \sum_{i<j} \int_{(v_1, \ldots , v_N) \in {\mathbb T}^n}  \\
& & \hspace{0cm}
\left( H(\hat v_{ij}^* v_i) \, \int_{(w,z) \in {\mathbb T}^2}  \Phi (v_1, \ldots , w \hat v_{ij}, \ldots ,  z \hat v_{ij}, \ldots , v_N) \, g(w) \, g(z) \, dw \, dz \right. \\
& & \hspace{-1cm}  \left. \phantom{\int_{{\mathbb T}^2}  \Phi } +  (1-H(\hat v_{ij}^* v_i)) \, \Phi (v_1, \ldots , v_i , \ldots ,  v_j , \ldots , v_N)  \right) \,  F_N (v_1, \ldots , v_N) \, dv_1 \ldots  dv_N  
\end{eqnarray*}
The change of variables $(w,z) \in {\mathbb T}^2 \to (v'_i,v'_j) \in {\mathbb T}^2$ defined by  $v'_i = w \hat v_{ij}$, $v'_j = z \hat v_{ij}$ is of jacobian unity. Therefore, we deduce that 
\begin{eqnarray*} 
& & \hspace{-1cm}  I_\Phi = \frac{2}{N(N-1)} \sum_{i<j} \left\{ \int_{(v_1, \ldots , v_N) \in {\mathbb T}^n} \, \int_{(v'_i,v'_j) \in {\mathbb T}^2} H(\hat v_{ij}^{*} v_i) \, \Phi (v_1, \ldots , v'_i, \ldots ,  v'_j, \ldots , v_N) \right.  \\
& & \hspace{5cm}   g(\hat v^{*}_{ij} v'_i ) \, g(\hat v^{*}_{ij} v'_j) \,  F_N (v_1, \ldots , v_N) \, dv_1 \ldots  dv_N \, dv'_i \, dv'_j  \\
& & \hspace{0cm} + \left.  \int_{(v_1, \ldots , v_N) \in {\mathbb T}^n}  (1-H(\hat v_{ij}^{*} v_i)) \, \Phi (v_1, \ldots , v_N) \,  F_N (v_1, \ldots , v_N) \, dv_1 \ldots  dv_N  \right\} 
\end{eqnarray*}
Using (\ref{eq_master_eq}) leads to the weak form (\ref{eq_mas_BBDG_weak}). 

Now exchanging the notations $(v_i,v_j)$ and $(v'_i,v'_j)$ in the first term, we have:
\begin{eqnarray*} 
& & \hspace{-1cm}  I_\Phi = \frac{2}{N(N-1)} \sum_{i<j} \left\{ \int_{(v_1, \ldots , v_N) \in {\mathbb T}^n} \, \int_{(v'_i,v'_j) \in {\mathbb T}^2} H(\hat v^{'*}_{ij} v'_i) \, \Phi (v_1, \ldots , v_i, \ldots ,  v_j, \ldots , v_N)\right. \\
& & \hspace{1cm}   g(\hat v^{'*}_{ij} v_i ) \, g(\hat v^{'*}_{ij} v_j) \,  F_N (v_1, \ldots , v'_i, \ldots ,  v'_j, \ldots , v_N) \, dv_1 \ldots  dv_N \, dv'_i \, dv'_j  \\
& & \hspace{0cm} + \left.  \int_{(v_1, \ldots , v_N) \in {\mathbb T}^n}  (1-H(\hat v^{*}_{ij} v_i)) \, \Phi (v_1,  \ldots , v_N) \,  F_N (v_1, \ldots , v_N) \, dv_1 \ldots  dv_N  \right\} 
\end{eqnarray*}
with 
$$ \hat v^{'*}_{ij} = \frac{v'_i + v'_j}{|v'_i + v'_j|} . $$
With (\ref{eq_cond_exp}), we deduce that
\begin{eqnarray*} 
& & \hspace{-1cm}  Q_N^* F_N (v_1, \ldots , v_N) = \frac{2}{N(N-1)} \sum_{i<j} \\
& & \hspace{-1cm} \left\{  \int_{(v'_i,v'_j) \in {\mathbb T}^2} H(\hat v^{'*}_{ij} v'_i) \,   g(\hat v^{'*}_{ij} v_i ) \, g(\hat v^{'*}_{ij} v_j) \,  F_N (v_1, \ldots , v'_i, \ldots ,  v'_j, \ldots , v_N)\, dv'_i \, dv'_j  \right. \\
& & \hspace{6cm} \left. \phantom{\int_{{\mathbb T}^2}}  +  (1-H(\hat v_{ij}^* v_i))  \,  F_N (v_1, \ldots , v_N)  \right\} 
\end{eqnarray*}
Now, we change variables $(v'_i,v'_j) \in {\mathbb T}^2$ to $(u,z) \in {\mathbb S}^1 \times {\mathbb S}^1_+$  such that $u = \hat v^{'}_{ij}$, $z = \hat v^{'*}_{ij} v'_i$ where ${\mathbb S}^1_+ = \{ z \in {\mathbb S}^1 \, | \, \mbox{Re} z \geq 0 \}$. The restriction $z \in {\mathbb S}^1_+$ comes from the observation that $\mbox{Re} (\hat v^{'*}_{ij} v'_i)$ is non-negative. Indeed 
$$ \mbox{Re} (\hat v^{'*}_{ij} v'_i) = \frac{1}{|v'_i + v'_j|} \mbox{Re} ((v^{'*}_i + v^{'*}_j) v'_i) = \frac{1}{|v'_i + v'_j|} ( 1 + \mbox{Re} (v^{'*}_j v'_i)) ,$$
and since $|v'_i| = |v'_j|=1$, we have $|\mbox{Re} (v^{'*}_j v'_i)| \leq 1$ and  $1 + \mbox{Re} (v^{'*}_j v'_i) \geq0$. Reciprocally, it is easy to see that for $(u,v) \in {\mathbb S}^1 \times {\mathbb S}^1_+$ there exists a unique pair $(v'_i,v'_j) \in {\mathbb T}^2$ such that $u = \hat v^{'}_{ij}$, $z = \hat v^{'*}_{ij} v'_i$. Indeed, clearly, $v'_i = uz$, and we have $u z^* {v'_j}^* = 1$, which gives $v'_j = uz^*$. Finally, using the phases, one sees that $dv'_i \, dv'_j = 2  \, du \, dz$. Using this change of variables, we have: 
\begin{eqnarray*} 
& & \hspace{-1cm}  Q_N^* F_N (v_1, \ldots , v_N) = \frac{2}{N(N-1)} \sum_{i<j} \\
& & \hspace{-1cm} \left\{  2 \int_{(u,z) \in {\mathbb S}^1 \times {\mathbb S}^1_+} H(z) \,   g(u^* v_i ) \, g(u^* v_j) \,  F_N (v_1, \ldots , uz, \ldots ,  uz^*, \ldots , v_N)\, du \, dz  \right. \\
& & \hspace{6cm} \left. \phantom{\int_{{\mathbb T}^2}}  +  (1-H(\hat v_{ij}^* v_i))  \,  F_N (v_1, \ldots , v_N)  \right\} 
\end{eqnarray*}
With definition (\ref{eq_master_eq}) for the master equation, we are led to (\ref{eq_mas_BBDG_strong}), which ends the proof of Proposition \ref{prop_mas_eq_BBDG}. \endproof

\medskip
\noindent
{\bf Proof of Proposition \ref{prop_hierarchy_BBDG}.} We start with the weak form (\ref{eq_hierarchy_BBDG_weak}). We take $\Phi$ only depending on $(v_1,\ldots,v_k)$ in (\ref{eq_mas_BBDG_weak}). We get 
\begin{eqnarray} 
& & \hspace{-1cm}  
\frac{\partial}{\partial t} \int_{{\mathbb T}^N} F_N(v_1, \ldots , v_N,t) \, \Phi(v_1, \ldots , v_k) \, dv_1 \ldots dv_N = \nonumber\\
& & \hspace{0.5cm}  
= \frac{\partial}{\partial t} \int_{{\mathbb T}^k} F_{N,k} (v_1, \ldots , v_k,t) \, \Phi(v_1, \ldots , v_k) \, dv_1 \ldots dv_k \nonumber \\
& & \hspace{2cm}  
= {\mathcal T}_1 +  {\mathcal T}_2  +  {\mathcal T}_3,
\label{eq_T1T2T3}
\end{eqnarray}
where ${\mathcal T}_1$ (resp. ${\mathcal T}_2$ ; resp. ${\mathcal T}_3$) collects the terms at the right-hand side of  (\ref{eq_mas_BBDG_weak}) corresponding to $i < j \leq k$ (resp. $i \leq k <  j $ ; resp. $k < i <  j $). 

We first note that ${\mathcal T}_3=0$. Indeed, with $k < i <  j $, since $\Phi$ only depends on $(v_1, \ldots, v_k)$ and $g$ is a probability, the curly bracket at the right-hand side of (\ref{eq_mas_BBDG_weak}) becomes: 
\begin{eqnarray*} 
& & \hspace{0cm} \int_{(v'_i,v'_j) \in {\mathbb T}^2}  g(\hat v^{*}_{ij} v'_i ) \, g(\hat v^{*}_{ij} v'_j)  \, \Phi (v_1, \ldots , v_k)  \, dv'_i \, dv'_j - \, \Phi (v_1, \ldots , v_k) 
=0 .
\end{eqnarray*}
Concerning the first term ${\mathcal T}_1$, we can just integrate $(v_{k+1}, \ldots ,v_N)$ out of $F_N$ and the result is a just the first term at the right-hand side of (\ref{eq_hierarchy_BBDG_weak}). Now, we focus on ${\mathcal T}_2$. By permutation invariance, all terms of the sum over $j$ are identical. We can collect them into one single term for say $j=k+1$. 
There are $N-k$ such terms. We can then integrate $(v_{k+2}, \ldots ,v_N)$ out of $F_N$ leading to an expression involving only $F_{n,k+1}$. Now, for the same reason as above, we have :
\begin{eqnarray*} 
& & \hspace{-1cm} 
\int_{(v'_i,v'_j) \in {\mathbb T}^2}  g(\hat v^{*}_{i \, k+1} v'_i ) \, g(\hat v^{*}_{i \, k+1} v'_j)  \, \Phi (v_1, \ldots , v'_i, \ldots , v_k)  \, dv'_i \, dv'_{k+1}  = \\
& & \hspace{4cm} 
= \int_{v'_i\in {\mathbb S}^1}  g(\hat v^{*}_{i \, k+1} v'_i )  \, \Phi (v_1, \ldots , v'_i, \ldots , v_k)  \, dv'_i  .
\end{eqnarray*}
Collecting these observations leads to the second term of at the right-hand side of (\ref{eq_hierarchy_BBDG_weak}). 

Now, we turn to the strong form (\ref{eq_hierarchy_BBDG_strong}). The first term at the right-hand side of (\ref{eq_hierarchy_BBDG_strong}) is obtained from the corresponding term of (\ref{eq_hierarchy_BBDG_weak}) in exactly the same way as in the proof of proposition \ref{prop_mas_eq_BBDG}. We skip the details. We now focus on the second term. We exchange notations between $v_i$ and $v'_i$ and we change the notation $v_{k+1}$ into $v'_{k+1}$. We get
\begin{eqnarray*} 
& & \hspace{-0.5cm}
\int_{(v_i, v_{k+1}, v'_i) \in {\mathbb T}^3} H(\hat v_{i \, k+1}^{*} v_i) g(\hat v^{*}_{i \, k+1} v'_i ) \, \Phi (v_1, \ldots , v'_i, \ldots , v_k) \\
& & \hspace{5cm}
  F_{N,k+1} (v_1, \ldots , v_i, \ldots, v_k, v_{k+1}) \,  dv_i \,  dv_{k+1}  \, dv'_i  = \\
& & \hspace{-1cm}
= \int_{(v_i, v'_{k+1}, v'_i) \in {\mathbb T}^3} H(\hat v^{'*}_{i \, k+1} v'_i) g(\hat v^{'*}_{i \, k+1} v_i ) \, \Phi (v_1, \ldots , v_i, \ldots , v_k) \\
& & \hspace{5cm}
  F_{N,k+1} (v_1, \ldots , v'_i, \ldots, v_k, v'_{k+1})  \,  dv_i \,  dv'_{k+1}  \, dv'_i  ,
\end{eqnarray*}
and changing to the variables $( v'_i,  v'_{k+1}) \in {\mathbb T}^2$ to $(u,z) \in {\mathbb S}^1 \times {\mathbb S}^1_+$  such that $u = \hat v^{'}_{i \, k+1}$, $z = \hat v^{'*}_{i \, k+1} v'_i$, we get 
\begin{eqnarray*} 
& & \hspace{-1cm}
\int_{(v'_{k+1}, v'_i) \in {\mathbb T}^2} H(\hat v^{'*}_{i \, k+1} v'_i) g(\hat v^{'*}_{i \, k+1} v_i )  \,  F_{N,k+1} (v_1, \ldots , v'_i, \ldots, v_k, v'_{k+1})  \,  dv'_{k+1}  \, dv'_i  = \\
& & \hspace{2cm}
= 2 \int_{(u,z) \in {\mathbb S}^1 \times {\mathbb S}^1_+} H(z) g(u^* v_i ) \, 
  F_{N,k+1} (v_1, \ldots , uz, \ldots, v_k, uz^*)  \,  du  \, dz .
\end{eqnarray*}
Inserting this expression into the second term at the right-hand side of  (\ref{eq_hierarchy_BBDG_weak}) leads to the corresponding term of (\ref{eq_hierarchy_BBDG_strong}) and ends the proof. \endproof

\subsection{CL master equation and hierarchy}
\label{sec_appA_proof_master_CLD}

{\bf Proof of Proposition \ref{prop_mas_eq_CLD}.} The Markov transition operator for the CL dynamics is clearly
\begin{eqnarray*} 
& & \hspace{-1cm}  Q_N \Phi =  \frac{2}{N(N-1)} \sum_{i<j} \int_{z \in {\mathbb S}^1}  \frac{1}{2} \big(\Phi (v_1, \ldots , v_i, \ldots ,  z v_i, \ldots , v_N) + \\
& & \hspace{5cm} \Phi (v_1, \ldots , z v_j, \ldots ,  v_j, \ldots , v_N) \big)\, g(z) \, dz .
\end{eqnarray*}
We then compute: 
\begin{eqnarray*} 
& & \hspace{-1cm}  I_\Phi := \int_{(v_1, \ldots , v_N) \in {\mathbb T}^n} Q_N \Phi (v_1, \ldots , v_N) \,  F_N (v_1, \ldots , v_N) \, dv_1 \ldots  dv_N  \\
& & \hspace{-0.3cm} = \frac{2}{N(N-1)} \sum_{i<j} \int_{(v_1, \ldots , v_N) \in {\mathbb T}^n}  
\int_{z \in {\mathbb S}^1}  \frac{1}{2} \big(\Phi (v_1, \ldots , v_i, \ldots ,  z v_i, \ldots , v_N) + \\
& & \hspace{1cm} \Phi (v_1, \ldots , z v_j, \ldots ,  v_j, \ldots , v_N) \big)\, g(z)  \, F_N (v_1, \ldots , v_N)\, dz \, dv_1 \ldots  dv_N 
\end{eqnarray*}
This leads to expression (\ref{eq_mas_CLD_weak}).
 
By the change of variables $v'_j = zv_i$ in the first term and $v'_i = zv_j$ in the second one, we get 
\begin{eqnarray*} 
& & \hspace{-1cm}  I_\Phi := \frac{2}{N(N-1)} \sum_{i<j} \frac{1}{2} \left\{ \phantom{\int_{{\mathbb T}^{n}}} \right. \\
& & \hspace{0cm} \int_{(v_1, \ldots , \hat v_j, \ldots , v_N) \in {\mathbb T}^{n-1}}  \int_{v'_j \in {\mathbb S}^1} \Phi (v_1, \ldots , v_i, \ldots ,  v'_j, \ldots , v_N) \, g(v_i^* v'_j ) \\
& & \hspace{5cm}[ F_N ]_{\hat j} (v_1, \ldots , \hat v_j, \ldots , v_N) dv_1 \ldots d \hat v_j \ldots  dv_N \,dv'_j   \\
& & \hspace{0cm}  + \int_{(v_1, \ldots , \hat v_i, \ldots , v_N) \in {\mathbb T}^{n-1}}  \int_{v'_i \in {\mathbb S}^1} \Phi (v_1, \ldots , v'_i, \ldots ,  v_j, \ldots , v_N) \, g(v_j^* v'_i ) \\
& & \hspace{4.5cm} \left. \phantom{\int_{{\mathbb T}^{n}}} [ F_N ]_{\hat i} (v_1, \ldots , \hat v_i, \ldots , v_N) dv_1 \ldots d \hat v_i \ldots  dv_N \,dv'_i \right\} .
\end{eqnarray*}
By exchanging the roles of the primed and unprimed variables, we get: 
\begin{eqnarray*} 
& & \hspace{-1cm}  I_\Phi := \frac{2}{N(N-1)} \sum_{i<j} \int_{(v_1, \ldots , v_N) \in {\mathbb T}^{n}} \Phi (v_1, \ldots , v_N) \, \frac{1}{2} \left\{ g(v_i^* v_j ) [ F_N ]_{\hat j} (v_1, \ldots , \hat v_j, \ldots , v_N) +  \right. \\
& & \hspace{6cm}  \left. + g(v_j^* v_i ) [ F_N ]_{\hat i} (v_1, \ldots , \hat v_i, \ldots , v_N) \right\} dv_1 \ldots dv_N 
\end{eqnarray*}
Therefore, using the symmetry of $g$, we find:
\begin{eqnarray*} 
& & \hspace{-1cm}  Q_N^* F_N (v_1, \ldots , v_N) = \frac{1}{N(N-1)} \sum_{i<j}   g(v_i^* v_j ) \left\{ [ F_N ]_{\hat j} (v_1, \ldots , \hat v_j, \ldots , v_N) +  \right. \\
& & \hspace{8cm}  \left. +  [ F_N ]_{\hat i} (v_1, \ldots , \hat v_i, \ldots , v_N) \right\} 
\end{eqnarray*}
In view of (\ref{eq_master_eq}), we find (\ref{eq_mas_CLD_strong}), which ends the proof of Proposition \ref{prop_mas_eq_CLD}. \endproof

\medskip
\noindent
{\bf Proof of Proposition \ref{prop_hierarchy_CLD}.} We start with the weak form (\ref{eq_hierarchy_CLD_weak}). We take $\Phi$ only depending on $(v_1,\ldots,v_k)$ in (\ref{eq_mas_CLD_weak}). We get a similar expression as (\ref{eq_T1T2T3}) where ${\mathcal T}_1$ (resp. ${\mathcal T}_2$ ; resp. ${\mathcal T}_3$) collects the terms at the right-hand side of  (\ref{eq_mas_CLD_weak}) corresponding to $i < j \leq k$ (resp. $i \leq k <  j $ ; resp. $k < i <  j $). Again, ${\mathcal T}_3 = 0$. Indeed, for $k < i <  j $, since $\Phi$ only depends on $(v_1, \ldots, v_k)$ and $g$ is a probability, the curly bracket at the right-hand side of (\ref{eq_mas_CLD_weak}) becomes,  
\begin{eqnarray*} 
& & \hspace{0cm}  
\int_{z \in {\mathbb S}^1}  \Phi (v_1, \ldots , v_k) \, g(z) \, dz  - \Phi(v_1, \ldots, v_k) = 0. 
\end{eqnarray*}

We can integrate $(v_{k+1}, \ldots ,v_N)$ out of $F_N$ in ${\mathcal T}_1$ and get:
\begin{eqnarray*} 
& & \hspace{-0.5cm}  
{\mathcal T}_1 = \frac{2 \nu}{N-1} \int_{(v_1, \ldots , v_k) \in {\mathbb T}^k}  \left[ \phantom{\sum_{i<j \leq k}} \right.  \nonumber \\
& & \hspace{0cm} 
\sum_{i<j \leq k} \left\{ \int_{z \in {\mathbb S}^1}  \frac{1}{2} \left( \Phi (v_1, \ldots , v_i, \ldots ,  z v_i, \ldots , v_k) +  \Phi (v_1, \ldots , z v_j, \ldots ,  v_j, \ldots , v_k) \phantom{\frac{1}{2}} \hspace{-0.3cm} \right) \, g(z) \, dz \nonumber \right. \\
& & \hspace{5cm} \left. \left. \phantom{\int_{z \in {\mathbb S}^1}}  - \Phi(v_1, \ldots, v_k) \right\} \, F_{N,k} (v_1, \ldots , v_k) \right] \, dv_1 \ldots  dv_k  ,
\end{eqnarray*}
However, we can integrate once more the first term in the curly bracket because $\Phi (v_1, \ldots , v_i$, $\ldots ,  z v_i, \ldots , v_k)$ does not depend on $v_j$ and similarly with $i$ and $j$ exchanged in the second term. By permutation invariance, this can be expressed in terms of $F_{N,k-1}$. This leads to the first term at the right-hand side of (\ref{eq_hierarchy_CLD_weak}). 

By  permutation invariance, ${\mathcal T}_2$ is the sum of $N-k$ copies of the term obtained by choosing $i=k+1$. We can integrate $(v_{k+2}, \ldots ,v_N)$ out of $F_N$ and obtain an expression involving $F_{N,k+1}$ only. Now, since $i \leq k <  j $ and $\Phi$ only depends on $(v_1, \ldots, v_k)$, the curly bracket at the right-hand side of (\ref{eq_mas_CLD_weak}) becomes,
\begin{eqnarray} 
& & \hspace{0cm} \int_{v_{k+1} \in {\mathbb S}^1} \left\{ \int_{z \in {\mathbb S}^1}  \frac{1}{2} \left( \Phi (v_1, \ldots , v_i, \ldots , v_k) +  \Phi (v_1, \ldots , z v_{k+1}, \ldots , v_k) \phantom{\frac{1}{2}} \hspace{-0.3cm} \right)\, g(z) \, dz \right. \nonumber \\
& & \hspace{3cm} \left. \phantom{\int_{z \in {\mathbb S}^1}  \frac{1}{2}}  - \Phi(v_1, \ldots, v_k) \right\}  \, F_{N,k+1} (v_1, \ldots , v_k, v_{k+1}) \, dv_{k+1} \nonumber  \\
& & \hspace{0cm} = \int_{v_{k+1} \in {\mathbb S}^1} \frac{1}{2} \left( \int_{z \in {\mathbb S}^1} \Phi (v_1, \ldots , z v_{k+1}, \ldots , v_k) \, g(z) \, dz \right. \nonumber \\
& & \hspace{3cm} \left. \phantom{\int_{z \in {\mathbb S}^1}  \frac{1}{2}}  - \Phi(v_1, \ldots, v_k) \right)  \, F_{N,k+1} (v_1, \ldots , v_k, v_{k+1}) \, dv_{k+1} .
\label{eq_curly}
\end{eqnarray}
Collecting these observations, we can express ${\mathcal T}_2$ as follows: 
\begin{eqnarray*} 
& & \hspace{-0.5cm}  
{\mathcal T}_2 = \frac{2 \nu}{N-1} \int_{(v_1, \ldots , v_k) \in {\mathbb T}^k}  \left[ \phantom{\sum_{i<j \leq k}} \right.  \nonumber \\
& & \hspace{0cm} 
(N-k) \sum_{i\leq k} \int_{v_{k+1} \in {\mathbb S}^1} \frac{1}{2} \left\{ \int_{z \in {\mathbb S}^1} \Phi (v_1, \ldots , z v_{k+1}, \ldots , v_k) \, g(z) \, dz \nonumber \right. \\
& & \hspace{2.5cm} \left. \left. \phantom{\int_{z \in {\mathbb S}^1}}  - \Phi(v_1, \ldots, v_k) \right\} \, F_{N,k+1} (v_1, \ldots , v_{k+1}) \, dv_{k+1} \right] \, dv_1 \ldots  dv_k  ,
\end{eqnarray*}
However, we can integrate once more the terms in the curly bracket since $\Phi (v_1, \ldots , z v_{k+1},$ $\ldots , v_k)$ does not depend on $v_i$ and similarly $\Phi(v_1, \ldots, v_k)$ does not depend on $v_{k+1}$. The resulting expression involves only $F_{N,k}$ and leads to the second term at the right-hand side of (\ref{eq_hierarchy_CLD_weak}). 

To establish the strong form (\ref{eq_hierarchy_CLD_strong}) from the weak one (\ref{eq_hierarchy_CLD_weak}) is easy. We perform the change of variables $v_j = z v_i$ (resp. $v_i = z v_j$ ; resp. $v_i = z v_k$) in the first integral of the first sum (resp. in the second integral of the first sum ; resp. in the first integral of the second sum) and use the symmetry of $g$. \endproof

\setcounter{equation}{0}
\section{Limit $N \to \infty$ in the rescaled hierarchies}
\label{sec_appB_proof_grazing}

\subsection{Preliminaries}
\label{sec_appB_preliminaries}

In this section, we state the following straightforward lemma, which will be used throughout the proofs below.

\begin{lemma}
The following estimate holds in the sense of distributions:
$$ g_\varepsilon = \delta + \frac{\sigma^2 \varepsilon^2}{2} \delta'' + o( \varepsilon^2), \quad \quad h_\varepsilon = \delta + \frac{\tau^2 \varepsilon^2}{2} \delta'' + O( \varepsilon^3).  $$
More precisely, assuming that the test function $\Phi \in C^3({\mathbb R}/(2 \pi {\mathbb Z}))$, we have: 
$$ \left| \frac{1}{\varepsilon^2} \left( \int_0^{2 \pi} g_\varepsilon \, \Phi \, \frac{d \theta}{2 \pi} - (\Phi(0) + \frac{\sigma^2 \varepsilon^2}{2} \Phi''(0) ) \right) \right|  \leq C \varepsilon \| \Phi \|_{3,\infty} , $$
and similarly for $h_\varepsilon$ with $\sigma$ replaced by $\tau$ and the integral is taken on the interval $[-\pi/2,\pi/2]$. 
We denote by $C^3({\mathbb R}/(2 \pi {\mathbb Z}))$ the space of $2 \pi$ periodic, three-times continuously differentiable functions on ${\mathbb R}$ and by $\| \Phi \|_{3,\infty}$ the supremum of all the derivatives of $\Phi$ up to order~$3$. 
\label{lem_delta}
\end{lemma}

\noindent
We will also use the following operators: 
\begin{eqnarray} 
& & \hspace{-1cm} \Delta_{ij} = \frac{\partial^2}{\partial \theta_i^2} +  \frac{\partial^2}{\partial \theta_j^2} 
\label{eq_def_Delta_ij}
\\ 
& & \hspace{-1cm} D_{ij} = \frac{\partial^2}{\partial \theta_i^2} +  \frac{\partial^2}{\partial \theta_j^2} - 2 \frac{\partial^2}{\partial \theta_i \, \partial \theta_j} = \left( \frac{\partial}{\partial \theta_i} - \frac{\partial}{\partial \theta_j} \right)^2,
\label{eq_def_D_ij}
\end{eqnarray} 
in addition to $\bar D_{ij}$ defined at eq. (\ref{eq_def_bar_D_ij}).

\subsection{Limit $N \to \infty$ in the rescaled BDG hierarchy}
\label{sec_appB_proof_grazing_BBDG}

{\bf Proof of Theorem \ref{thm_hierarchy_grazing_BBDG}.} We first write the weak form of the hierarchy (\ref{eq_hierarchy_BBDG_weak}), introducing the phases (\ref{eq_phase_not}) and the rescaled noise distribution $g_\varepsilon$ and biasing function $h_\varepsilon$ (\ref{eq_rescale_h}), as well as the rescaled time (which will still be denoted by $t$ for the sake of simplicity). Now, we use the additive group notations of the phases for the multiplication of two elements of ${\mathbb S}^1$. We have 
\begin{eqnarray} 
& & \hspace{-0.7cm}  \frac{\partial}{\partial t} \int_{[0,2\pi]^k} F_{N,k}(\theta_1, \ldots , \theta_k,t) \, \Phi(\theta_1, \ldots , \theta_k) \, \frac{d\theta_1 \ldots d\theta_k}{(2 \pi)^k} = \frac{\nu }{(N-1) \varepsilon^2} \int_{(\theta_1, \ldots , \theta_k) \in [0,2\pi]^k} \left[ \phantom{\int_{[0,2\pi]^2}} \hspace{-0.5cm} \right. \nonumber \\
& & \hspace{-0.5cm}
 \sum_{i<j\leq k} h_\varepsilon(\theta_{\frac{i-j}{2}}) \left\{  \int_{(\theta'_i,\theta'_j) \in [0,2\pi]^2}  g_\varepsilon(\theta'_i - \hat \theta_{ij} ) \, g_\varepsilon(\theta'_j - \hat \theta_{ij} ) \, \Phi (\theta_1, \ldots , \theta'_i, \ldots ,  \theta'_j, \ldots , \theta_k) \, \frac{d\theta'_i \, d\theta'_j}{(2 \pi)^2}  \right. \nonumber  \\
& & \hspace{1cm} \left. \phantom{\int_{(\theta'_i,\theta'_j) \in [0,2\pi]^2}} - \Phi (\theta_1, \ldots , \theta_k) \, \right\}   F_{N,k} (\theta_1, \ldots , \theta_k) \nonumber \\
& & \hspace{-0.5cm}
+ (N-k) \int_{\theta_{k+1} \in [0,2\pi]} \sum_{i\leq k} h_\varepsilon(\theta_{\frac{i-(k+1)}{2}}) \left\{  \int_{\theta'_i \in [0,2\pi]}  g_\varepsilon(\theta'_i - \hat \theta_{i \, k+1}) \, \Phi (\theta_1, \ldots , \theta'_i, \ldots , \theta_k) \, \frac{d\theta'_i}{2 \pi}   \right. \nonumber  \\
& & \hspace{1cm} \left. \left. \phantom{\int_{\theta'_i \in [0,2\pi]}} - \Phi (\theta_1, \ldots , \theta_k) \, \right\}   F_{N,k+1} (\theta_1, \ldots , \theta_k, \theta_{k+1}) \, \frac{d\theta_{k+1}}{2 \pi} \right] \, \frac{d\theta_1 \ldots  d\theta_k}{(2 \pi)^k}
,
\label{eq_hierarchy_BBDG_weak_phases}
\end{eqnarray}
for any continuous test function $\Phi(\theta_1, \ldots , \theta_k)$ on $[0,2\pi]^k$. 
In this formula, $\theta_{\frac{i-j}{2}} = (\theta_i - \theta_j)/2$ or $(\theta_i - \theta_j)/2 + \pi$  with $-\pi/2 < \theta_{\frac{i-j}{2}} \leq \pi/2$ (modulo $2 \pi$) is the phase of the quantity $\hat v_{ij}^{*} v_i$, while $\hat \theta_{ij}$ is the phase of $\hat v_{ij}$, i.e. $\hat \theta_{ij} = (\theta_i + \theta_j)/2$ or $(\theta_i + \theta_j)/2 + \pi$ and $\theta_i - \hat \theta_{ij} = \hat \theta_{ij} - \theta_j = \theta_{\frac{i-j}{2}} \in [ - \pi/2, \pi/2]$ (modulo $2 \pi$). 

We first investigate the asymptotics of (\ref{eq_hierarchy_BBDG_weak_phases}) when $\varepsilon \to 0$ with fixed $N$. Take $\Phi \in C^3({\mathbb R}/(2 \pi {\mathbb Z}))^k$. We first prove the following formula: 
\begin{eqnarray} 
& & \hspace{-1cm}  \frac{\partial}{\partial t} \int_{[0,2\pi]^k} F_{N,k}(\theta_1, \ldots , \theta_k,t) \, \Phi(\theta_1, \ldots , \theta_k) \, \frac{d\theta_1 \ldots d\theta_k}{(2 \pi)^k} =\nonumber \\
& & \hspace{-1cm}
= \frac{\nu}{N-1} \left[ \, \, \left\{ \int_{(\theta_1, \ldots , \theta_k) \in [0,2\pi]^k} \sum_{i<j\leq k} \left\{ \sigma^2 (\Delta_{ij} \Phi) \, F_{N,k} + \,  \right. \right. \right. \nonumber \\ 
& & \hspace{0.5cm} 
\left. \phantom{\sum_{j}}  \left.  + \tau^2 \left[ \Phi\, (D_{ij} F_{N,k}) -  D_{ij} (F_{N,k} \Phi) \right] \right\} \delta(\theta_i - \theta_j) \, \frac{d\theta_1 \ldots  d\theta_k}{(2 \pi)^k} + o(\varepsilon^2) \right\}\nonumber \\
& & \hspace{-0.2cm} + (N-k) \left\{  \int_{(\theta_1, \ldots , \theta_{k+1}) \in [0,2\pi]^{k+1}} \sum_{i\leq k} \big\{ \sigma^2 (\Delta_{i \, k+1} \Phi)  \, F_{N,k+1} + \tau^2 \left[ \Phi\, (D_{i\, k+1} F_{N,k+1}) \right.  \right.  \nonumber \\
& & \hspace{1cm} 
\left. \left. \phantom{\sum_{j}} \left. - D_{i\, k+1} (F_{N,k+1} \Phi) \right] \big\} \, \delta(\theta_i - \theta_{k+1}) \, \frac{d\theta_1 \ldots  d\theta_{k+1}}{(2 \pi)^{k+1}} + o(\varepsilon^2) \right\} \, \,  \right],
\label{eq_hierarchy_BBDG_weak_asymptotics}
\end{eqnarray}
Moreover, thanks to the assumptions of the theorem, the $o(\varepsilon^2)$ terms are estimated by $C \varepsilon \| \Phi \|_{3,\infty}$ with $C$ independent of $N$ on any finite time interval $T$.  

\medskip
\noindent
{\bf Proof of (\ref{eq_hierarchy_BBDG_weak_asymptotics}).} We start by considering the first sum in (\ref{eq_hierarchy_BBDG_weak_phases}), where $i<j\leq k$. We use Lemma \ref{lem_delta} and get: \begin{eqnarray*} 
& & \hspace{-1cm}  
\int_{(\theta'_i,\theta'_j) \in ({\mathbb R}/(2 \pi {\mathbb Z}))^2} \Phi (\theta_1, \ldots , \theta'_i, \ldots ,  \theta'_j, \ldots , \theta_k) \, g_\varepsilon(\theta'_i - \hat \theta_{ij}) \, g_\varepsilon(\theta'_j - \hat \theta_{ij}) \,  \frac{d\theta'_i}{2 \pi} \,  \frac{d\theta'_j}{2 \pi} =  \\
& & \hspace{1.cm} 
= \left( \Phi + \frac{\sigma^2 \varepsilon^2}{2}  \left( \frac{\partial^2}{\partial \theta_i^2} +  \frac{\partial^2}{\partial \theta_j^2} \right) \Phi \right) (\theta_1, \ldots , \hat \theta_{ij}, \ldots ,  \hat \theta_{ij} , \ldots , \theta_k) + o(\varepsilon^2) .
\end{eqnarray*}
Now, using this formula, we compute the integral: 
\begin{eqnarray*} 
& & \hspace{-1cm} K_{ij} := \int_{(\theta_i,\theta_j) \in ({\mathbb R}/(2 \pi {\mathbb Z}))^2} h_\varepsilon( \theta_{\frac{i-j}{2}}) \\
& & \hspace{-1cm} \left\{ \int_{(\theta'_i,\theta'_j) \in ({\mathbb R}/(2 \pi {\mathbb Z}))^2} \Phi (\theta_1, \ldots , \theta'_i, \ldots ,  \theta'_j, \ldots , \theta_k) \, g_\varepsilon(\theta'_i - \hat \theta_{ij}) \, g_\varepsilon(\theta'_j - \hat \theta_{ij}) \,  \frac{d\theta'_i}{2 \pi} \,  \frac{d\theta'_j}{2 \pi} \right. \\
& & \hspace{3cm} \left. \phantom{\int_{(\theta'_i,\theta'_j) \in ({\mathbb R}/(2 \pi {\mathbb Z}))^2}} -  \Phi (\theta_1, \ldots , \theta_k) \right\} F_{N,k} (\theta_1, \ldots , \theta_k,t) \, \frac{d\theta_i}{2 \pi}  \frac{d\theta_j}{2 \pi} .   
\end{eqnarray*}
The change of variables $(\theta_i,\theta_j) \in  ({\mathbb R}/(2 \pi{\mathbb Z}))^2  \to (\hat \theta_{ij},\theta_{\frac{i-j}{2}}) \in {\mathbb R}/(2 \pi{\mathbb Z}) \times ]-\pi/2, \pi/2]$, is one-to-one onto and we have 
\begin{eqnarray*} 
& & \hspace{-1cm} K_{ij} = \int_{\hat \theta_{ij} \in {\mathbb R}/(2 \pi{\mathbb Z})}  \left( \Phi + \frac{\sigma^2 \varepsilon^2}{2}  \left( \frac{\partial^2}{\partial \theta_i^2} +  \frac{\partial^2}{\partial \theta_j^2} \right) \Phi \right) (\theta_1, \ldots , \hat \theta_{ij}, \ldots ,  \hat \theta_{ij} , \ldots , \theta_k) \\
& & \hspace{-0.5cm}  \int_{\theta_{\frac{i-j}{2}} \in [-\pi/2,\pi/2]}  h_\varepsilon( \theta_{\frac{i-j}{2}}) F_{N,k} (\theta_1, \ldots , \hat \theta_{ij} + \theta_{\frac{i-j}{2}} , \ldots, \hat \theta_{ij} - \theta_{\frac{i-j}{2}}, \ldots,  \theta_k,t) \, \frac{d\theta_{\frac{i-j}{2}}}{\pi} \, \frac{d \hat \theta_{ij}}{2 \pi} \\
& & \hspace{-0.5cm}  - \int_{\hat \theta_{ij} \in {\mathbb R}/(2 \pi{\mathbb Z})}  \int_{\theta_{\frac{i-j}{2}} \in [-\pi/2,\pi/2]} h_\varepsilon( \theta_{\frac{i-j}{2}}) \\
& & \hspace{2cm}
(\Phi F_{N,k}) (\theta_1, \ldots , \hat \theta_{ij} + \theta_{\frac{i-j}{2}} , \ldots, \hat \theta_{ij} - \theta_{\frac{i-j}{2}}, \ldots,  \theta_k) \, \frac{d\theta_{\frac{i-j}{2}}}{\pi} \, \frac{d \hat \theta_{ij}}{2 \pi}  + o(\varepsilon^2) .
\end{eqnarray*}
Using again Lemma \ref{lem_delta}, we find:
\begin{eqnarray*}
& & \hspace{-1cm} K_{ij} = \frac{\sigma^2 \varepsilon^2}{2} \int_{\hat \theta_{ij} \in {\mathbb R}/(2 \pi{\mathbb Z})}    \left[ \left( \left( \frac{\partial^2}{\partial \theta_i^2} +  \frac{\partial^2}{\partial \theta_j^2} \right) \Phi \right) \, F_{N,k} \right] (\theta_1, \ldots , \hat \theta_{ij} , \ldots, \hat \theta_{ij} , \ldots,  \theta_k,t) \, \frac{d \hat \theta_{ij}}{2 \pi} \\
& & \hspace{-1cm} + \frac{\tau^2 \varepsilon^2}{2} \int_{\hat \theta_{ij} \in {\mathbb R}/(2 \pi{\mathbb Z})}  \Phi (\theta_1, \ldots , \hat \theta_{ij}, \ldots ,  \hat \theta_{ij} , \ldots , \theta_k) \\
& & \hspace{4cm} \frac{\partial^2}{\partial \theta^2} [ F_{N,k} (\theta_1, \ldots , \hat \theta_{ij} + \theta , \ldots, \hat \theta_{ij} - \theta , \ldots,  \theta_k,t) ]_{\theta = 0} \, \frac{d \hat \theta_{ij}}{2 \pi} \\
& & \hspace{-1cm}  - \frac{\tau^2 \varepsilon^2}{2}  \int_{\hat \theta_{ij} \in {\mathbb R}/(2 \pi{\mathbb Z})} \frac{\partial^2}{\partial \theta^2}   [ (\Phi F_{N,k}) (\theta_1, \ldots , \hat \theta_{ij} + \theta , \ldots, \hat \theta_{ij} - \theta, \ldots,  \theta_k) ]_{\theta = 0} \frac{d \hat \theta_{ij}}{2 \pi}  \\
& & \hspace{-1cm} + o(\varepsilon^2) 
\end{eqnarray*}
Using definitions (\ref{eq_def_Delta_ij}) and (\ref{eq_def_D_ij}), we get:
\begin{eqnarray*} 
& & \hspace{-1cm}
K_{ij} =  \varepsilon^2 \int_{\hat \theta_{ij} \in {\mathbb R}/(2 \pi{\mathbb Z})} \left[ \frac{\sigma^2}{2} (\Delta_{ij} \Phi) F_{N,k} + \right. \\
& & \hspace{0cm} \left. + \frac{\tau^2}{2} [ \Phi D_{ij} F_{N,k} - D_{ij}(\Phi F_{N,k}) ] \right](\theta_1, \ldots , \hat \theta_{ij} , \ldots, \hat \theta_{ij}, \ldots,  \theta_k,t) \, \frac{d \hat \theta_{ij}}{2 \pi} + o(\varepsilon^2)
\end{eqnarray*} 
We can write
\begin{eqnarray*} 
& & \hspace{-1cm} K_{ij} = \varepsilon^2 \int_{\hat \theta_{ij} \in {\mathbb R}/(2 \pi{\mathbb Z})} \int_{\theta_{\frac{i-j}{2}} \in [-\pi/2,\pi/2]}
\left[ \frac{\sigma^2}{2} (\Delta_{ij} \Phi) F_{N,k} + \frac{\tau^2}{2} [ \Phi D_{ij} F_{N,k} \right. \\
& & \hspace{0.cm} \left. \phantom{\frac{\sigma^2}{2}} - D_{ij}(\Phi F_{N,k}) ] \right](\theta_1, \ldots , \theta_i , \ldots, \theta_j, \ldots,  \theta_k,t) \, \delta (\theta_{\frac{i-j}{2}}) \, \frac{d \hat \theta_{ij}}{2 \pi} \, \frac{d \theta_{\frac{i-j}{2}}}{\pi} + o(\varepsilon^2)
\end{eqnarray*} 
Going back to variables $(\theta_i,\theta_j)$, we find, since $\delta (\theta_{\frac{i-j}{2}}) = 2 \delta (\theta_i - \theta_j)$: 
\begin{eqnarray*} 
& & \hspace{-1cm} K_{ij} = 2 \varepsilon^2  \int_{(\theta_i, \theta_j) \in ({\mathbb R}/(2 \pi{\mathbb Z}))^2} 
\left[ \frac{\sigma^2}{2} (\Delta_{ij} \Phi) F_{N,k} + \right. \\
& & \hspace{3cm} \left. + \frac{\tau^2}{2} [ \Phi D_{ij} F_{N,k} - D_{ij}(\Phi F_{N,k}) ] \right]\, \delta (\theta_i - \theta_j) \, \frac{d \theta_i}{2 \pi} \, \frac{d \theta_j}{2 \pi} + o(\varepsilon^2) 
\end{eqnarray*} 
This gives the first term of (\ref{eq_hierarchy_BBDG_weak_asymptotics}). The second term is found using exactly the same type of computation. The details are omitted. This ends the proof of (\ref{eq_hierarchy_BBDG_weak_asymptotics}). 

\medskip
\noindent
{\bf End of proof of Theorem \ref{thm_hierarchy_grazing_BBDG}.} Now, using Green's formula, we get either the weak formulation of (\ref{eq_BBGKY_BBDG_grazing}): 
\begin{eqnarray} 
& & \hspace{-1cm}  \frac{\partial}{\partial t} \int_{[0,2\pi]^k} F_{N,k}(\theta_1, \ldots , \theta_k,t) \, \Phi(\theta_1, \ldots , \theta_k) \, \frac{d\theta_1 \ldots d\theta_k}{(2 \pi)^k} =\nonumber \\
& & \hspace{-1cm}
= \frac{\nu}{N-1} \left[ \, \, \left\{ \int_{(\theta_1, \ldots , \theta_k) \in [0,2\pi]^k} \sum_{i<j\leq k}  F_{N,k} \, \left\{ \sigma^2 (\Delta_{ij} \Phi)  \, \delta(\theta_i - \theta_j)  + \,  \right. \right. \right. \nonumber \\ 
& & \hspace{0.5cm} 
\left. \phantom{\sum_{j}}  \left.  + \tau^2 \left[ \, D_{ij}( \Phi \delta(\theta_i - \theta_j)) -   \Phi ( D_{ij} \delta(\theta_i - \theta_j))  \right] \right\}  \, \frac{d\theta_1 \ldots  d\theta_k}{(2 \pi)^k} + o(\varepsilon^2) \right\}\nonumber \\
& & \hspace{-0.7cm} + (N-k) \left\{  \int_{(\theta_1, \ldots , \theta_{k+1}) \in [0,2\pi]^{k+1}} \sum_{i\leq k} F_{N,k+1} \big\{ \sigma^2 ( \Delta_{i \, k+1} \Phi )  \, \delta(\theta_i - \theta_{k+1}) +   \right.  \nonumber \\
& & \hspace{-1.5cm} 
\left. \left. \phantom{\sum_{j}} + \tau^2 \left[  D_{i\, k+1} ( \Phi \delta(\theta_i - \theta_{k+1})) -  \Phi ( D_{i\, k+1} \delta(\theta_i - \theta_{k+1}) ) \right] \big\} \,  \, \frac{d\theta_1 \ldots  d\theta_{k+1}}{(2 \pi)^{k+1}} + o(\varepsilon^2) \right\} \, \,  \right],
\label{eq_hierarchy_BBDG_weak_asymptotics_2}
\end{eqnarray}
or the strong form: 
\begin{eqnarray} 
& & \hspace{-1cm}  \frac{\partial F_{N,k}}{\partial t} (\theta_1, \ldots , \theta_k,t) \,  =\nonumber \\
& & \hspace{-1cm}
= \frac{\nu}{N-1} \left[ \, \, \left\{  \sum_{i<j\leq k} \left\{ \sigma^2 \Delta_{ij} ( F_{N,k} \delta(\theta_i - \theta_j)) \,  + \,  \right. \right. \right. \nonumber \\ 
& & \hspace{0.5cm} 
\left. \phantom{\sum_{j}}  \left.  + \tau^2 \left[  (D_{ij} F_{N,k}) \, \delta(\theta_i - \theta_j) - F_{N,k} D_{ij} ( \delta(\theta_i - \theta_j)) \right] \right\}  \,  + o(\varepsilon^2) \right\}\nonumber \\
& & \hspace{-1cm} + (N-k) \left\{   \sum_{i\leq k} \int_{\theta_{k+1} \in [0,2\pi]} \big\{ \sigma^2 \Delta_{i \, k+1} ( F_{N,k+1} \delta(\theta_i - \theta_{k+1})) +  \right.  \nonumber \\
& & \hspace{-1.8cm} 
\left. \left. \phantom{\sum_{j}}  + \tau^2 \left[  (D_{i\, k+1} F_{N,k+1}) \delta(\theta_i - \theta_{k+1})  -  F_{N,k+1} D_{i\, k+1} (\delta(\theta_i - \theta_{k+1}) ) \right] \big\} \, \frac{d \theta_{k+1}}{2 \pi} \,  + o(\varepsilon^2) \right\}  \right].
\label{eq_hierarchy_BBDG_strong_asymptotics}
\end{eqnarray}

Now, letting $N \to \infty$ and $\varepsilon = \varepsilon_N \to 0$, we find that the leading order term in (\ref{eq_hierarchy_BBDG_weak_asymptotics_2}) or (\ref{eq_hierarchy_BBDG_strong_asymptotics}) is the second sum because of the factor $N -k$. Thanks to the hypotheses of the theorem, the other terms in the weak form (\ref{eq_hierarchy_BBDG_weak_asymptotics_2}) are either $o(\varepsilon_N^2)$ or $O(1/N)$ or $o(\varepsilon_N^2/N)$ multiplied by $\|\Phi\|_{3,\infty}$ and constants which are independent of $N$ and of time. Therefore, passing to the limit in the weak form or equivalently in the strong form (\ref{eq_hierarchy_BBDG_strong_asymptotics}) in the sense of distributions is allowed and leads to
\begin{eqnarray} 
& & \hspace{-0.5cm}  \frac{\partial F_{[k]}}{\partial t} (\theta_1, \ldots , \theta_k,t) \,  = \nu  \sum_{i\leq k} \int_{\theta_{k+1} \in [0,2\pi]} \big\{ \sigma^2 \Delta_{i \, k+1} ( F_{[k+1]} \delta(\theta_i - \theta_{k+1})) + \nonumber \\
& & \hspace{1cm} 
+ \tau^2 \left[  (D_{i\, k+1} F_{[k+1]}) \delta(\theta_i - \theta_{k+1})  -  F_{[k+1]} D_{i\, k+1} (\delta(\theta_i - \theta_{k+1}) ) \right] \big\} \frac{d \theta_{k+1}}{2 \pi}  .
\label{eq_hierarchy_BBDG_strong_asymptotics_3}
\end{eqnarray}
It remains to show that
\begin{eqnarray} 
& & \hspace{-1cm}  
\int_{\theta_{k+1} \in [0,2\pi]} \big\{ \sigma^2 \Delta_{i \, k+1} ( F_{[k+1]} \delta(\theta_i - \theta_{k+1})) + \nonumber \\
& & \hspace{-0.cm} 
+ \tau^2 \left[  (D_{i\, k+1} F_{[k+1]}) \delta(\theta_i - \theta_{k+1})  -  F_{[k+1]} D_{i\, k+1} (\delta(\theta_i - \theta_{k+1}) ) \right] \big\} \frac{d \theta_{k+1}}{2 \pi}  = \nonumber\\
& & \hspace{4.cm}  
= (\sigma^2 - \tau^2) (\bar D_{i \, k+1} F_{[k+1]} ) (\theta_1, \ldots, \theta_i , \ldots , \theta_k, \theta_i). 
\label{eq_curly_2}
\end{eqnarray}
We have:
\begin{eqnarray*} 
& & \hspace{-1cm}  
\Delta_{i \, k+1} ( F_{[k+1]} \delta(\theta_i - \theta_{k+1})) = (\Delta_{i \, k+1} F_{[k+1]}) \delta(\theta_i - \theta_{k+1}) + 2 F_{[k+1]} \delta''(\theta_i - \theta_{k+1}) + \\
& & \hspace{7cm}  
+ 2 \left( \frac{\partial F_{[k+1]}}{\partial \theta_i} - \frac{\partial F_{[k+1]}}{\partial \theta_{k+1}} \right) \delta'(\theta_i - \theta_{k+1}),   
\end{eqnarray*} 
where $\delta'$ and $\delta''$ denote the first and second derivatives of the Dirac delta. Using Green's formula, we deduce that 
\begin{eqnarray} 
& & \hspace{-1cm} 
\int_{\theta_{k+1} \in [0,2\pi]}  \Delta_{i \, k+1} ( F_{[k+1]} \delta(\theta_i - \theta_{k+1})) \frac{d \theta_{k+1}}{2 \pi} = \nonumber \\
& & \hspace{-0.5cm} 
= \left( ( \frac{\partial ^2}{\partial \theta_i^2} + \frac{\partial ^2}{\partial \theta_{k+1}^2} ) + 2 \frac{\partial ^2}{\partial \theta_{k+1}^2} + 2 \frac{\partial}{\partial \theta_{k+1}} ( \frac{\partial}{\partial \theta_i}  - \frac{\partial}{\partial \theta_{k+1}} )  \right) F_{[k+1]} (\theta_1, \ldots, \theta_i, \ldots, \theta_k, \theta_i) \nonumber \\
& & \hspace{-0.5cm} 
= (\bar D_{i \, k+1} F_{[k+1]} ) (\theta_1, \ldots, \theta_i , \ldots , \theta_k, \theta_i).
\label{eq_Delta_ik+1}
\end{eqnarray}
Now, we clearly have
\begin{eqnarray*} 
& & \hspace{-1cm}  
\int_{\theta_{k+1} \in [0,2\pi]}  (D_{i\, k+1} F_{[k+1]}) \delta(\theta_i - \theta_{k+1})  \frac{d \theta_{k+1}}{2 \pi}  = \nonumber\\
& & \hspace{3.cm}  
=  \left( \frac{\partial ^2}{\partial \theta_i^2} - 2 \frac{\partial ^2}{\partial \theta_i \, \partial \theta_{k+1}}+ \frac{\partial ^2}{\partial \theta_{k+1}^2} \right)  F_{[k+1]} (\theta_1, \ldots, \theta_i, \ldots, \theta_k, \theta_i). 
\end{eqnarray*}
On the other hand, since $  D_{i\, k+1} (\delta(\theta_i - \theta_{k+1})) = 4 \delta''(\theta_{k+1} - \theta_i) $, Green's formula leads to 
\begin{eqnarray*} 
& & \hspace{-1cm}  
\int_{\theta_{k+1} \in [0,2\pi]} F_{[k+1]} D_{i\, k+1} (\delta(\theta_i - \theta_{k+1}) ) \frac{d \theta_{k+1}}{2 \pi} = 4 \frac{\partial ^2 F_{[k+1]}}{\partial \theta_{k+1}^2}   (\theta_1, \ldots, \theta_i, \ldots, \theta_k, \theta_i). 
\end{eqnarray*}
But, by permutation invariance, we have 
\begin{eqnarray*} 
& & \hspace{-1cm}  
\frac{\partial ^2 F_{[k+1]}}{\partial \theta_{k+1}^2}  (\theta_1, \ldots, \theta_i, \ldots, \theta_k, \theta_i) = \frac{\partial ^2 F_{[k+1]}}{\partial \theta_i^2}  (\theta_1, \ldots, \theta_i, \ldots, \theta_k, \theta_i). 
\end{eqnarray*}
Therefore, 
\begin{eqnarray} 
& & \hspace{-1cm}  
\int_{\theta_{k+1} \in [0,2\pi]}  \big\{ (D_{i\, k+1} F_{[k+1]}) \delta(\theta_i - \theta_{k+1}) -   F_{[k+1]} D_{i\, k+1} (\delta(\theta_i - \theta_{k+1}) ) \big\} \frac{d \theta_{k+1}}{2 \pi}  = \nonumber\\
& & \hspace{6.cm}  
=  -  \bar D_{i\, k+1}  F_{[k+1]} (\theta_1, \ldots, \theta_i, \ldots, \theta_k, \theta_i). 
\label{eq_ik+1}
\end{eqnarray}
With (\ref{eq_Delta_ik+1}) and (\ref{eq_ik+1}), we easily verify (\ref{eq_curly_2}). Now, inserting (\ref{eq_curly_2}) into (\ref{eq_hierarchy_BBDG_strong_asymptotics_3}) leads to (\ref{eq_BBGKY_BBDG_grazing}) and ends the proof of Theorem \ref{thm_hierarchy_grazing_BBDG} \endproof

\subsection{Limit $N \to \infty$ in the rescaled CL hierarchy}
\label{sec_appB_proof_grazing_CLD}

{\bf Proof of Theorem \ref{thm_hierarchy_grazing_CLD}.} 
We first write the weak form of the hierarchy (\ref{eq_hierarchy_CLD_weak}), introducing the phases (\ref{eq_phase_not}) and the rescaled noise distribution $g_\varepsilon$, as well as the rescaled time (which will still be denoted by $t$ for the sake of simplicity). We have 
\begin{eqnarray} 
& & \hspace{-0.5cm}  
\frac{\partial}{\partial t} \int_{[0,2\pi]^N} F_{N,k}(\theta_1, \ldots , \theta_k,t) \, \Phi(\theta_1, \ldots , \theta_k) \, \frac{d\theta_1 \ldots d\theta_k}{(2 \pi)^k} = \nonumber \\
& & \hspace{0.5cm}  
= \frac{2 \nu}{(N-1) \varepsilon^2}  \left[ \sum_{i<j \leq k} \int_{(\theta_1, \ldots ,\hat \theta_i, \ldots, \hat \theta_j, \ldots, \theta_k) \in [0,2\pi]^{k-2}} \left\{
\phantom{\sum_{i<j \leq k}} \right. \right.   \nonumber \\
& & \hspace{-0.5cm} 
\frac{1}{2} \int_{(\theta_i,\varphi) \in [0,2\pi]^2} \Phi (\theta_1, \ldots , \theta_i, \ldots , \theta_i + \varphi, \ldots , \theta_k) F_{N,k-1} (\theta_1, \ldots , \theta_i, \ldots, \hat \theta_j, \ldots ,  \theta_k) \, g_\varepsilon(\varphi) \, \frac{d\varphi \, d\theta_i}{(2 \pi)^2}
\nonumber \\
& & \hspace{-0.5cm} 
+ \frac{1}{2} \int_{(\theta_j, \varphi) \in [0,2\pi]^2}  \Phi (\theta_1, \ldots , \theta_j + \varphi, \ldots ,  \theta_j, \ldots , \theta_k) F_{N,k-1} (\theta_1, \ldots, \hat \theta_i, \ldots , \theta_j, \ldots ,  \theta_k) \, g_\varepsilon(\varphi) \, \frac{d\varphi \, d\theta_j}{(2 \pi)^2} \nonumber \\
& & \hspace{1.5cm} \left. - \int_{(\theta_i,\theta_j) \in [0,2\pi]^2} \Phi(\theta_1, \ldots, \theta_k)  \, F_{N,k} (\theta_1, \ldots , \theta_k) \, \frac{d\theta_i \, d\theta_j}{(2 \pi)^2}\right\} \frac{ d \theta_1 \ldots d \hat \theta_i \ldots  d \hat \theta_j \ldots  d \theta_k}{(2 \pi)^{k-2}} \nonumber \\
& & \hspace{0.5cm} 
+ (N-k) \sum_{i\leq k} \int_{(\theta_1, \ldots ,\hat \theta_i, \ldots, \theta_k) \in [0,2\pi]^{k-1}} \frac{1}{2} \left\{ \phantom{\int_{(\theta_{k+1},\varphi) \in [0,2\pi]^2}} \right. \nonumber \\
& & \hspace{0cm} 
\int_{(\theta_{k+1},\varphi) \in [0,2\pi]^2} \Phi (\theta_1, \ldots ,  \theta_{k+1} + \varphi, \ldots , \theta_k) \, F_{N,k} (\theta_1, \ldots ,\hat \theta_i, \ldots, \theta_k, \theta_{k+1}) \, g_\varepsilon(\varphi) \, \frac{d\varphi \, d\theta_{k+1}}{(2 \pi)^2} \nonumber \\
& & \hspace{1.5cm} \left. \left.  - \int_{\theta_i \in [0,2\pi]}  \Phi(\theta_1, \ldots, \theta_k)  \, F_{N,k} (\theta_1, \ldots , \theta_k) \,  \right\} \frac{ d \theta_1 \ldots d \hat \theta_i \ldots  d \theta_k}{(2 \pi)^{k-1}}  \right] .
\label{eq_hierarchy_CLD_weak_phases}
\end{eqnarray}

Using Lemma \ref{lem_delta}, we have:
\begin{eqnarray*} 
& & \hspace{-1cm}   
\int_{\varphi \in [0,2 \pi]^2} \Phi (\theta_1, \ldots , \theta_i, \ldots ,  \theta_i + \varphi, \ldots , \theta_k)  \, g_\varepsilon(\varphi) \, \frac{d \varphi}{2 \pi} = \\
& & \hspace{-0cm}   
= \left( \Phi + \frac{\sigma^2 \varepsilon^2}{2} \frac{\partial^2
    \Phi}{\partial \theta_j^2} \right)  (\theta_1, \ldots , \theta_i,
\ldots ,  \theta_i , \ldots , \theta_k)  +  O(\varepsilon^3) \\
& & \hspace{-0cm}   
= \int_{\theta_j \in [0,2 \pi]} \left( \Phi + \frac{\sigma^2 \varepsilon^2}{2} \frac{\partial^2 \Phi}{\partial \theta_j^2} \right) (\theta_1, \ldots , \theta_i, \ldots ,  \theta_j , \ldots , \theta_k) \, \delta(\theta_j - \theta_i) \, \frac{d \theta_j}{2 \pi}
  +  O(\varepsilon^3) ,
\end{eqnarray*}
and similarly with $i$ and $j$ exchanged. In the same way, we have: 
\begin{eqnarray*} 
& & \hspace{-1cm}   
\int_{\varphi \in [0,2 \pi]^2} \Phi (\theta_1, \ldots , \theta_{k+1}+ \varphi, \ldots ,  \theta_k)  \, g_\varepsilon(\varphi) \, \frac{d \varphi}{2 \pi} = \\
& & \hspace{-0cm}   
= \int_{\theta_i \in [0,2 \pi]} \left( \Phi + \frac{\sigma^2 \varepsilon^2}{2} \frac{\partial^2 \Phi}{\partial \theta_i^2} \right) (\theta_1, \ldots , \theta_i, \ldots , \theta_k) \, \delta(\theta_i - \theta_{k+1}) \, \frac{d \theta_i}{2 \pi}
  + O(\varepsilon^3) .
\end{eqnarray*}
Moreover, the $O(\varepsilon^3)$ term is of the form $C \varepsilon^3 \| \Phi \|_{3,\infty}$ where $C$ is independent of $N$. Therefore, we have: 
\begin{eqnarray} 
& & \hspace{-0.5cm}  
\frac{\partial}{\partial t} \int_{[0,2\pi]^N} F_{N,k}(\theta_1, \ldots , \theta_k,t) \, \Phi(\theta_1, \ldots , \theta_k) \, \frac{d\theta_1 \ldots d\theta_k}{(2 \pi)^k} = \nonumber \\
& & \hspace{-0.5cm}  
= \frac{2 \nu}{(N-1) \varepsilon^2}  \left[ \sum_{i<j \leq k}  \left\{ \int_{(\theta_1, \ldots, \theta_k) \in [0,2\pi]^{k}} \left( \Phi (\theta_1, \ldots, \theta_k) \frac{1}{2} \left[ F_{N,k-1} (\theta_1, \ldots , \theta_i, \ldots, \hat \theta_j, \ldots ,  \theta_k) + \right. \right. \right.  \right. \nonumber \\
& & \hspace{4.5cm} 
\left. + F_{N,k-1} (\theta_1, \ldots, \hat \theta_i, \ldots , \theta_j, \ldots ,  \theta_k) \right]  \, \delta(\theta_j - \theta_i) \nonumber \\
& & \hspace{1.5cm}  
+  \frac{\varepsilon^2 \sigma^2}{4} \left[ \frac{\partial^2 \Phi}{\partial \theta_j^2} (\theta_1, \ldots, \theta_k) \, F_{N,k-1} (\theta_1, \ldots , \theta_i, \ldots, \hat \theta_j, \ldots ,  \theta_k) + \right.  \nonumber \\
& & \hspace{3cm} 
\left. + \frac{\partial^2 \Phi}{\partial \theta_i^2} (\theta_1, \ldots, \theta_k)  \, F_{N,k-1} (\theta_1, \ldots, \hat \theta_i, \ldots , \theta_j, \ldots ,  \theta_k) \right]  \, \delta(\theta_j - \theta_i) \nonumber \\
& & \hspace{4.5cm} \left. \left. \phantom{\frac{1}{2}}-  \Phi(\theta_1, \ldots, \theta_k)  \, F_{N,k} (\theta_1, \ldots , \theta_k)  \right) \frac{ d \theta_1 \ldots  d \theta_k}{(2 \pi)^k} + O(\varepsilon^3) \right\} \nonumber \\
& & \hspace{0.5cm} 
+ (N-k) \sum_{i\leq k} \left\{ \int_{(\theta_1, \ldots , \theta_k) \in [0,2\pi]^k} \frac{1}{2}   \left( \phantom{\int_{\theta_{k+1} \in [0,2\pi]}} \right. \right. \nonumber \\
& & \hspace{0cm} 
\int_{\theta_{k+1} \in [0,2\pi]} \left( \Phi + \frac{\varepsilon^2 \sigma^2}{2} \frac{\partial^2 \Phi}{\partial \theta_i^2} \right) (\theta_1, \ldots , \theta_k) \, F_{N,k} (\theta_1, \ldots ,\hat \theta_i, \ldots, \theta_k, \theta_{k+1}) \, \delta(\theta_i - \theta_{k+1})  \, \frac{d\theta_{k+1}}{2 \pi} \nonumber \\
& & \hspace{2cm} \left. \left. \left. \phantom{\int_{\theta_{k+1} \in [0,2\pi]}}  - \Phi(\theta_1, \ldots, \theta_k)  \, F_{N,k} (\theta_1, \ldots , \theta_k) \,  \right)   \frac{ d \theta_1 \ldots  d \theta_k}{(2 \pi)^k} + O(\varepsilon^3) \right\} \, \,   \right] .
\label{eq_hierarchy_CLD_weak_phases_2}
\end{eqnarray}
where, since $F_{N,k}$ is a probability, the $O(\varepsilon^3)$ terms are, as previously, of the form $C \varepsilon^3 \| \Phi \|_{3,\infty}$ and are uniform in $N$ and time. However, by permutation invariance, we have 
\begin{eqnarray*} 
& & \hspace{-0.5cm} 
\int_{(\theta_1, \ldots , \theta_k) \in [0,2\pi]^k}   \Big( \int_{\theta_{k+1} \in [0,2\pi]}  \Phi (\theta_1, \ldots , \theta_k) \, F_{N,k} (\theta_1, \ldots ,\hat \theta_i, \ldots, \theta_k, \theta_{k+1}) \, \delta(\theta_i - \theta_{k+1})  \, \frac{d\theta_{k+1}}{2 \pi}  \nonumber \\
& & \hspace{6cm}   - \Phi(\theta_1, \ldots, \theta_k)  \, F_{N,k} (\theta_1, \ldots , \theta_k) \,  \Big)   \frac{ d \theta_1 \ldots  d \theta_k}{(2 \pi)^k}  = \\
& & \hspace{-0.5cm} 
= \int_{(\theta_1, \ldots , \theta_k) \in [0,2\pi]^k}   \Big( \Phi (\theta_1, \ldots , \theta_k) \, F_{N,k} (\theta_1, \ldots ,\hat \theta_i, \ldots, \theta_k, \theta_i) \,  \nonumber \\
& & \hspace{6cm}   - \Phi(\theta_1, \ldots, \theta_k)  \, F_{N,k} (\theta_1, \ldots , \theta_k) \,  \Big)   \frac{ d \theta_1 \ldots  d \theta_k}{(2 \pi)^k}  = \\
& & \hspace{-0.5cm} 
= \int_{(\theta_1, \ldots , \theta_k) \in [0,2\pi]^k}   \Big( \Phi (\theta_1, \ldots , \theta_k) \, F_{N,k} (\theta_1, \ldots , \theta_k) \,  \nonumber \\
& & \hspace{6cm}  - \Phi(\theta_1, \ldots, \theta_k)  \, F_{N,k} (\theta_1, \ldots , \theta_k) \,  \Big)   \frac{ d \theta_1 \ldots  d \theta_k}{(2 \pi)^k} = 0 .
\end{eqnarray*}
Therefore, (\ref{eq_hierarchy_CLD_weak_phases_2}) yields
\begin{eqnarray*} 
& & \hspace{-0.5cm}  
\frac{\partial}{\partial t} \int_{[0,2\pi]^N} F_{N,k}(\theta_1, \ldots , \theta_k,t) \, \Phi(\theta_1, \ldots , \theta_k) \, \frac{d\theta_1 \ldots d\theta_k}{(2 \pi)^k} = \nonumber \\
& & \hspace{-0.5cm}  
= \frac{2 \nu}{(N-1) \varepsilon^2}  \left[ \sum_{i<j \leq k}  \left\{ \int_{(\theta_1, \ldots, \theta_k) \in [0,2\pi]^{k}} \left( \Phi (\theta_1, \ldots, \theta_k) \frac{1}{2} \left[ F_{N,k-1} (\theta_1, \ldots , \theta_i, \ldots, \hat \theta_j, \ldots ,  \theta_k) + \right. \right. \right.  \right. \nonumber \\
& & \hspace{4.5cm} 
\left. + F_{N,k-1} (\theta_1, \ldots, \hat \theta_i, \ldots , \theta_j, \ldots ,  \theta_k) \right]  \, \delta(\theta_j - \theta_i) \nonumber \\
& & \hspace{1.5cm}  
+  \frac{\varepsilon^2 \sigma^2}{4} \left[ \frac{\partial^2 \Phi}{\partial \theta_j^2} (\theta_1, \ldots, \theta_k) \, F_{N,k-1} (\theta_1, \ldots , \theta_i, \ldots, \hat \theta_j, \ldots ,  \theta_k) + \right.  \nonumber \\
& & \hspace{3cm} 
\left. + \frac{\partial^2 \Phi}{\partial \theta_i^2} (\theta_1, \ldots, \theta_k)  \, F_{N,k-1} (\theta_1, \ldots, \hat \theta_i, \ldots , \theta_j, \ldots ,  \theta_k) \right]  \, \delta(\theta_j - \theta_i) \nonumber \\
& & \hspace{4.5cm} \left. \left. \phantom{\frac{1}{2}}-  \Phi(\theta_1, \ldots, \theta_k)  \, F_{N,k} (\theta_1, \ldots , \theta_k)  \right) \frac{ d \theta_1 \ldots  d \theta_k}{(2 \pi)^k} + O(\varepsilon^3) \right\} \nonumber \\
& & \hspace{-0.5cm} \left.
+ \frac{(N-k) \varepsilon^2 \sigma^2}{4} \sum_{i\leq k} \left\{ \int_{(\theta_1, \ldots , \theta_k) \in [0,2\pi]^k} \frac{\partial^2 \Phi}{\partial \theta_i^2}  (\theta_1, \ldots , \theta_k) \, F_{N,k}  (\theta_1, \ldots , \theta_k) \,  \frac{ d \theta_1 \ldots  d \theta_k}{(2 \pi)^k} + O(\varepsilon) \right\} \right] . 
\end{eqnarray*}

Now, linking $\varepsilon$ and $N$ by (\ref{eq_scaling_CLD}), the $O(1)$ terms as $N \to \infty$ are on the first, second, fifth and sixth lines of the right-hand side. We can also drop the $O(\varepsilon)$ and $O(\varepsilon^3)$ remainders. All the terms which are dropped are estimated by $C \varepsilon \| \Phi \|_{3,\infty}$ uniformly in $N$ and time, thanks to the property that $F_{N,k}$ is a probability and to the assumptions of the theorem. Therefore, $F_{[k]}$ satisfies
\begin{eqnarray} 
& & \hspace{-0.5cm}  
\frac{\partial}{\partial t} \int_{[0,2\pi]^N} F_{[k]}(\theta_1, \ldots , \theta_k,t) \, \Phi(\theta_1, \ldots , \theta_k) \, \frac{d\theta_1 \ldots d\theta_k}{(2 \pi)^k} = \nonumber \\
& & \hspace{-0.5cm}  
= \nu  \left[ \sum_{i<j \leq k} \int_{(\theta_1, \ldots, \theta_k) \in [0,2\pi]^{k}} \Phi (\theta_1, \ldots, \theta_k) \left(  \left[ F_{[k-1]} (\theta_1, \ldots , \theta_i, \ldots, \hat \theta_j, \ldots ,  \theta_k) + \right. \right. \right.  \nonumber \\
& & \hspace{-0.5cm} 
\left. \left.+ F_{[k-1]} (\theta_1, \ldots, \hat \theta_i, \ldots , \theta_j, \ldots ,  \theta_k)  \right]  \, \delta(\theta_j - \theta_i) - 2 F_{[k]} (\theta_1, \ldots , \theta_k)  \right) \frac{ d \theta_1 \ldots  d \theta_k}{(2 \pi)^k} \nonumber \\
& & \hspace{0.5cm} \left.
+ \frac{\sigma^2}{2} \sum_{i\leq k} \int_{(\theta_1, \ldots , \theta_k) \in [0,2\pi]^k} \frac{\partial^2 \Phi}{\partial \theta_i^2}  (\theta_1, \ldots , \theta_k) \, F_{[k]}  (\theta_1, \ldots , \theta_k) \,  \frac{ d \theta_1 \ldots  d \theta_k}{(2 \pi)^k} \right] . 
\label{eq_hierarchy_CLD_weak_phases_3}
\end{eqnarray}
Now, using Green's formula, Eq. (\ref{eq_hierarchy_CLD_weak_phases_3}) appears as the weak form of (\ref{eq_BBGKY_CLD_grazing}), which ends the proof of Theorem \ref{thm_hierarchy_grazing_CLD}. \endproof


\end{document}